\documentclass[10pt,a4paper,oneside]{article}

\pdfoutput=1

\usepackage{caption}
\usepackage{mdframed}

\usepackage{graphicx}
\usepackage{pstricks}
\usepackage{diagbox} 
\usepackage[flushleft]{threeparttable}
\usepackage{makecell,booktabs}
\usepackage{multirow}
\usepackage{multicol}
\usepackage{xcolor,colortbl,hhline}
\definecolor{Gray}{gray}{0.80}
\definecolor{LightGray}{gray}{0.90}
\usepackage{dcolumn}
\usepackage{bm}
\usepackage{algorithmic,pstricks}
\usepackage{algorithm}
\usepackage{amsmath}
\usepackage{amsthm}
\usepackage{amssymb,latexsym,amsfonts}
\usepackage{tikz}
\usepackage{subfigure}
\usepackage{txfonts} 
\usepackage{comment}
\usepackage{bbm}
\usepackage{mathtools}
\usepackage{cite}
\usepackage{fancyhdr}

\usepackage[toc,page]{appendix} 

\usepackage{authblk}

\title{Noise-induced periodicity in a frustrated network of interacting diffusions}

\author[1]{Elisa Marini}
\author[2]{Luisa Andreis}
\author[3]{Francesca Collet}
\author[1,4]{Marco Formentin}
\affil[1]{{\small Dipartimento di Matematica ``Tullio Levi-Civita", Università degli Studi di Padova, Via Trieste 63, 35121 Padova, Italy}}
\affil[2]{{\small Dipartimento di Matematica, Politecnico di Milano, Via Edoardo Bonardi 9, 20133 Milano, Italy}}
\affil[3]{{\small Dipartimento di Informatica, Università degli Studi di Verona, Strada le Grazie 15, 37134 Verona, Italy}}
\affil[4]{{\small Padova Neuroscience Center, Via Giuseppe Orus 2, 35131 Padova, Italy}}

 \newtheorem{thm}{Theorem}[section]

 \theoremstyle{definition}
 
 \theoremstyle{remark}
 \newtheorem{remark}[thm]{Remark} 
 
 \numberwithin{equation}{section}

\allowdisplaybreaks

\let\svthefootnote\thefootnote

\begin{document}

%
%
%
%
%
%
%
%
%

\maketitle

\begin{abstract}
We investigate the emergence of a collective periodic behavior in a frustrated network of interacting diffusions. Particles are divided into two communities depending on their mutual couplings. On the one hand, both intra-population interactions are positive; each particle wants to conform to the average position of the particles in its own community. On the other hand, inter-population interactions have different signs: the particles of one population want to conform to the average position of the particles of the other community, while the particles in the latter want to do the opposite. We show that this system features the phenomenon of noise-induced periodicity: in the infinite volume limit, in a certain range of interaction strengths, although the system has no periodic behavior in the zero-noise limit, a moderate amount of noise may generate an attractive periodic law.

\let\thefootnote\relax
\addtocounter{footnote}{-1}\let\thefootnote\svthefootnote 

\vspace{0.3cm}


\end{abstract}






	

\section{Introduction}
Robust periodic behaviors are frequently encountered in life sciences and are indeed one of the most commonly observed self-organized dynamics. For instance,  spontaneous brain activity exhibits rhythmic oscillations called alpha and beta waves \cite{ermentrout2010mathematical}. From a theoretical standpoint, the mechanism driving the emergence of periodic behaviors in such systems is poorly understood. For example, neurons neither have any tendency to behave periodically on their own, nor are subject to any periodic forcing; nevertheless, they organize to produce a regular motion perceived at the macroscopic scale \cite{wallace2011emergent}. Various models of large families of interacting particles showing self-sustained oscillations have been proposed; we refer the reader to \cite{DiLo,AlMi2019,AlMi2021a,AnTo2018,CDFT2021,CoFo2019,GiPo15,DaPGiRe14, DFP2021,dfr,LuPo2020,LuPo2021,FeFoNe09}, where possible mechanisms leading to a rhythmic behavior are discussed and many related references are given.\\
Here we mention two mechanisms - which are of interest to us - capable to induce or enhance periodic behaviors in stochastic systems with many degrees of freedom. The first one is noise. The role of the noise is twofold: on the one hand, it can lead to oscillatory laws in systems of nonlinear diffusions whose deterministic counterparts do not display any periodic behavior \cite{Sc1985a,Sc1985b}; on the other hand, it can facilitate the transition from an equilibrium solution to macroscopic self-organized oscillations \cite{CoDaPFo15,LG-ONS-G04,Tou2012}.\\
The second mechanism is the topology of the interaction network. It has been recently pointed out in \cite{DiLo,Tou2019,FMD2016} that a specific network structure may favor the emergence of collective rhythms. In particular, in \cite{Tou2019,FMD2016}, the large volume dynamics of a two-population generalization of the mean field Ising model is considered.
The system is shown to undergo a transition from a disordered phase, where the magnetizations of both populations fluctuate around zero, to a phase in which they both display a macroscopic regular rhythm. Such a transition is driven by inter- and intra-population interactions of different strengths and signs leading to dynamical frustration.\\
In the present paper we combine the two mechanisms described above and we design a toy model of frustratedly interacting diffusions that shows noise-induced periodicity, in the sense that periodic oscillations appear for an intermediate amount of noise. The peculiar feature of the model under consideration is that the structure of the interaction network depends on the noise in that it is the noise that switches on the interaction terms, thus leading to periodic dynamics.

\section{Description of the model and outline of the results}\label{sec:model:results}

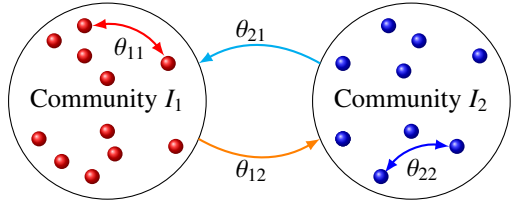
\begin{figure}[h!]\centering
\begin{tikzpicture}
\draw (1.,1.5) circle (1.3cm);
\draw (5.,1.5) circle (1.3cm);
\node (3) at (1.,1.5) {$\mbox{Community $I_1$}$};
\node (4) at (5.,1.5) {$\mbox{Community $I_2$}$};
\draw [-latex,orange,thick, bend right] (2.2,1.) to (3.8,1.);
\draw [-latex,cyan,thick, bend right] (3.8,2) to (2.2,2);
\draw [latex-latex,blue,thick, bend right]  (5.5,.9) to (4.65,0.55);
\draw [latex-latex,red,thick, bend left]  (.8,2.5) to (1.75,2.1);
\draw (5.15,0.6) node {$\theta_{22}$};
\draw (1.3,2.2) node {$\theta_{11}$};
\draw (2.9,0.58) node {$\theta_{12}$};
\draw (2.9,2.44) node {$\theta_{21}$};
\shade[ball color=blue]  (4.1,1.) circle (.1cm);
\shade[ball color=blue]  (4.6,0.5) circle (.1cm);
\shade[ball color=blue]  (4.1,1.) circle (.1cm);
\shade[ball color=blue]  (5.6,0.9) circle (.1cm);
\shade[ball color=blue]  (5.,1.1) circle (.1cm);
\shade[ball color=blue]  (4.1,1.) circle (.1cm);
\shade[ball color=blue]  (4.7,2.5) circle (.1cm);
\shade[ball color=blue]  (4.9,1.9) circle (.1cm);
\shade[ball color=blue]  (5.1,2.3) circle (.1cm);
\shade[ball color=blue]  (5.9,2.1) circle (.1cm);
\shade[ball color=blue]  (4.1,2.) circle (.1cm);
\shade[ball color=red]  (.4,0.7) circle (.1cm);
\shade[ball color=red]  (.8,0.5) circle (.1cm);
\shade[ball color=red]  (.1,1.) circle (.1cm);
\shade[ball color=red]  (1.9,0.9) circle (.1cm);
\shade[ball color=red]  (1.,1.1) circle (.1cm);
\shade[ball color=red]  (1.1,.8) circle (.1cm);
\shade[ball color=red]  (.7,2.5) circle (.1cm);
\shade[ball color=red]  (1.,1.8) circle (.1cm);
\shade[ball color=red]  (.3,2.3) circle (.1cm);
\shade[ball color=red]  (.7,2.1) circle (.1cm);
\shade[ball color=red]  (1.8,2.) circle (.1cm);
\end{tikzpicture}
\caption{\footnotesize{A schematic representation of the interaction network. Particles are divided into two communities, $I_1$ and $I_2$. Ignoring inter-population interactions, each community taken alone is a mean field system with interaction strength $\theta_{i i}$ ($i=1$ or $2$). When we couple the two communities, population $I_1$ (resp. $I_2$) influences the dynamics of population $I_2$ (resp. $I_1$) through the average position of its particles with strength $\theta_{21}$ (resp. $\theta_{12}$).}}
\label{CW}
\end{figure}

\noindent  Let us consider a system of $N$ diffusive particles on $\mathbb{R}$. We divide the $N$ particles into two disjoint communities of sizes $N_1$ and $N_2$ respectively and we denote by $I_1$ (resp. $I_2$) the set of sites belonging to the first (resp. second) community. 
In this setting, we indicate with $\left(x^{(N)}_{j}(t)\right)_{j=1, \dots, N_{1}}$ the positions at time $t$ of the particles of population $I_1$ and with $\left(y^{(N)}_{j}(t)\right)_{j=1, \dots, N_{2}}$ the positions at time $t$ of the particles of  population $I_2$, so that 
\[ 
\mathbf{z}^{(N)}(t) =  \Big( \, \overbracket{x^{(N)}_1(t), x^{(N)}_2(t), \dots, x^{(N)}_{N_1}(t)}^{\text{\scriptsize Community $I_1$}}, \, \overbracket{y^{(N)}_{1}(t), y^{(N)}_{2}(t), \dots, y^{(N)}_{N_2}(t)}^{\text{\scriptsize Community $I_2$}} \, \Big)
\] 
represents the state of the whole system at time $t$. The basic feature of our model is that the strength of the interaction between particles depends on the community they belong to: 
$\theta_{11}$ and $\theta_{22}$ tune the interaction between sites of the same community, whereas $\theta_{12}$ and $\theta_{21}$ control the coupling strength between particles of different groups.  In fact, we construct the network of interacting diffusions visualized in Fig.~\ref{CW}. 
A crucial feature for the system to show periodic behavior is \emph{frustration of the network}, i.e. the  inter-community interactions must have opposite signs.

Now we introduce the microscopic dynamics we are interested in. Let
\[
m^{(N)}_{1} (t) := \frac{1}{N_1} \sum_{j =1}^{N_{1}} x^{(N)}_j(t) \quad \text{ and } \quad m^{(N)}_{2} (t) := \frac{1}{N_2} \sum_{j =1}^{N_{2}} y^{(N)}_j(t)
\]
be the empirical means of the positions of the particles in populations $I_1$ and $I_2$, respectively, at time~$t$. Moreover, denote by $\alpha := \frac{N_1}{N}$ the fraction of sites belonging to the first group. Then, omitting time dependence for notational convenience, the interacting particle system we are going to study reads 

\begin{align}\label{micro:dyn}
dx^{(N)}_j &=\left(-\left(x^{(N)}_j\right)^3+x^{(N)}_j\right)dt-\alpha\:\theta_{11}\left(x^{(N)}_j-m^{(N)}_{1}\right)dt\nonumber\\
&-\left(1-\alpha\right)\theta_{12}\left(x^{(N)}_j-m^{(N)}_{2}\right)dt+\sigma dw_j, & \mbox{ for } j = 1, \dots, N_{1},\nonumber\\
&&&\nonumber\\
dy^{(N)}_j &=\left(-\left(y^{(N)}_j\right)^3+y^{(N)}_j\right)dt-\alpha\:\theta_{21}\left(y^{(N)}_j-m^{(N)}_{1}\right)dt\nonumber\\
&-\left(1-\alpha\right)\theta_{22}\left(y^{(N)}_j-m^{(N)}_{2}\right)dt+\sigma dw_{N_{1}+j}, &\mbox{ for } j = 1, \dots, N_{2},
\end{align}
where $\left(w_j(t); t \geq 0\right)_{j=1,\dots,N}$ are $N$ independent copies of a standard Brownian motion. Here \mbox{$\sigma \geq 0$} is the parameter that tunes the amount of noise in the system, since the diffusion coefficient is the same for each coordinate. 

\begin{remark} 
Existence and uniqueness of a strong solution to \eqref{micro:dyn} can be established via the Khasminskii criterion \cite{Kha12,MeTw93}: by taking the norm-like function
\[
V\left(\mathbf{z}^{(N)}\right) = \frac{1}{N_1} \sum_{i=1}^{N_1} \left[ \frac{\left(x_i^{(N)}\right)^4}{4} + \frac{\left(x_i^{(N)}\right)^2}{2} \right] + \frac{1}{N_2} \sum_{i=1}^{N_2} \left[ \frac{\left(y_i^{(N)}\right)^4}{4} + \frac{\left(y_i^{(N)}\right)^2}{2} \right],
\]
one obtains an inequality of the form $\mathcal{L}V\left(\mathbf{z}^{(N)}\right) \leq k \left[1 + V\left(\mathbf{z}^{(N)}\right)\right]$, for some $k>0$, with $\mathcal{L}$ the infinitesimal generator of diffusion \eqref{micro:dyn}.
\end{remark}

\noindent Notice that in system \eqref{micro:dyn} the two groups of particles interact only through their empirical means. This makes our model mean field and, in particular, when $\theta_{11}=\theta_{22}=\theta_{12}=\theta_{21}=\theta > 0$, the system of equations \eqref{micro:dyn} reduces to the mean field interacting diffusions considered in \cite{Daw83}. In a general setting, all the interaction parameters can be either positive or negative allowing both cooperative/conformist and uncooperative/anti-conformist interactions. In the present paper, we focus on the case $\theta_{11}>0$, $\theta_{22}>0$ and  $\theta_{12} \theta_{21}<0$. Moreover, without loss of generality, we make the specific choice $\theta_{12}>0$ and $\theta_{21}<0$, which means that particles in $I_1$ tend to conform to the average particle position of community $I_2$, whereas particles in $I_2$ prefer to differ from the average particle position of community $I_1$ (see  Eq.~\eqref{micro:dyn}).

Numerical simulations of system \eqref{micro:dyn} with large $N$ show that $m^{(N)}_1(t)$ and $m^{(N)}_2(t)$  display an oscillatory behavior in appropriate regions of the parameter space (see Section~\ref{sec:numerical:simulations}). This led us to investigate the thermodynamic limit of our system of interacting diffusions, that is, the limit when the number of particles goes to infinity. It is known that solutions of SDEs like~\eqref{micro:dyn} cannot have a time-periodic law, as these solutions are either positive recurrent, null recurrent or transient; see \cite{Sc1985a, Sc1985b} and references therein. However, the mean field interaction in~\eqref{micro:dyn} has a peculiar feature. When the interaction is of this type, at any time $t$, the empirical average of the particle positions in~\eqref{micro:dyn} is expected to converge, as the number of particles goes to infinity, to a limit given by the solution of a \emph{nonlinear SDE}. Nonlinear SDEs are SDEs where the coefficients  depend on the law of the solution itself and, in contrast with systems like~\eqref{micro:dyn}, such nonlinear SDEs might have solutions with periodic law, see~\cite{Sc1985b}.  Therefore, the oscillations in the trajectories of $m_1^N(t)$ and $m_2^N(t)$ shown by simulations can be theoretically explained via the thermodynamic limit of the system.

We outline here the main results presented in the sequel. We follow an approach similar to the one adopted in \cite{CoDaPFo15}.
\begin{enumerate}
\item 
In Section~\ref{sec:prop:chaos} we prove that, starting from i.i.d. initial conditions, independence propagates in time when taking the infinite volume limit. In particular, as $N$ grows large, the time evolution of a pair of representative particles, one for each population, is described by the limiting dynamics
\begin{align}\label{prop:chaos}
dx &= \left[-x^3+x-\alpha\theta_{11}\left(x-\mathbb{E}[x]\right)-(1-\alpha)\theta _{12}\left(x-\mathbb{E}[y]\right)\right]dt+\sigma dw_1 \nonumber\\[.1cm]
dy &= \left[-y^3+y-\alpha\theta_{21}\left(y-\mathbb{E}[x]\right)-(1-\alpha)\theta_{22}\left(y-\mathbb{E}[y]\right)\right]dt+\sigma dw_2,
\end{align}
where notation $\mathbb{E}$ stands for the expectation with respect to the probability measure \mbox{$Q(t;\cdot)=\text{Law}(x(t),y(t))$}, for every $t \in [0,T]$, and $(w_1(t); 0 \leq t \leq T)$ and $(w_2(t); 0 \leq t \leq T)$ are two independent standard Brownian motions. \\
In particular, we show that, for all $T > 0$ and for all $t\in [0,T]$, any random vector of the form $\left(x^{(N)}_{i_1}(t),\dots,x^{(N)}_{i_{k_1}}(t), y^{(N)}_{j_1}(t),\dots, y^{(N)}_{j_{k_2}}(t)\right)$ converges in distribution, as $N$ goes to infinity, to a vector $(x_1(t), \dots, x_{k_1}(t),$ $y_1(t), \dots, y_{k_2}(t))$, whose entries are independent random variables such that $x_i(t)$ ($i=1, \dots, k_1$) are copies of the solution to the first equation in \eqref{prop:chaos} and $y_i(t)$ ($i=1, \dots, k_2$) are copies of the solution to the second equation in \eqref{prop:chaos}. This is usually referred to as the phenomenon of  propagation of chaos. See~\cite{Luisa2018} for a proof in a general framework of weakly interacting diffusions with jumps. Notice that our model is a two-population version of the model in Section 4 in \cite{Luisa2018}, as here there are no jumps and the drift term in \eqref{prop:chaos} satisfies their Assumption~3.
\item Being nonlinear, system \eqref{prop:chaos} is a good candidate for having a solution with periodic law. It is however very hard to gain insight into its long-time behavior or to find periodic solutions as the problem is infinite dimensional, due to the presence of nonlinearity and noise. As a first step, in Section~\ref{sec:zero:noise:dyn} we study the limiting system \eqref{prop:chaos} in the absence of noise and, in particular, we argue that oscillatory behaviors are not observed when $\sigma=0$. This remains true for small values of $\sigma>0$ in some parameter regimes. See Section \ref{sec:numerical:simulations} for details.
\item In Sections~\ref{sec:FP} and \ref{sec:gauss:approx} we tackle system \eqref{prop:chaos} with noise. We show that, in the presence of an appropriate amount of noise, the limiting positions of representative particles of the two populations evolve approximately as a pair of independent Gaussian processes (\emph{small-noise Gaussian approximation}). This reduces the problem to a finite dimensional one, since we provide the explicit (deterministic) equations for the mean and variance of those processes. The dynamical system describing the time evolution of the means and the variances has a Hopf bifurcation and, as a consequence, in a certain range of the noise intensity, it has a limit cycle as a long-time attractor, implying that the laws of the previously mentioned Gaussian processes are periodic. Thus, the small-noise Gaussian approximation gives a good qualitative description of the emergence of the self-sustained oscillations observed for system \eqref{micro:dyn} (see Section~\ref{sec:numerical:simulations}). 
\end{enumerate}

Intuitively, the mechanism behind the emergence of periodicity in our system is similar to the one in \cite{FMD2016} and can be described as follows. Imagine to start with two independent communities, that is, particles evolve according to system \eqref{micro:dyn} with $\theta_{12}=\theta_{21}=0$. When the intra-population interaction strengths $\theta_{11}$ and $\theta_{22}$ are large enough, each population tends to its own rest state, that one may guess to be (close to) one of the minima of the double well potential $V(x)=\frac{x^4}{4}-\frac{x^2}{2}$ (see \cite{Daw83}). The key aspect, which we believe makes the model under consideration interesting, is that linking the two populations together within an interaction network with $\theta_{12}\theta_{21}<0$ is not enough for periodic behaviors to appear. Dynamical frustration and, in turn, oscillations arise only when the noise intensity is large enough, as the interaction terms in  system \eqref{micro:dyn} are switched on by the noise. Indeed, when $\sigma=0$ and all the particles in a same population share the same initial condition, the system is attracted to an equilibrium point where $x_j^{(N)}=y_k^{(N)}=m_1^{(N)}=m_2^{(N)}$ ($j=1, \dots, N_1$; $k=1, \dots, N_2$) - see Fig.~\ref{fig:small_vs_big_noise} - and, thus, the interaction terms vanish. It follows that  
the zero-noise dynamics does not display any periodic behavior. On the contrary, if $\sigma$ is positive and sufficiently large, particles do not get stuck at equilibrium points, as diffusion is enhanced, and the interaction terms start playing a role, generating dynamical frustration. The two populations form now a frustrated pair of systems where the rest state of the first is not compatible with the rest position of the second. As a consequence, the  dynamics does not settle down to a fixed equilibrium and keeps oscillating. Therefore, the noise is responsible for the emergence of a stable rhythm (see Section~\ref{sec:M2}). This feature is the hallmark of the phenomenon of noise-induced periodicity.

\section{Noise-induced periodicity: numerical study}\label{sec:numerical:simulations}

\begin{figure}[h!p]
    \captionsetup[subfigure]{labelformat=empty}
    \centering
    \begin{center}
        \subfigure{\includegraphics[width=0.325\textwidth]{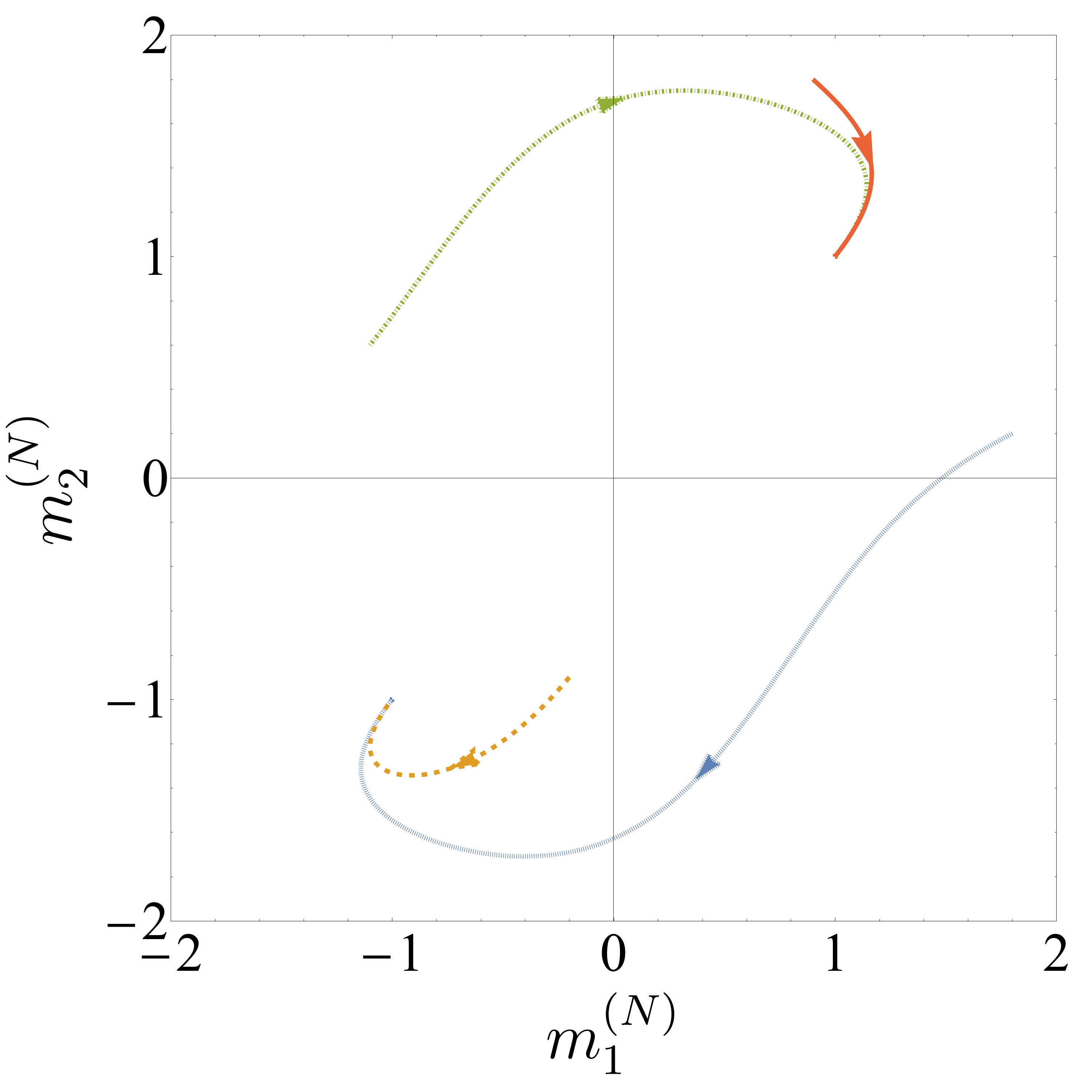}}
        \subfigure{\includegraphics[width=0.325\textwidth]{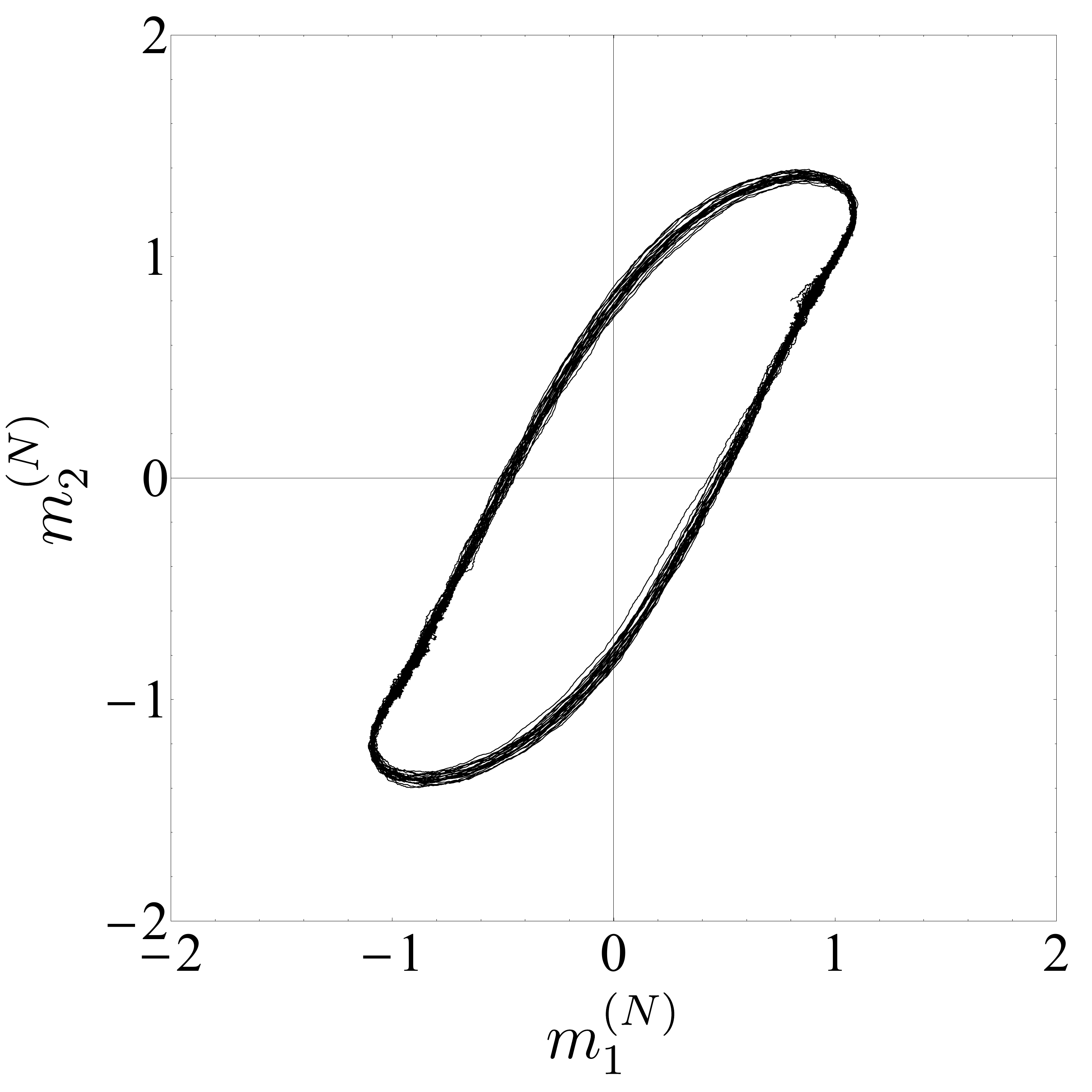}}
        \subfigure{\includegraphics[width=0.325\textwidth]{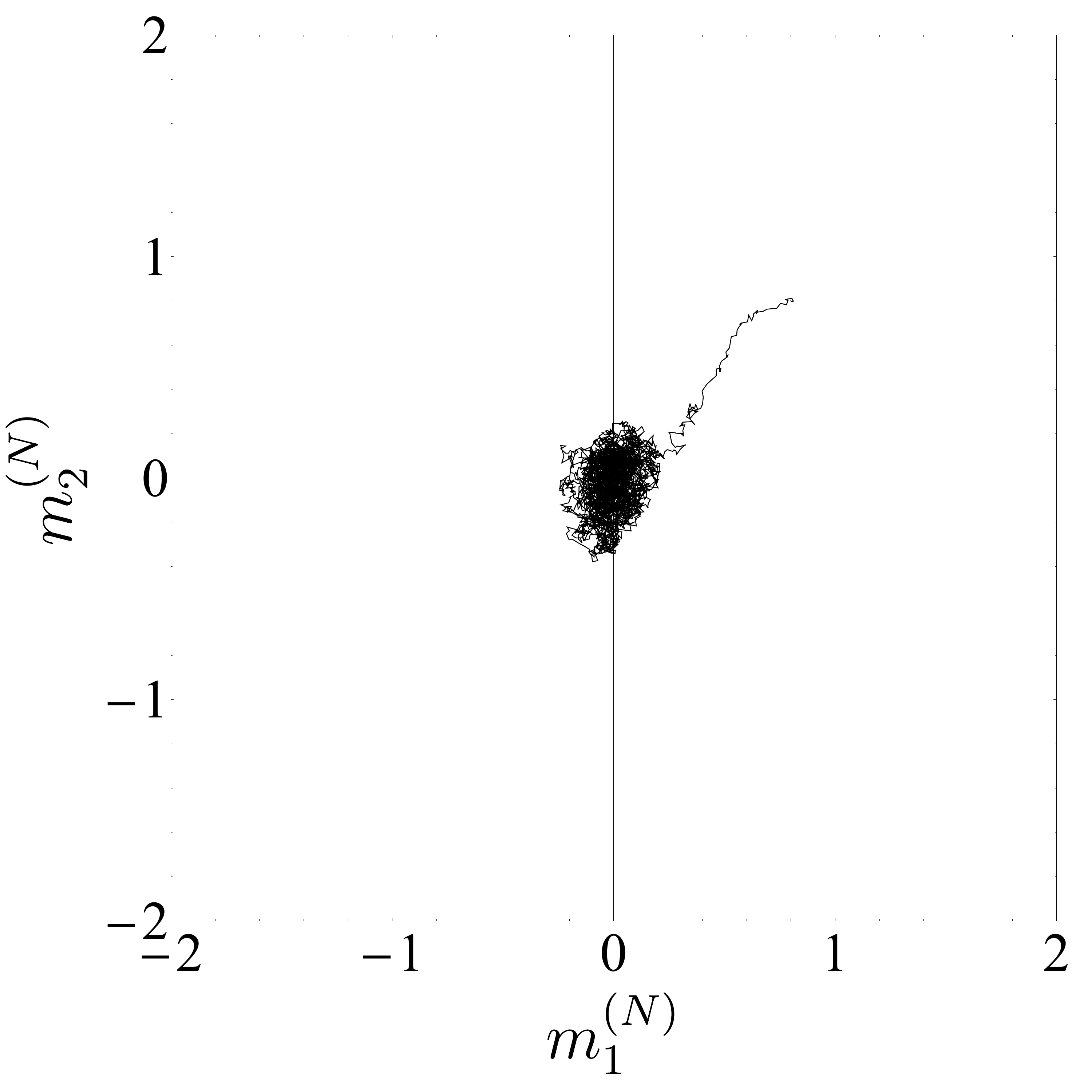}}
    \end{center}
    
    \vspace{5pt}
    \begin{center}
        \subfigure{\includegraphics[width=0.325\textwidth]{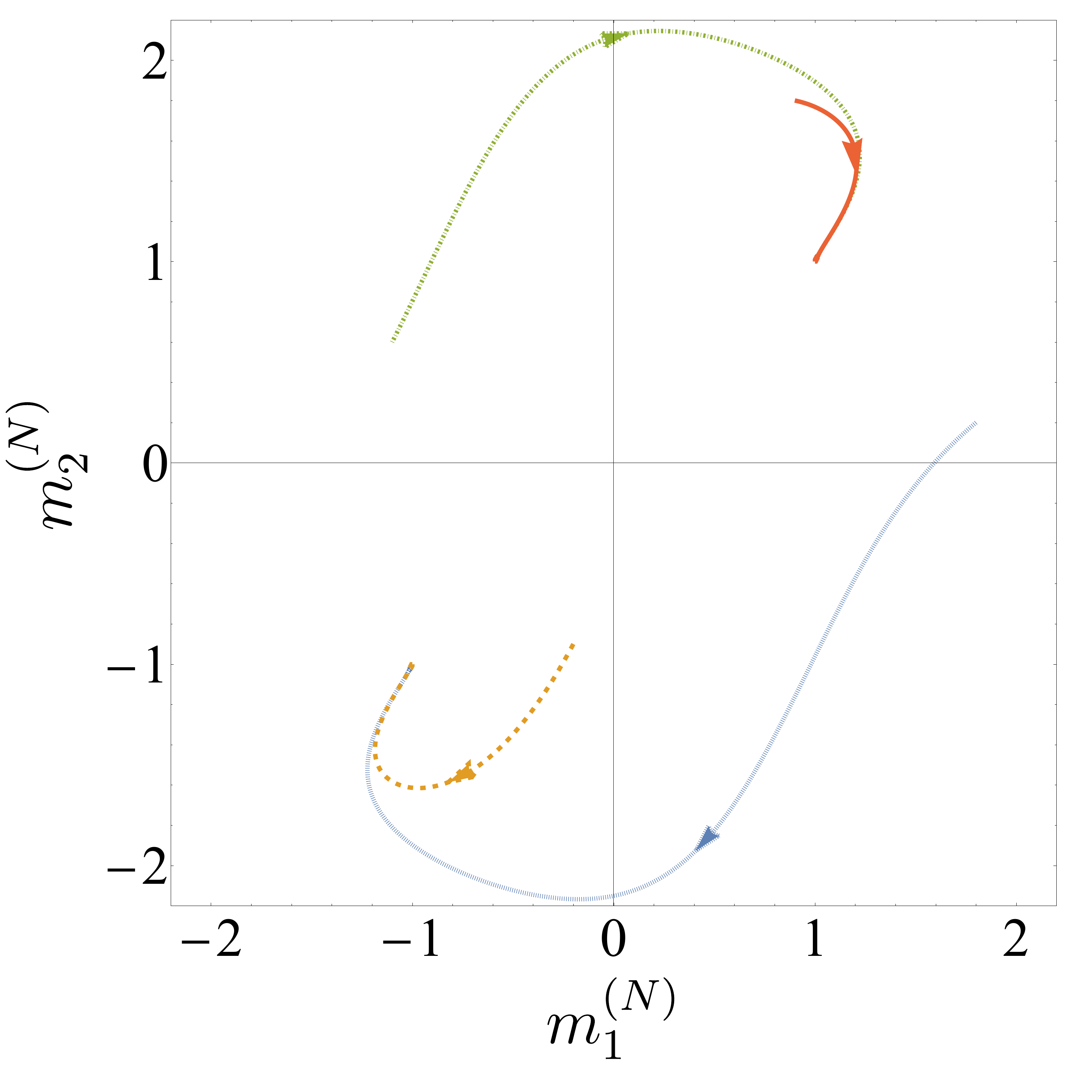}}
        \subfigure{\includegraphics[width=0.325\textwidth]{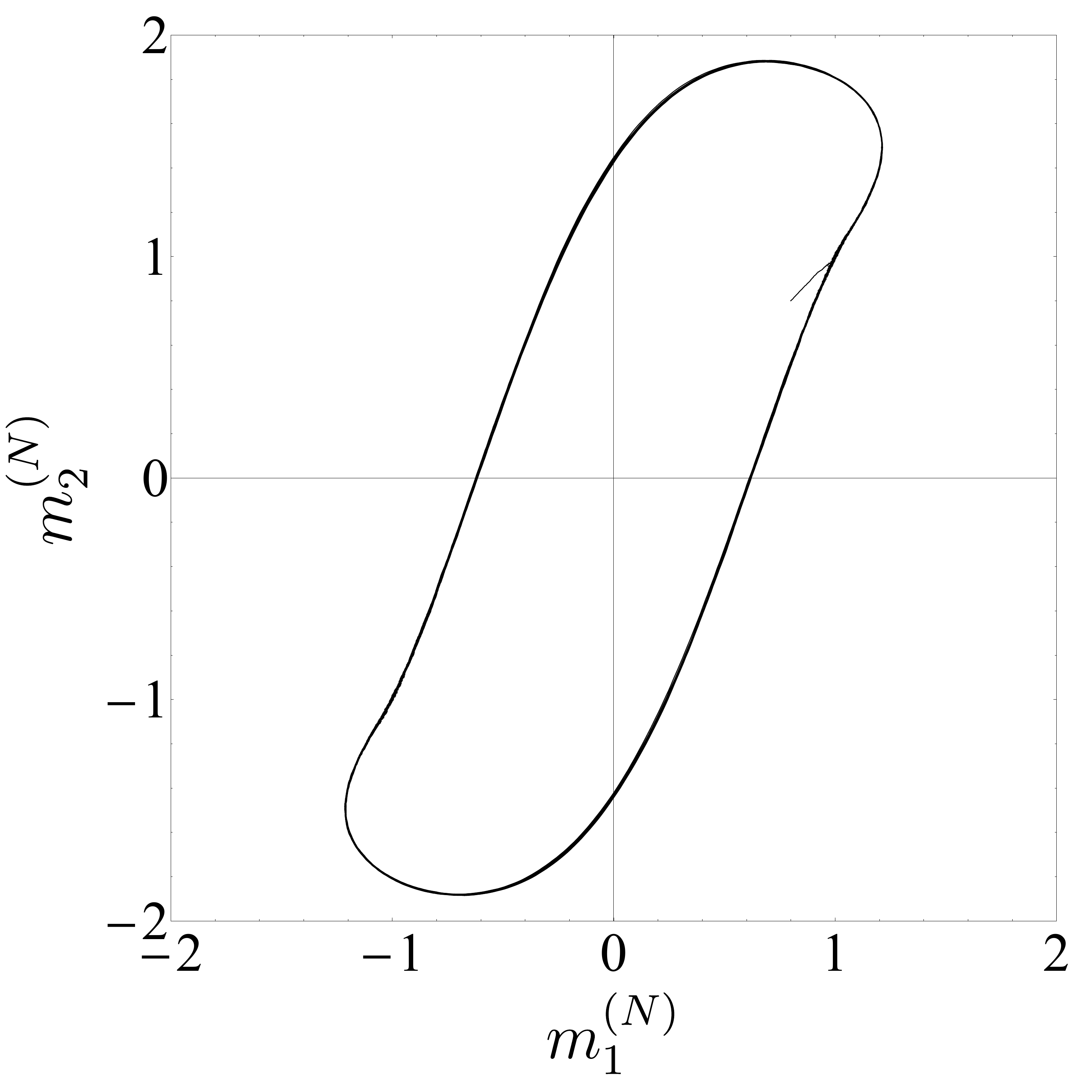}}
        \subfigure{\includegraphics[width=0.325\textwidth]{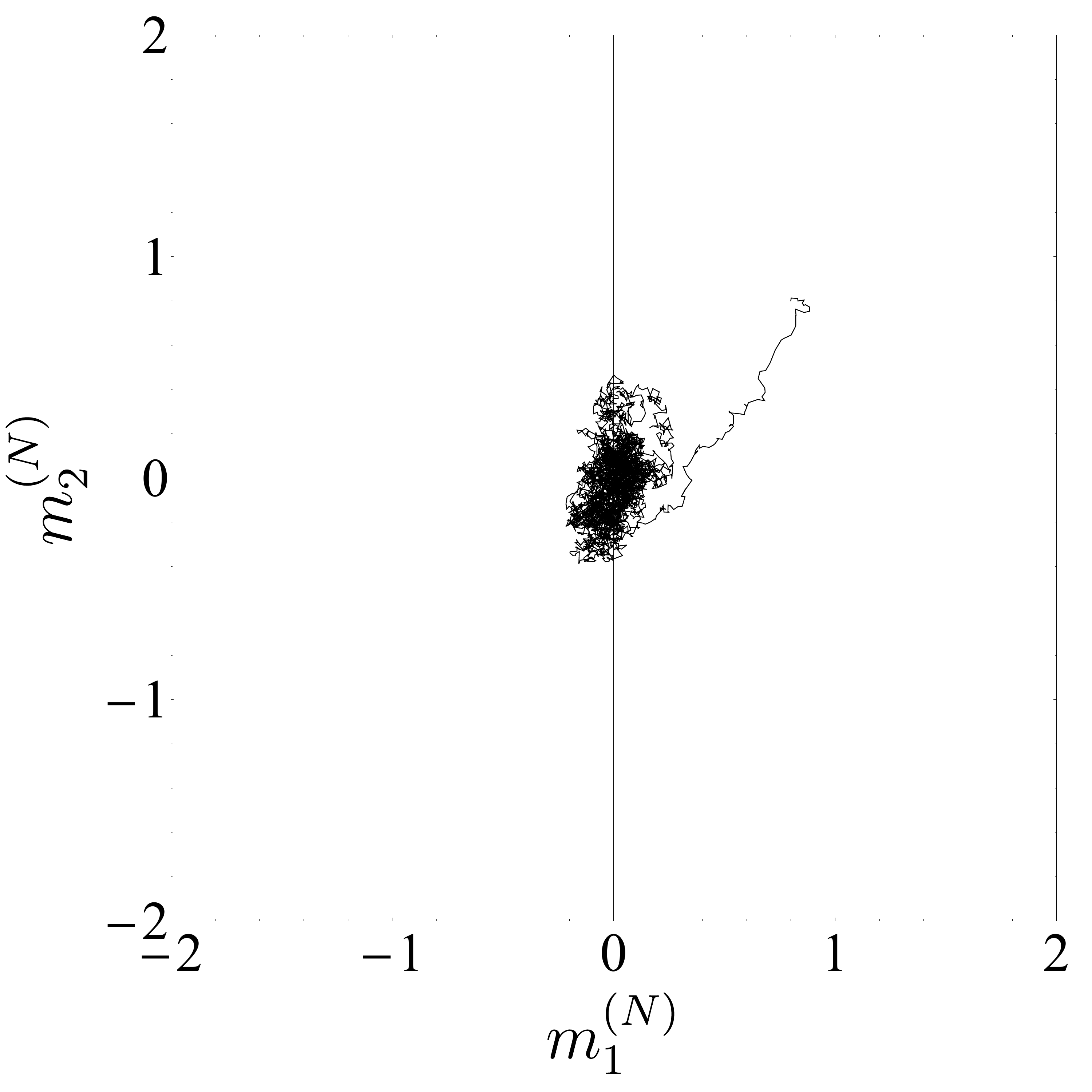}}
    \end{center}
    
    \vspace{5pt}
    \begin{center}
        \subfigure{\includegraphics[width=0.325\textwidth]{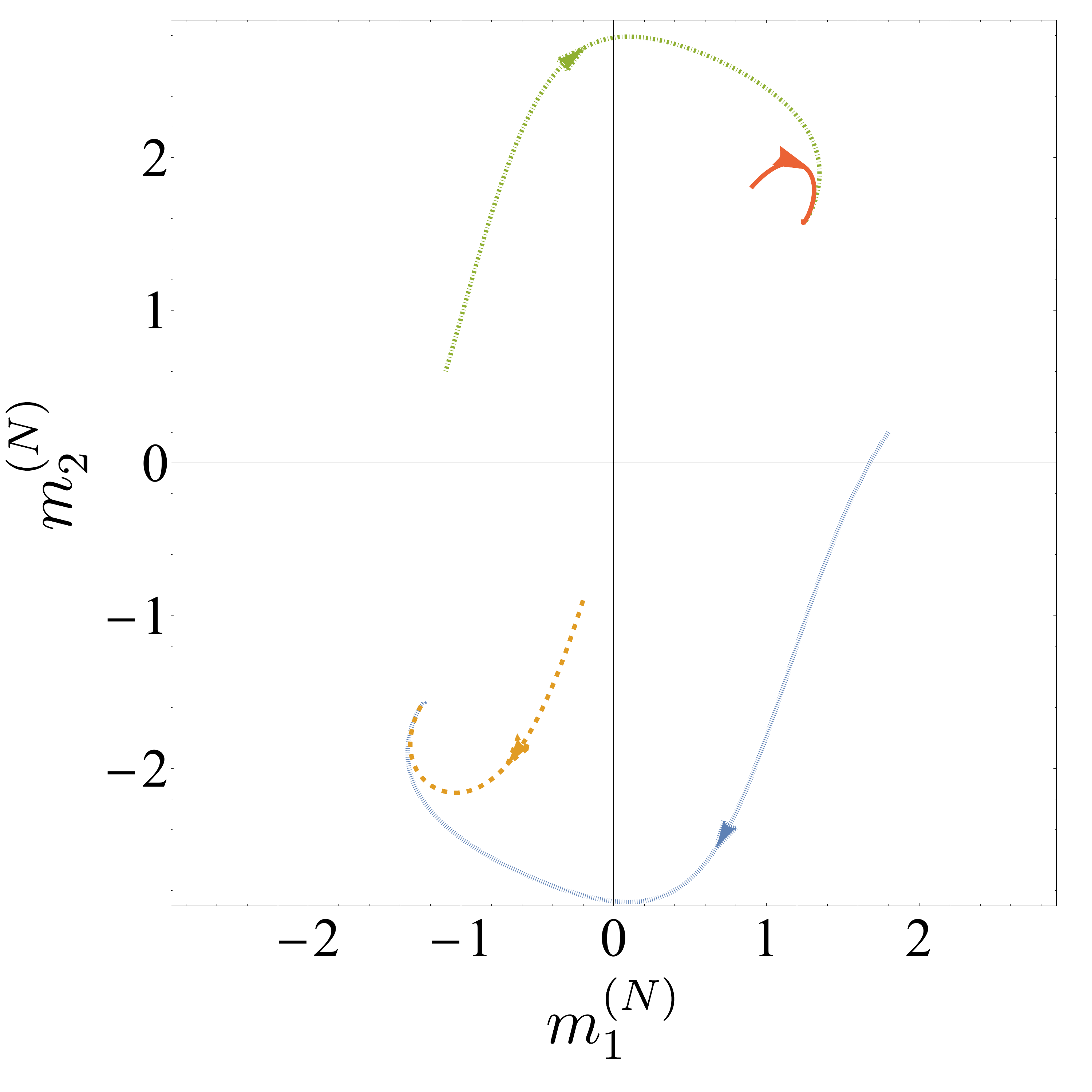}}
        \subfigure{\includegraphics[width=0.325\textwidth]{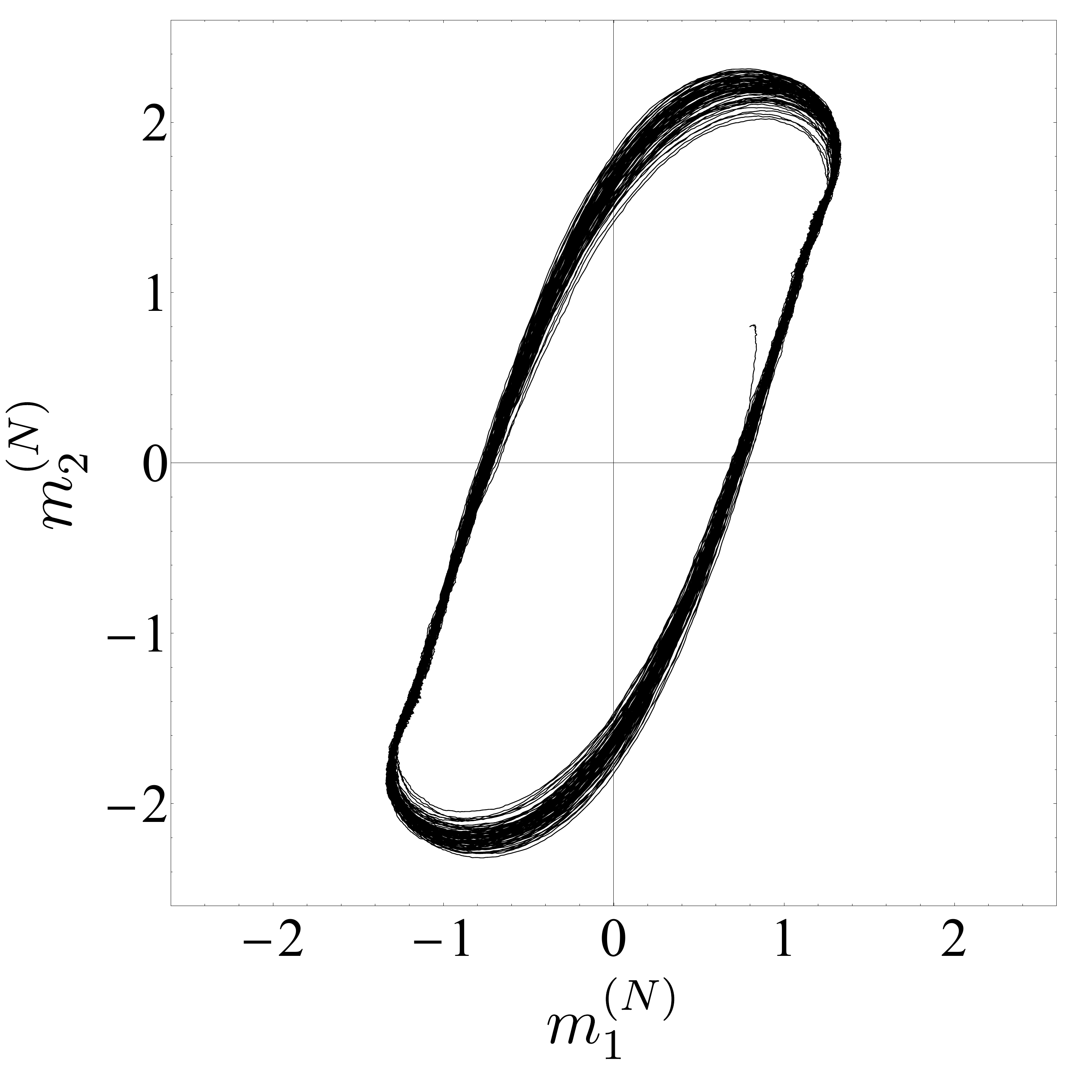}}
        \subfigure{\includegraphics[width=0.325\textwidth]{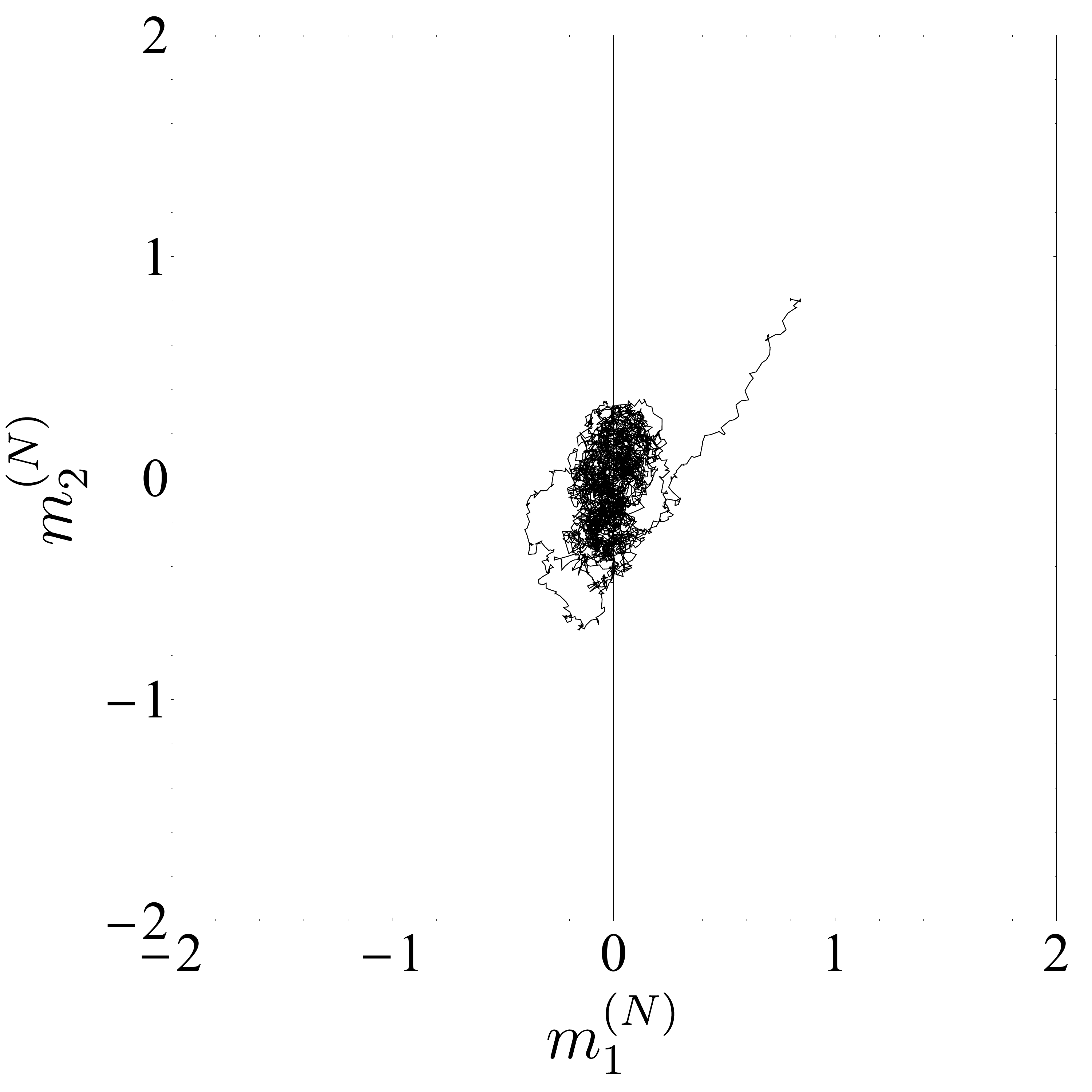}}
    \end{center}
    \caption{\footnotesize{Trajectories of $\left(m_1^{(N)}(t),m_2^{(N)}(t)\right)$ obtained with numerical simulations of system \eqref{micro:dyn}, in the absence of noise (first column), in the presence of an intermediate amount of noise (second column) and of a high-intensity noise (third column). In all cases, we considered $10^6$ iterations with a time-step $dt=0.005$, $1000$ particles, $\alpha = 0.5$, $\theta_{11}=\theta_{22}=8$. From top to bottom: $A-1<B<A+2$, in particular, $A=2$ and $B=2.5$; $B=A+2$, in particular, $A=2$ and $B=4$; $B>A+2$, in particular, $A=2$ and $B=7$.\\
    We see that, during a time interval of the same length (namely, $10^6$ iterations), when the intensity of the noise is below a certain threshold (first column, $\sigma=0$ in all the three panels) no periodic behavior arises in any of the three considered cases and the system ends up in one of the stable equilibria. On the contrary, when the intensity of the noise is large (third column, $\sigma=5$ in all the three panels), the zero-mean Brownian disturbance dominates and the trajectories resemble random excursions around the origin. Whenever the amount of noise is intermediate (second column, from top to bottom: $\sigma=0.5$, $\sigma=0.1$ and $\sigma=0.6$), self-sustained oscillations appear; for further details about this scenario see Fig.~\ref{fig:PhDiaNoise}. 
 }}
    \label{fig:small_vs_big_noise}
\end{figure}

\begin{figure}[h!p]
\captionsetup[subfigure]{labelformat=empty}
    \centering
    \begin{center}
    \subfigure{\includegraphics[width=0.33\textwidth]{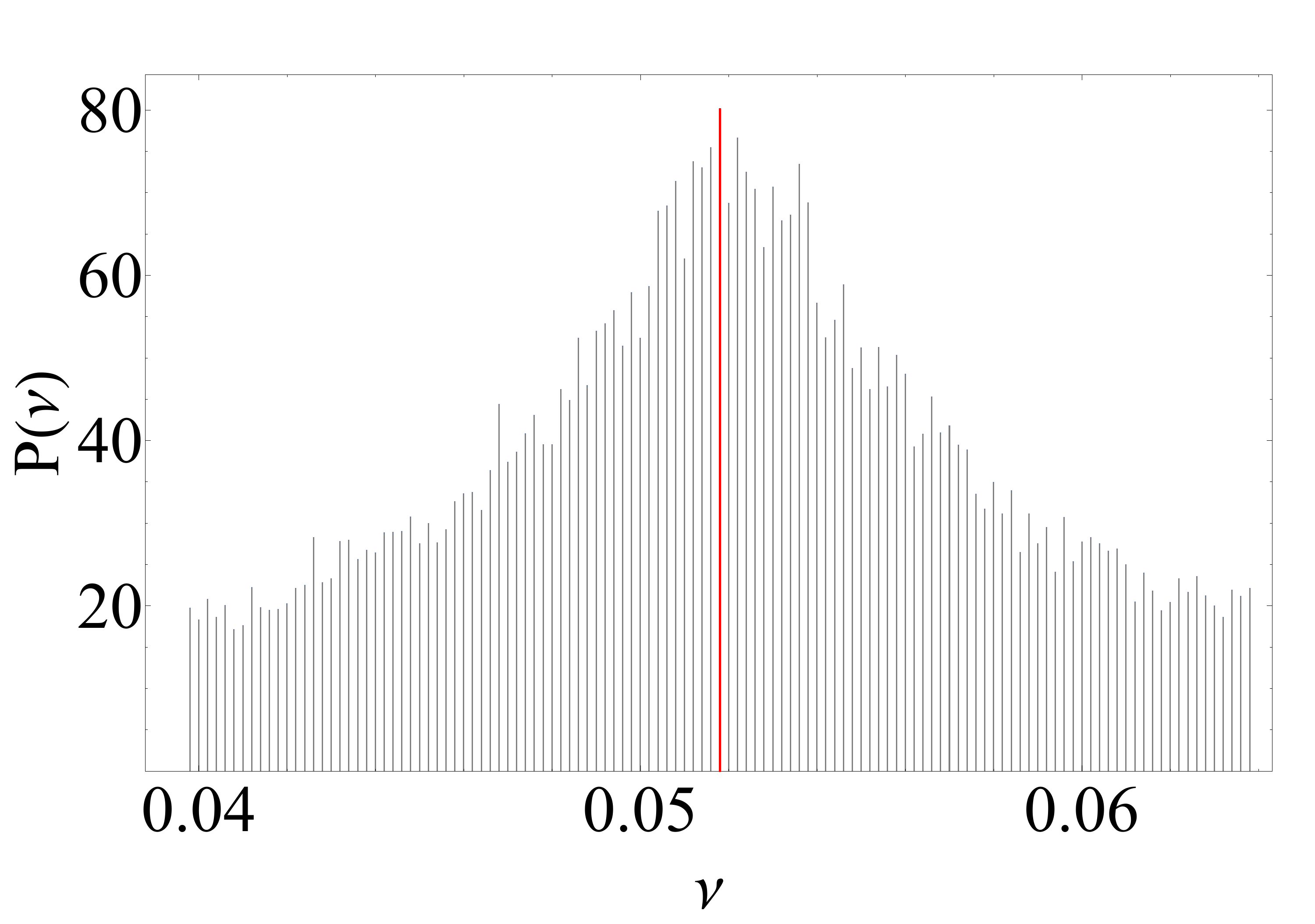}}
    \subfigure{\includegraphics[width=0.33\textwidth]{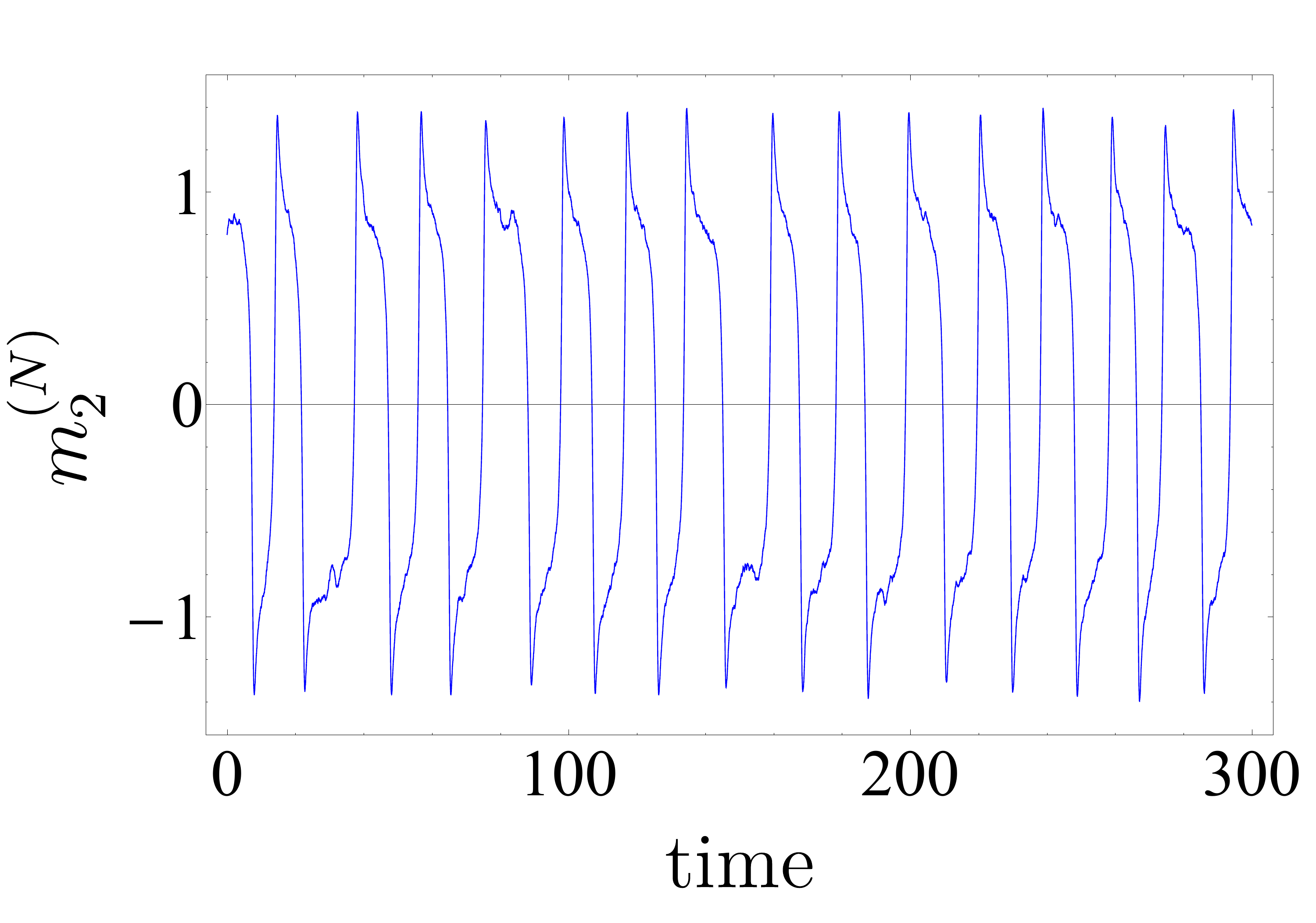}}
    \subfigure{\includegraphics[width=0.32\textwidth]{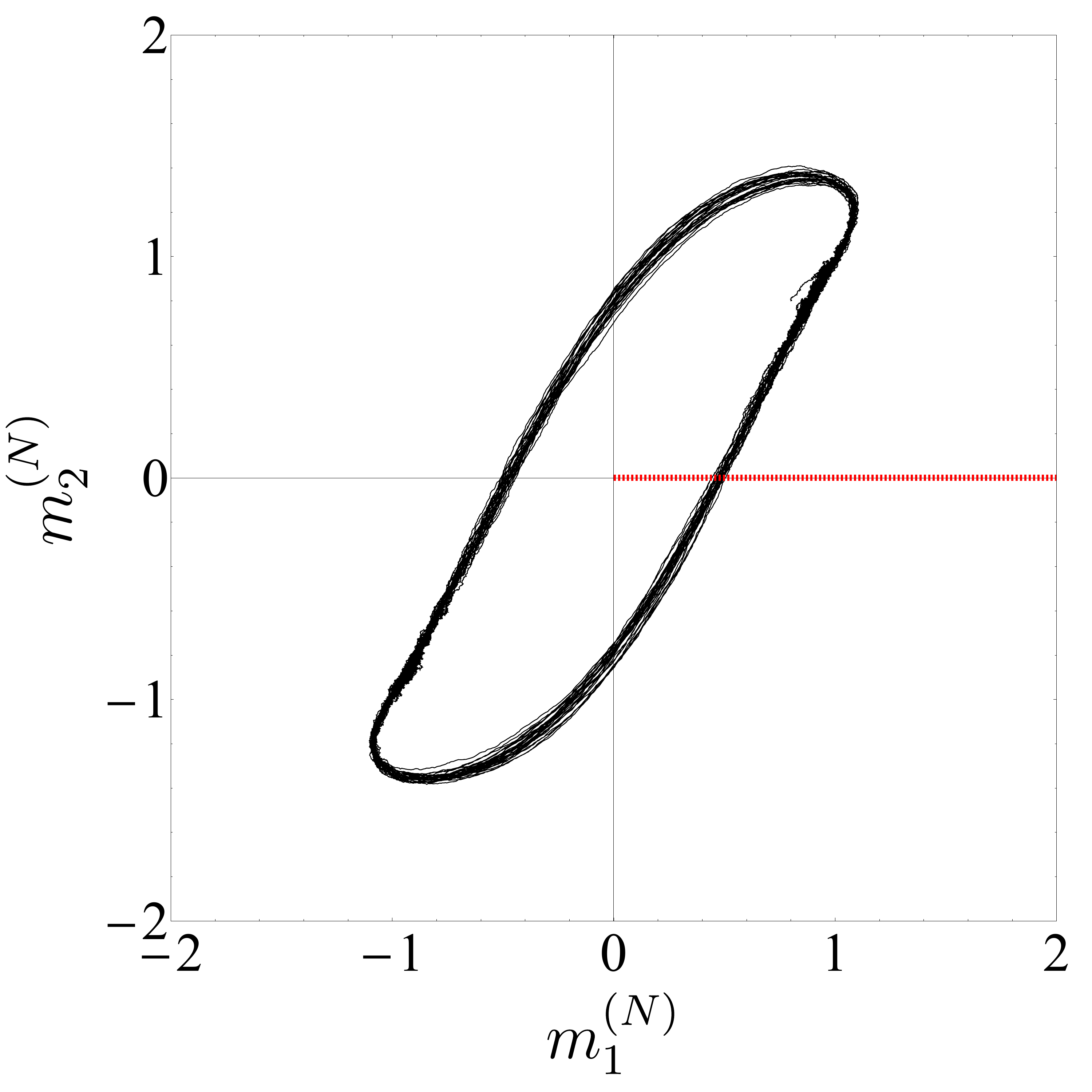}}
    \end{center}
    
    \vspace{5pt}
    \begin{center}
    \subfigure{\includegraphics[width=0.33\textwidth]{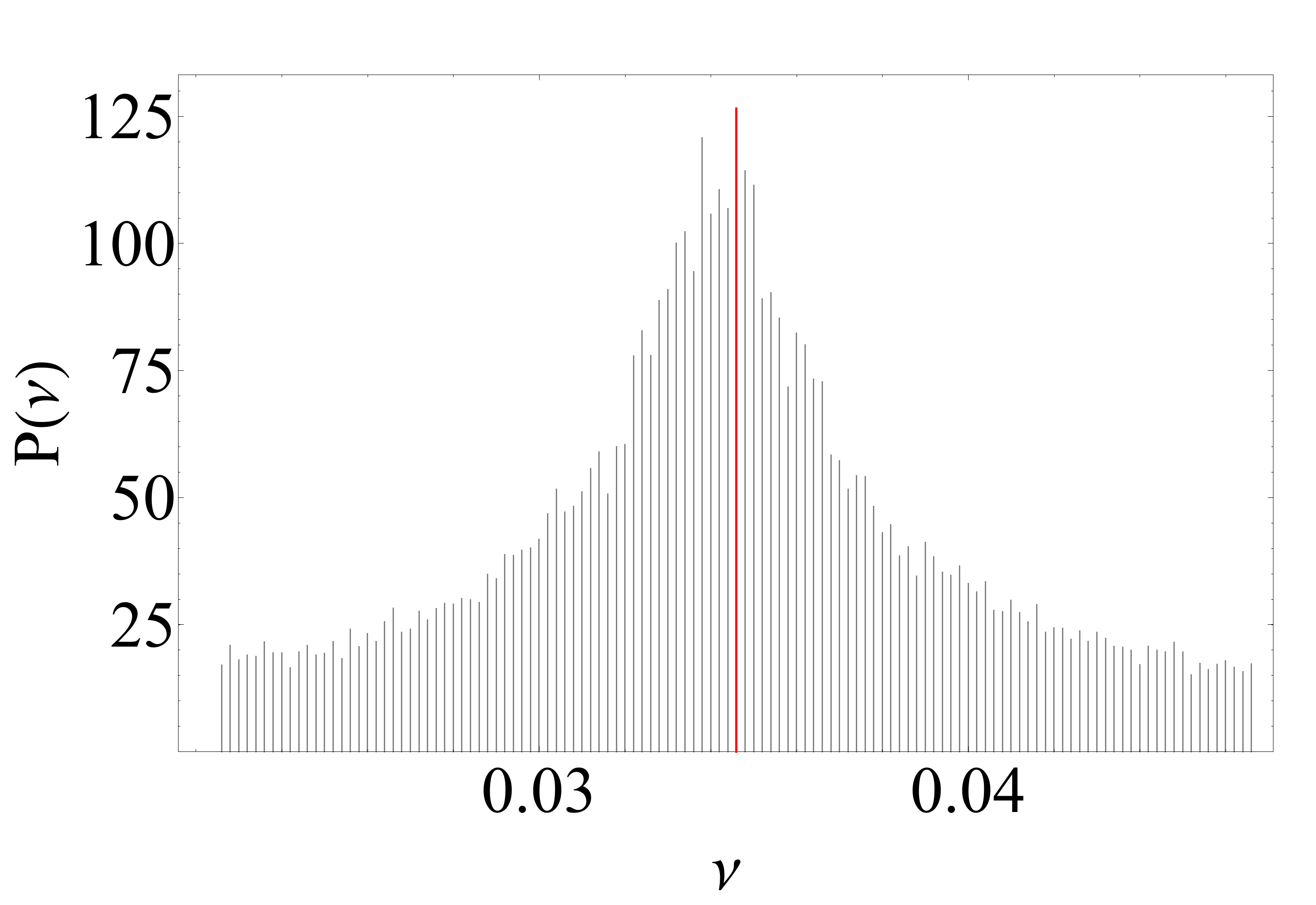}}
    \subfigure{\includegraphics[width=0.33\textwidth]{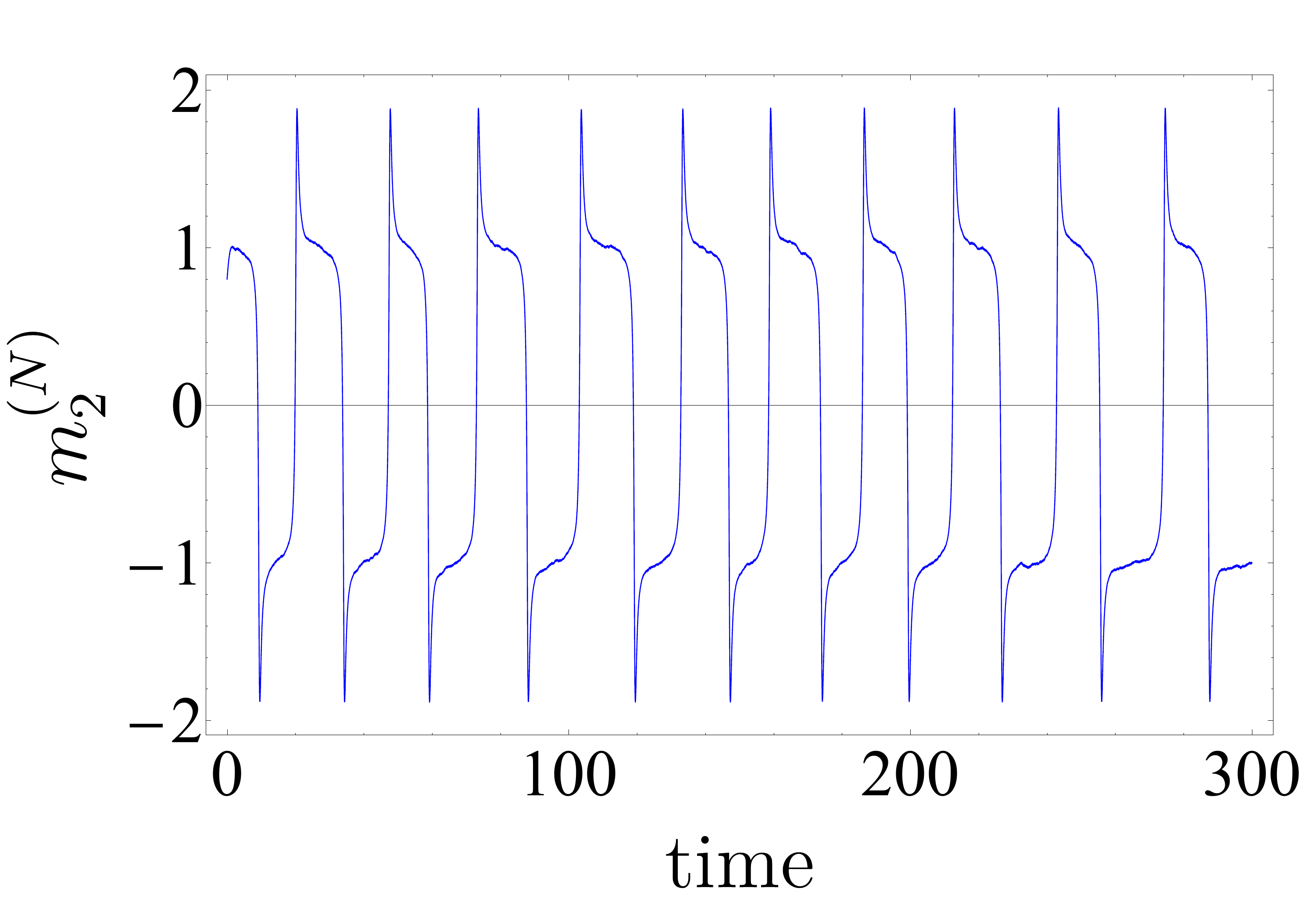}}
    \subfigure{\includegraphics[width=0.32\textwidth]{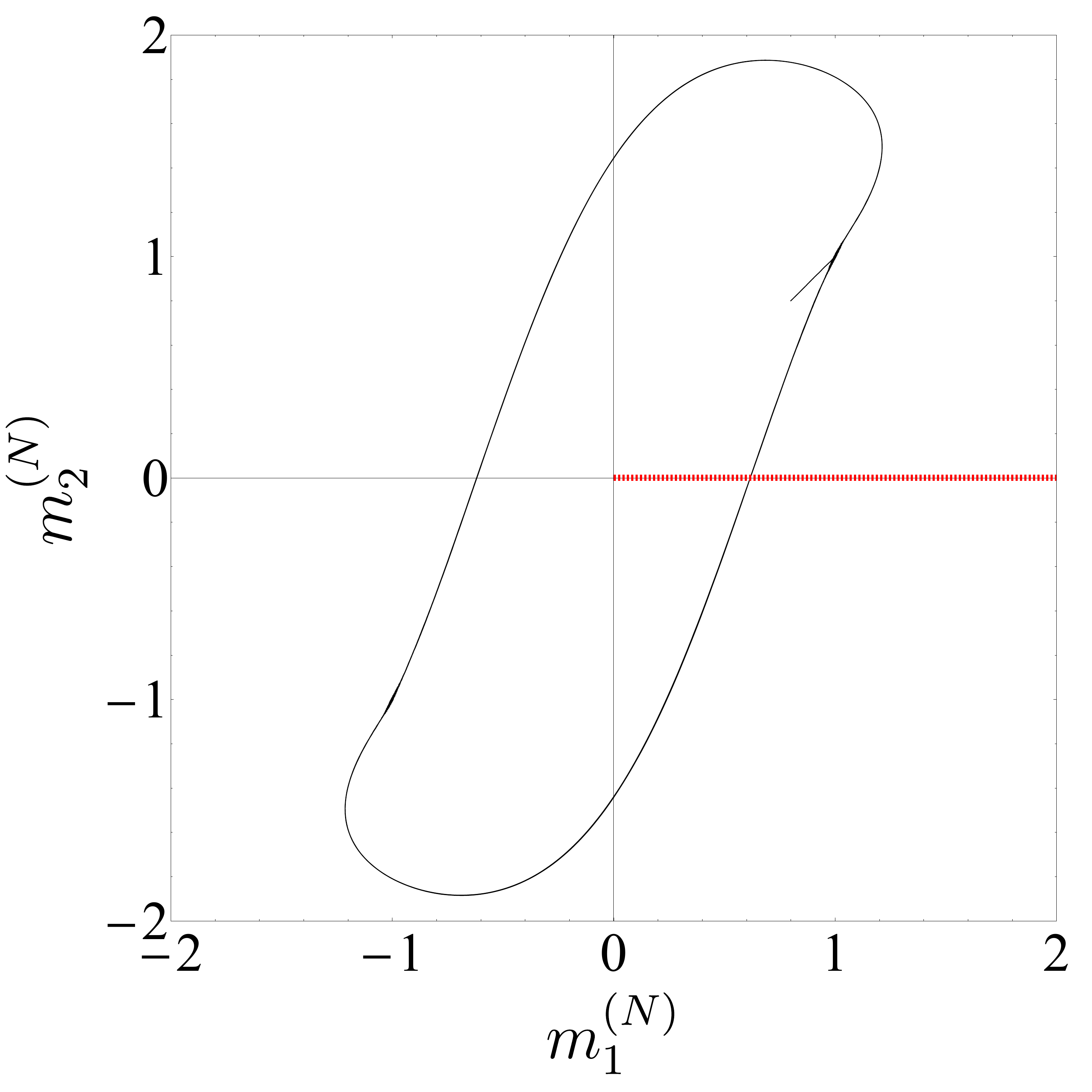}}
    \end{center}
    
    \vspace{5pt}
    \begin{center}
    \subfigure{\includegraphics[width=0.33\textwidth]{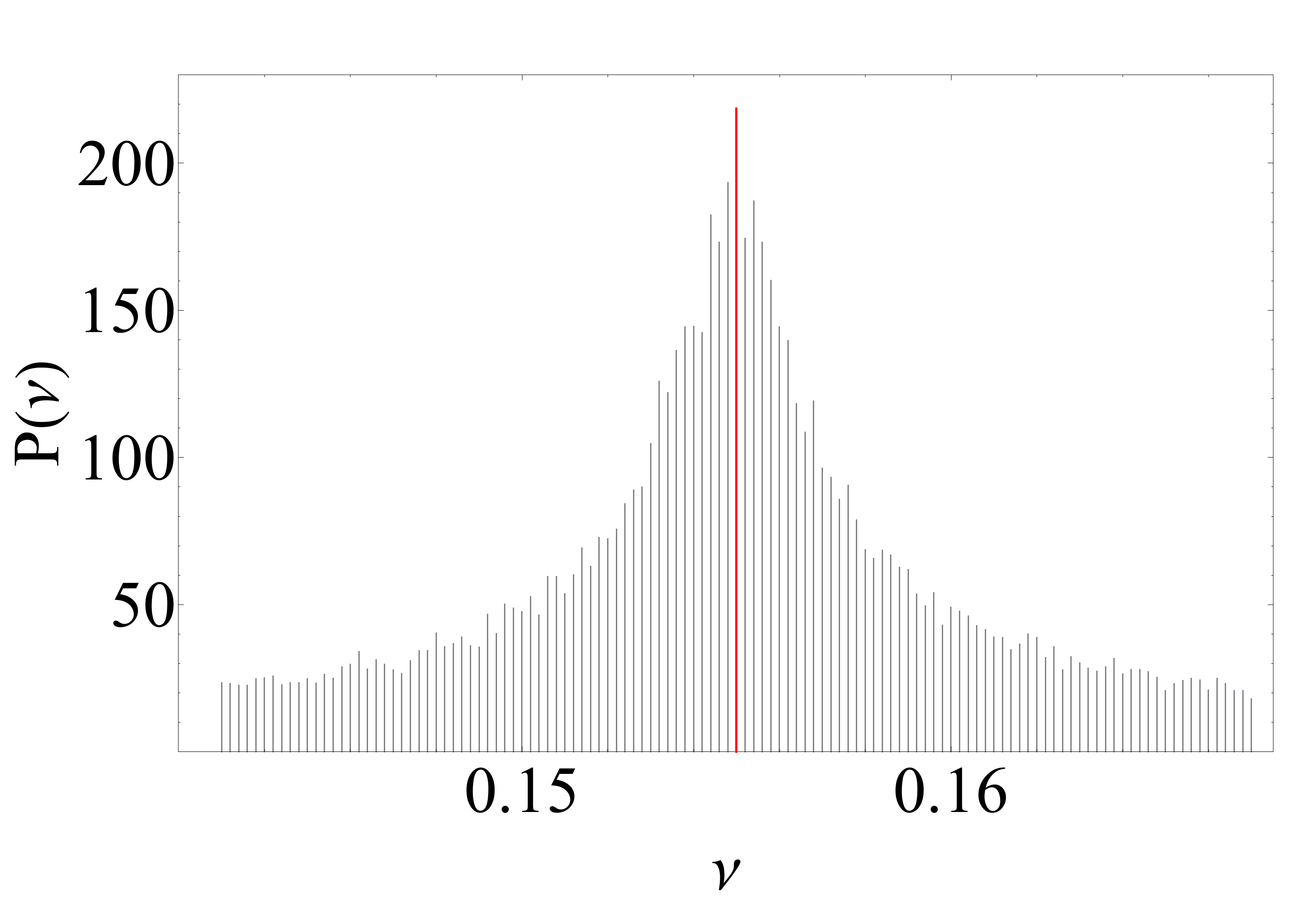}}
    \subfigure{\includegraphics[width=0.33\textwidth]{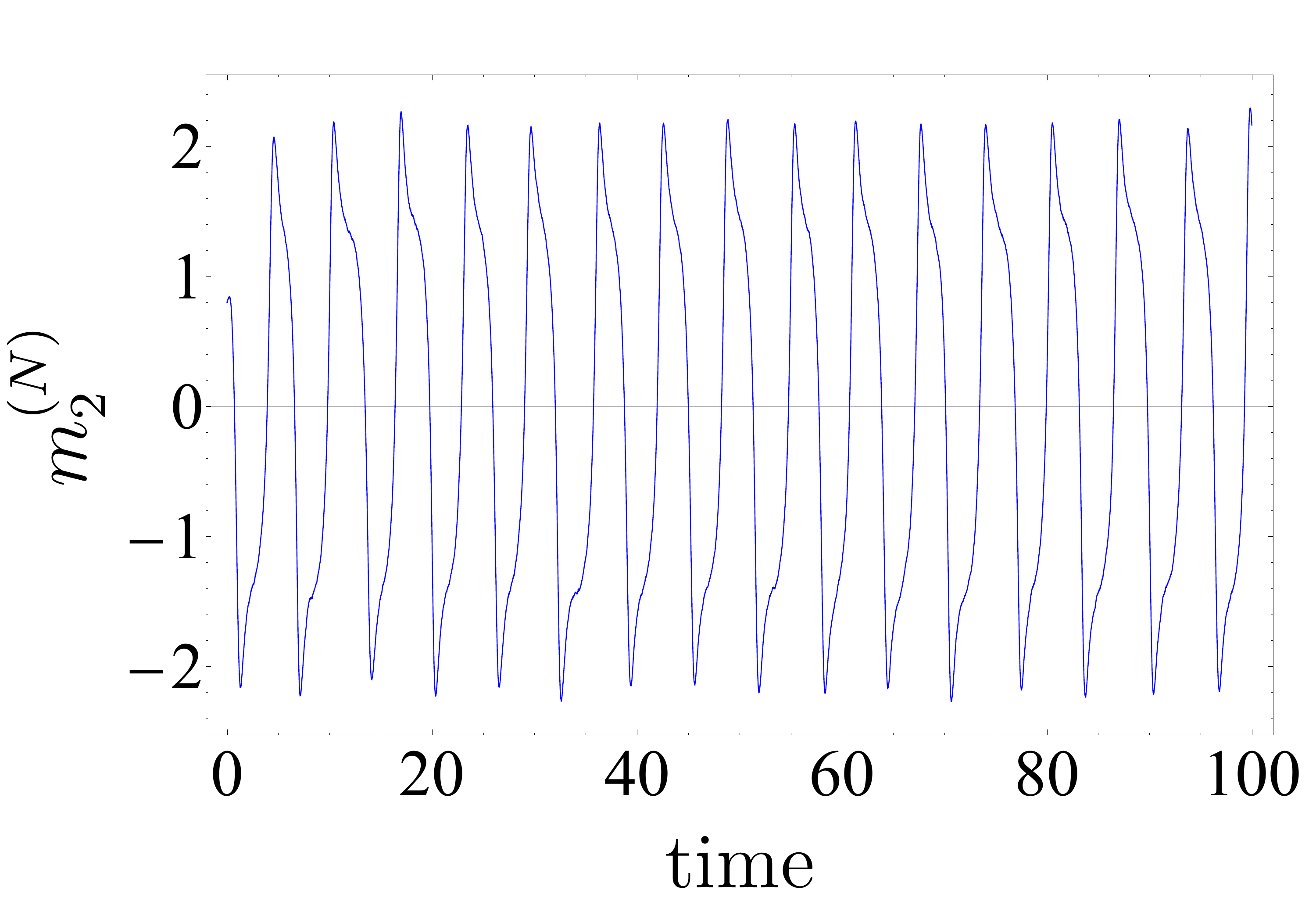}}
    \subfigure{\includegraphics[width=0.32\textwidth]{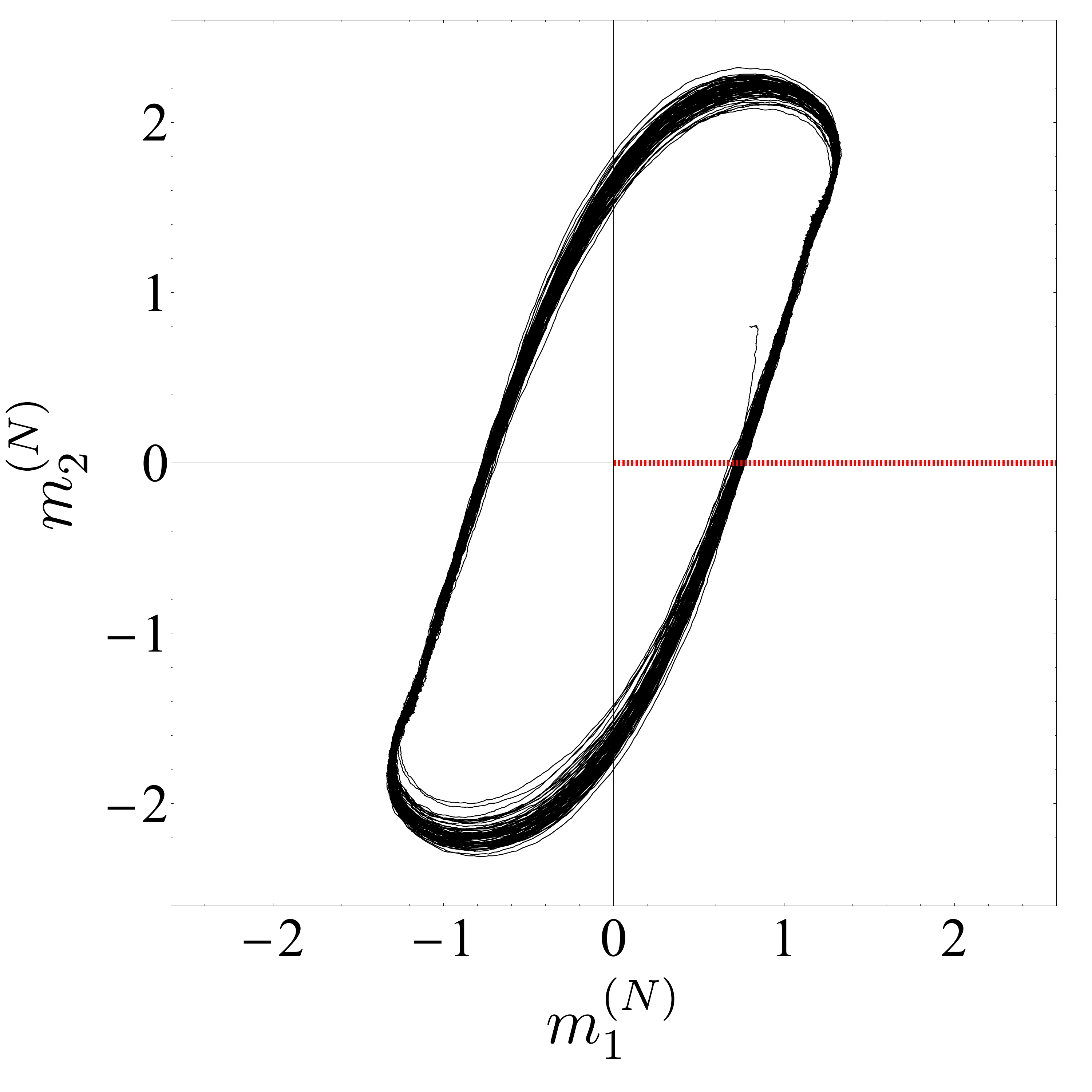}}
    \end{center}
    \caption{\footnotesize{Analysis of the period of the trajectories of $\left(m_1^{(N)}(t),m_2^{(N)}(t)\right)$ obtained via numerical simulations of system \eqref{micro:dyn}, in the presence of an intermediate amount of noise. In all cases, we considered $10^6$ iterations with a time-step $dt=0.005$, $1000$ particles, $\alpha = 0.5$, $\theta_{11}=\theta_{22}=8$. From top to bottom: $A-1<B<A+2$, in particular, $A=2$ and $B=2.5$, with $\sigma=0.5$; $B=A+2$, in particular, $A=2$ and $B=4$, with $\sigma=0.1$; $B>A+2$, in particular, $A=2$ and $B=7$, with $\sigma=0.6$.\\
    In the first column, we plotted the relevant spectral region of the averaged modulus of the discrete Fourier transform $P\left(\nu\right)$ of $m^{(N)}_2$ against the frequencies $\nu$. 
    For these figures we employed the Fourier function of Mathematica 
    applied to a trajectory of $m^{(N)}_2$ over $10^{6}$ steps and averaged the obtained spectrum over $M=50$ simulations. The average periods in the three cases were obtained as the reciprocals of the frequencies highlighted by the red peaks. 
    In the second column, we plotted the time evolution of $m^{(N)}_2$.
    The third column shows a trajectory of $\left(m^{(N)}_1(t),m^{(N)}_2(t)\right)$.  
There, red dashed horizontal lines mark the Poincar{\'e} sections we employed for the computation of the average period.}}
\label{fig:PhDiaNoise}
\end{figure}

In this section, we present numerical simulations of the finite-size system \eqref{micro:dyn}, aimed at giving evidences of the phenomenon of noise-induced periodicity.\\
In the setting introduced in Section~\ref{sec:model:results}, we ran several simulations of \eqref{micro:dyn} for different choices of $\sigma$ and several values of the interaction strengths. In all cases, we performed simulations with $10^{6}$ iterations with time-step $dt = 0.005$ for a system of $1000$ particles equally divided between the two populations ($\alpha=0.5$). All particles in the same population were given the same initial condition. We fixed $\theta_{11}=\theta_{22}=8$ and let $A\coloneqq \left(1-\alpha\right)\theta_{12}>0$ and $B\coloneqq -\alpha \theta_{21}>0$ vary. The results are displayed in Fig.~\ref{fig:small_vs_big_noise}, Fig.~\ref{fig:PhDiaNoise} and Table~\ref{tab:comparison_fft_poincare}, where also the specific values we employed  for $A$, $B$ and $\sigma$ are reported. The choices of the parameters are discussed in more detail in Section~\ref{sec:M2}, as they correspond to different regimes of the limiting noiseless dynamics (i.e., system \eqref{prop:chaos} with $\sigma=0$), namely, $A-1<B<A+2$, $B=A+2$ and $B>A+2$.\\
We observe the following: 
\begin{enumerate}
\item If $\sigma=0$ the system is attracted to a fixed point (see the first column of Fig.~\ref{fig:small_vs_big_noise}). Numerical evidences support the idea that, in the regimes $A-1 < B < A+2$ and $B > A+2$, this behavior persists for small $\sigma>0$.
\item When the intensity of the noise is tuned to an intermediate range of values, an oscillatory behavior is observed in the $\left(m^{(N)}_{1}, m^{(N)}_{2}\right)$ plane throughout the duration of the simulation, suggesting the presence of a periodic law (see the second column of Fig.~\ref{fig:small_vs_big_noise}). Thus, our model seems to exhibit noise-induced periodicity. This phenomenon, which at the best of our knowledge lacks a full theoretical comprehension, can be loosely described in the following terms: an intermediate amount of noise may create/stabilize some attractors and destabilize others. In our case it seems that the noise destabilizes (some of the) fixed points and generates a stable rhythmic behavior of the empirical averages of the particle positions of the two communities.\\
We would like to mention that, in the regime $A=B+2$, an arbitrarily small value of $\sigma>0$ seems to be sufficient to induce periodicity.
\item Letting $\sigma \gg 1$ completely alters the dynamics that essentially becomes a Brownian motion (see the third column of Fig.~\ref{fig:small_vs_big_noise}).
\end{enumerate}
\begin{table}[H] 
    \centering
    \begin{tabular}{|c|c|c|c|}
    \hline
     \cellcolor{LightGray} Noise & \cellcolor{LightGray} Coupling strengths & \cellcolor{LightGray} Period of Poincar{\'e} map & \cellcolor{LightGray} Period fft  \\
      \hline\hline
     $\sigma=0.5$ & $A=2$, $B=2.5$ & $19.35\pm 0.18$ & $19.31 \pm 0.77$\\
     $\sigma=0.1$ & $A=2$, $B=4$ & $29.34 \pm 0.29$ & $28.90 \pm 0.92$\\
     $\sigma=0.6$ & $A=2$, $B=7$ & $6.45 \pm 0.01$ & $6.45 \pm 1.28$\\
     \hline
    \end{tabular}
    \caption{\footnotesize{Period of the rhythmic oscillations of system \eqref{micro:dyn} in the $\left(m^{(N)}_1,m^{(N)}_2\right)$ plane in the various regimes and in the presence of an intermediate amount of noise. First and second column: values of the parameters. Third column: period of $t \mapsto \left(m^{(N)}_1(t),m^{(N)}_2(t)\right)$ obtained by computing the average passage time from positive zero to negative zero of $m^{(N)}_2$; i.e., the average return time to the Poincar{\'e} section $\left\{m^{(N)}_2=0, m^{(N)}_1>0\right\}$. Fourth column: period of $t \mapsto m_2^{(N)}(t)$  estimated from the Fourier spectra considering a sampling period equal to $dt=0.005$.}}
    \label{tab:comparison_fft_poincare}
\end{table}

In Fig.~\ref{fig:PhDiaNoise} and Table~\ref{tab:comparison_fft_poincare} the oscillatory behavior emerging in system \eqref{micro:dyn} is analyzed further. We computed the average return time of the system to the Poincar{\'e} section $\left\{m^{(N)}_2=0, m^{(N)}_1>0\right\}$ and its standard deviation, in the various regimes. These are reported in the third column of Table~\ref{tab:comparison_fft_poincare}. The Poincar\'e section is plotted as a red line in Fig.~\ref{fig:PhDiaNoise}. In addition, we computed the discrete Fourier transform, averaged over $M=50$ simulations, for the average particle position of the second population, $m^{(N)}_2$. From the peak of the Fourier transform we recovered the period of the trajectory of $m^{(N)}_2(t)$. The average period and its standard deviation are reported in the fourth column of Table~\ref{tab:comparison_fft_poincare} for different values of the parameters. 

\section{Propagation of chaos and small-noise approximation}
\label{sec:M2}

In this section we give our main results. We begin with a propagation of chaos statement, allowing to get the macroscopic description \eqref{prop:chaos} of our system. Then, we analyze the noiseless version of the macroscopic dynamics and we show the absence of limit cycles as attractors. Finally, in a small-noise regime, we derive a Gaussian approximation of the infinite volume evolution \eqref{prop:chaos} that displays an oscillatory behavior.

\subsection{Propagation of chaos}\label{sec:prop:chaos}
\noindent  Propagation of chaos claims that, as $N\to\infty$, the evolution of each particle remains independent of the evolution of any finite subset of the others. This is coherent with the fact that individual units interact only through the empirical means of the two populations, over which the influence of a finite number of particles becomes negligible when taking the infinite volume limit. In our case the limiting evolution of a pair of representative particles, one for each population, is the process $((x(t),y(t)), 0 \leq t \leq T)$ described by the stochastic differential equation~\eqref{prop:chaos}.

Under the assumptions $\mathbb{E}[x(0)]<\infty$ and $\mathbb{E}[y(0)]<\infty$, it is easy to prove that system \eqref{prop:chaos} has a unique strong solution (see Theorem~\ref{thm:prop:chaos:existence} in Appendix~\ref{appendixA}). Moreover, by a coupling argument, we obtain the following theorem. 

\begin{thm}\label{thm:prop:caos}
Fix $T>0$. Let $\left(\left(x^{(N)}_{1}(t),\dots,x^{(N)}_{N_1}(t), y^{(N)}_{1}(t),\dots,y^{(N)}_{N_2}(t)\right), 0 \leq t \leq T \right)$ be the solution to Eq.~\eqref{micro:dyn} with an initial condition satisfying the following requirements:
\begin{itemize}
\item the collection $\left(x^{(N)}_{1}(0), \dots, x^{(N)}_{N_1}(0), y^{(N)}_{1}(0), \dots, y^{(N)}_{N_2}(0)\right)$ is a family of independent random variables.
\item the random variables $\left(x^{(N)}_{1}(0), \dots, x^{(N)}_{N_1}(0)\right)$ (resp.  $\left(y^{(N)}_{1}(0), \dots, y^{(N)}_{N_2}(0)\right)$) are identically distributed with law $\lambda_x$ (resp. $\lambda_y$).
We assume that $\lambda_x$ and $\lambda_y$ have finite second moment.
\item the random variables $x^{(N)}_{j}(0)$ and $y^{(N)}_{k}(0)$ are independent of the Brownian motions $\left(w_{i}(t), 0 \leq t \leq T\right)_{i=1,\dots,N}$ for all $j=1,\dots, N_1$ and $k=1,\dots,N_2$.
\end{itemize}
Moreover, let $\left(\left(x_{1}(t), \dots, x_{N_1}(t), y_{1}(t), \dots, y_{N_2}(t)\right), 0 \leq t \leq T \right)$ be the process whose entries are independent and such that $(x_j(t), 0 \leq t \leq T)_{j=1, \dots, N_1}$ (resp. $(y_k(t), 0 \leq t \leq T)_{k=1, \dots, N_2}$) are  copies of the solution to the first (resp. second) equation in \eqref{prop:chaos}, 
with the same initial conditions and the same Brownian motions used to define system \eqref{micro:dyn}. Here, \lq\lq the same\rq\rq\, means component-wise equality.\\  
Define the index sets $\mathcal{I}=\{i_1,\ldots,i_{k_1}\}\subseteq \{1,\ldots,N_1\}$, with $|\mathcal{I}|=k_1$, and  $\mathcal{J}=\{j_1,\ldots,j_{k_2}\}\subseteq \{1,\ldots,N_2\}$, with $|\mathcal{J}|=k_2$.
Then, we have
\begin{equation}\label{prop_chaos_on_average}
   \lim_{N\to+\infty}\mathbb{E}\left[\sup_{t\in [0,T]}\left|\mathbf{z}^{(N)}_{k_1,k_2}(t) -\mathbf{z}_{k_1,k_2}(t) \right| \right] = 0,
\end{equation}
with $\left|\mathbf{z}\right|$ the $\ell^1$-norm of a vector $\mathbf{z}$, $\mathbf{z}^{(N)}_{k_1,k_2}(t) = \left(x^{(N)}_{i_1}(t), \dots, x^{(N)}_{i_{k_1}}(t), y^{(N)}_{j_1}(t), \dots, y^{(N)}_{j_{k_2}}(t)\right)$ and \(\mathbf{z}_{k_1,k_2}(t) =\left(x_{1}(t),\dots,x_{k_1}(t), y_{1}(t),\dots,y_{k_2}(t)\right)\).
\end{thm}

The proof of Theorem~\ref{thm:prop:caos} is postponed to Appendix~\ref{appendix:B}. Recall that the convergence in Theorem~\ref{thm:prop:caos} implies, for $t \in [0,T]$, convergence in distribution of any finite-dimensional vector $\mathbf{z}^{(N)}_{k_1,k_2}(t)$ to $\mathbf{z}_{k_1,k_2}(t)$.

\subsection{Analysis of the zero-noise dynamics}\label{sec:zero:noise:dyn}
In this section we consider system \eqref{prop:chaos} with $\sigma=0$.
Notice that, in the zero-noise version of \eqref{prop:chaos}, 
the terms $\alpha\theta_{11}\left(x- \mathbb{E}[x]\right)$ and $(1-\alpha)\theta_{22}\left(y-\mathbb{E}[y]\right)$ are both zero.  Thus, setting 
\[
A\coloneqq \left(1-\alpha\right)\theta_{12}>0 \qquad \text{ and } \qquad B\coloneqq -\alpha \theta_{21}>0,
\]
system \eqref{prop:chaos} reduces to
\begin{align}\label{micro:dyn:nonoise:AB}
\dot{x} &= -x^3 +x -A\left(x-y\right)\nonumber\\
\dot{y} &= -y^3 +y -B\left(x-y\right).
\end{align}

We make the following assumption at this point. We will focus on the case
\begin{itemize}
\item[\textbf{(H)}] $A>1$ and $B>A-1$.
\end{itemize}
The reason for this choice is that in this parameter regime one can obtain an analytic characterization of  
the phase portrait of system \eqref{micro:dyn:nonoise:AB}, still displaying a 
rich variety of cases. 
The central concern in the subsequent sections will be the investigation of    
the conditions under which noise-induced periodicity occurs. 
  
\begin{figure}[h!]
    \centering
    \subfigure[\label{fig:PhDiaNoNoiseNoInternalInteractions:first}]{\includegraphics[width=0.325\textwidth]{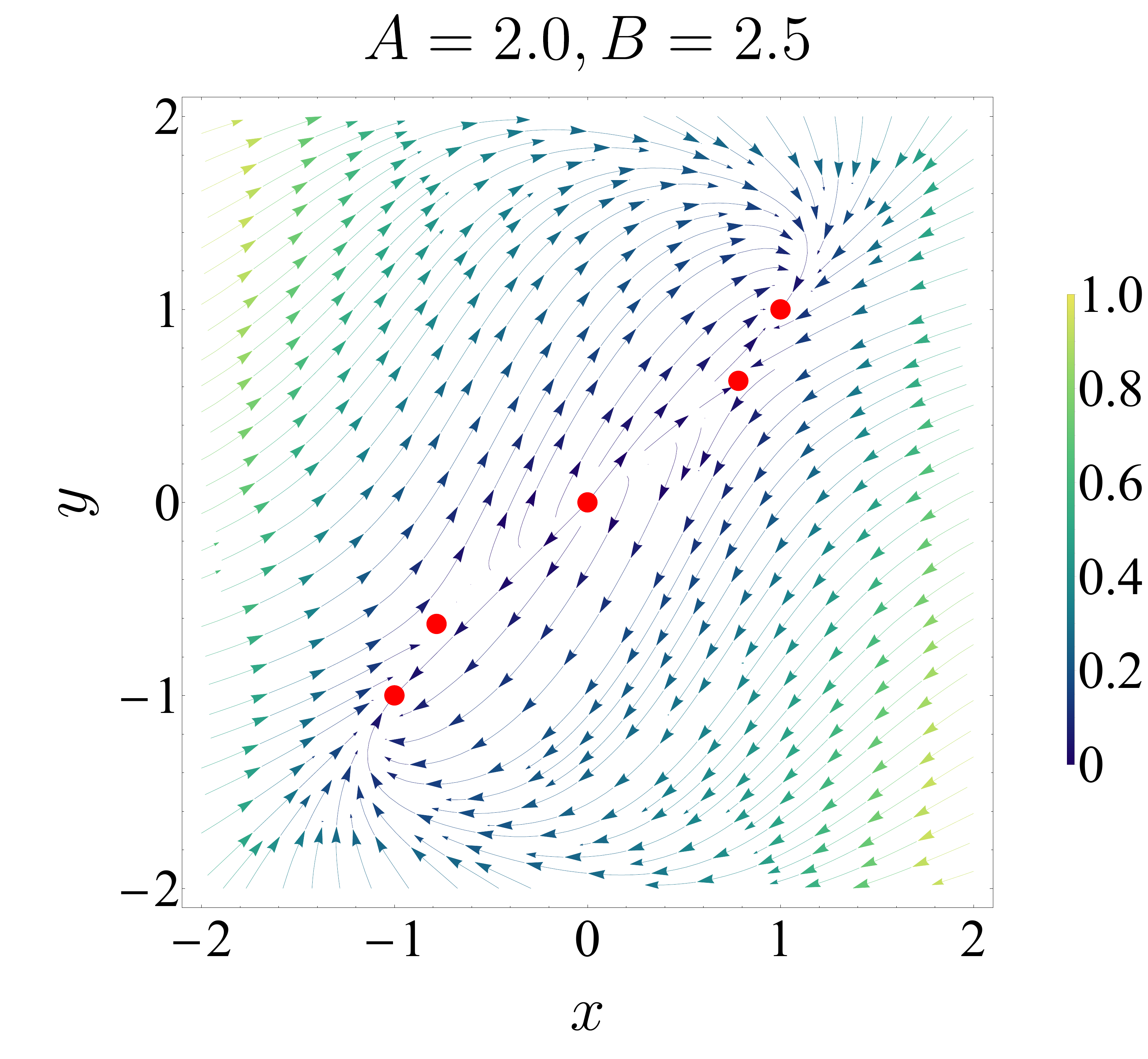}}
    \subfigure[\label{fig:PhDiaNoNoiseNoInternalInteractions:second}]{\includegraphics[width=0.325\textwidth]{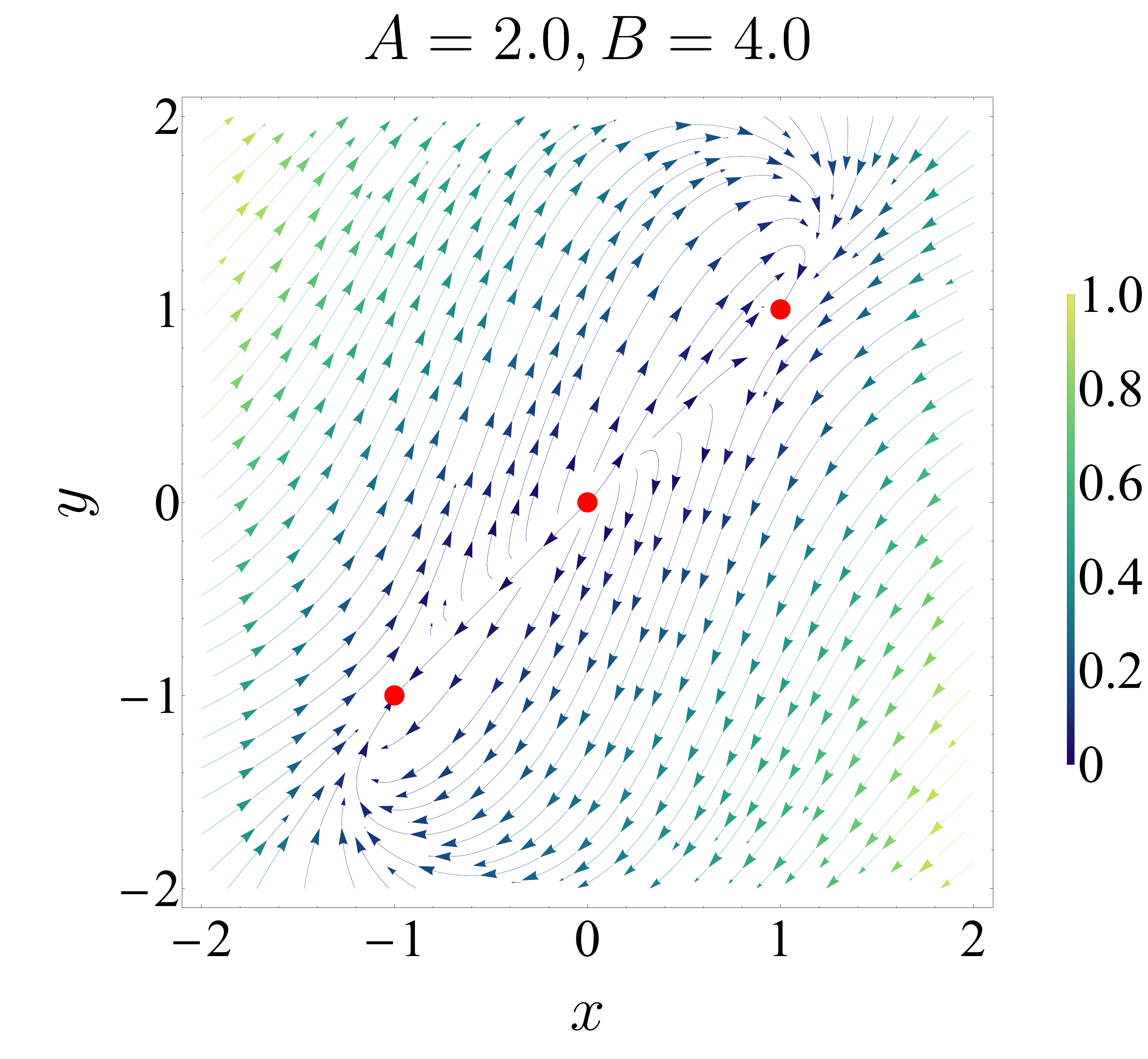}}
    \subfigure[\label{fig:PhDiaNoNoiseNoInternalInteractions:third}]{\includegraphics[width=0.325\textwidth]{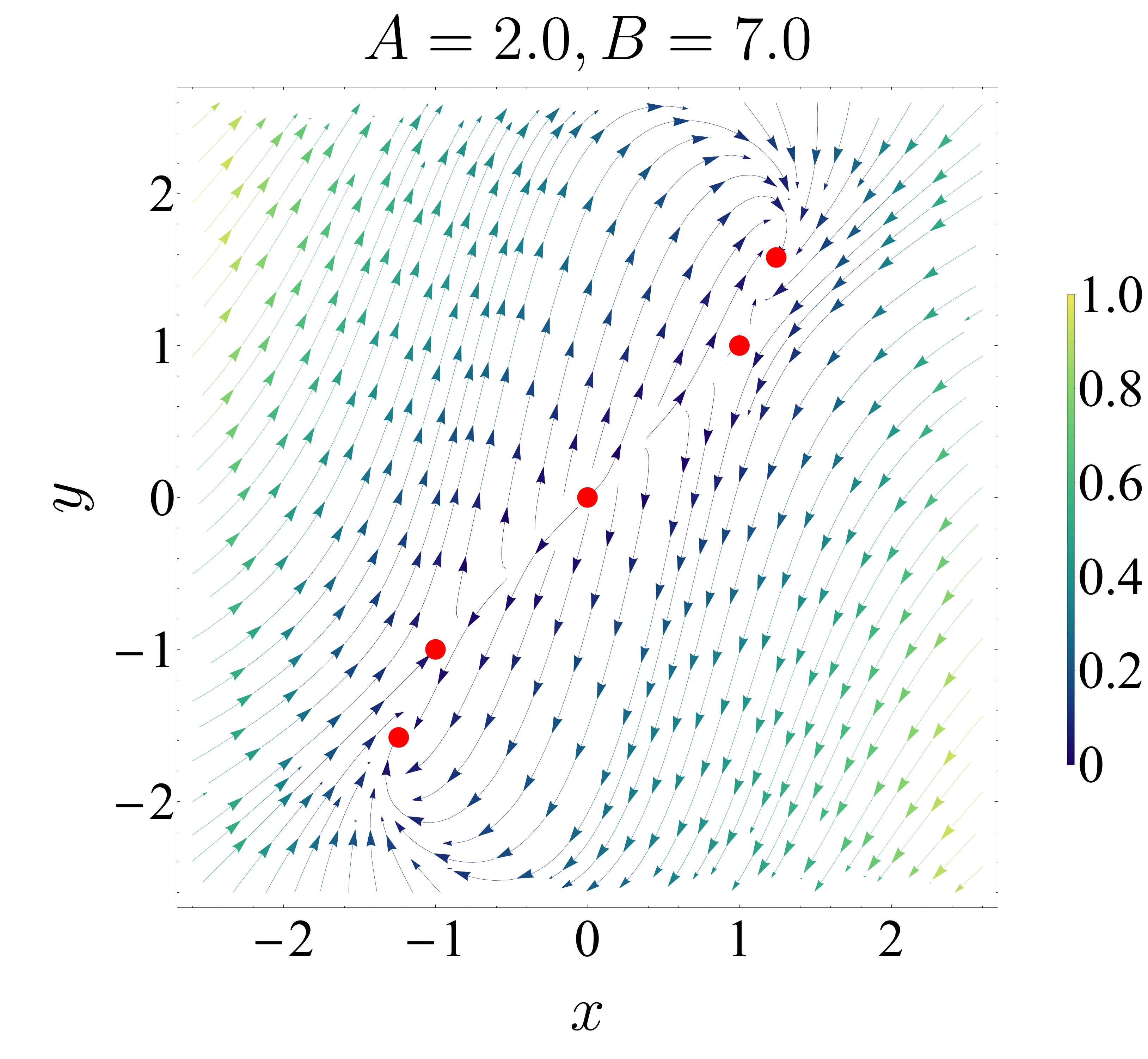}}
    \caption{\footnotesize{Phase portraits of system \eqref{micro:dyn:nonoise:AB} for diverse values of $A$ and $B$. 
\subref{fig:PhDiaNoNoiseNoInternalInteractions:first} Case $A-1<B<A+2$ with $A=2$ and $B=2.5$. Fixed points: $\left(0,0\right)$ is an unstable node, $\pm\left(1,1\right)$ are stable nodes and $\pm\left(0.78,0.63\right)$ (numerically obtained coordinates) are saddle points.
\subref{fig:PhDiaNoNoiseNoInternalInteractions:second} Case $B=A+2$ with $A=2$ and $B=4$. Fixed points: $\left(0,0\right)$ is an unstable node and $\pm\left(1,1\right)$ have a negative and a zero eigenvalue. 
\subref{fig:PhDiaNoNoiseNoInternalInteractions:third} Case $B>A+2$ with $A=2$ and $B=7$. Fixed points: $\left(0,0\right)$ is an unstable node, $\pm\left(1,1\right)$ are saddle points and $\pm\left(1.24, 1.58\right)$ (numerically obtained coordinates) are stable spirals. 
Red dots mark the equilibria. 
Streamline colors correspond to the magnitude of the vector field scaled to $[0,1]$ (relative magnitude).
A detailed analysis of the nature of the fixed points in the three regimes can be found in Appendix~\ref{appendix:C}.}} 
\label{fig:PhDiaNoNoiseNoInternalInteractions}
\end{figure}

To this end, we studied the location and the nature of the fixed points of system \eqref{micro:dyn:nonoise:AB} by varying $A$ and $B$ under the regime given by hypothesis \textbf{(H)} and checked that no local bifurcation generating limit cycles occurs. Unfortunately, the global analysis of the system turns out to be very involved and we are able to exclude the existence of limit cycles only by numerical evidences (see Fig.~\ref{fig:PhDiaNoNoiseNoInternalInteractions}).

System \eqref{micro:dyn:nonoise:AB} admits the following equilibria:

\begin{itemize}
    \item  The fixed points $\left(0,0\right)$ and $\pm\left(1,1\right)$ are present for any value of $A$ and $B$. However, their nature changes depending on the parameters. More specifically,
    \begin{itemize}
        \item when  $A-1<B<A+2$, $\left(0,0\right)$ is an unstable node and  $\pm\left(1,1\right)$ are stable nodes. 
        \item when $B=A+2$, $\left(0,0\right)$ is an unstable node and $\pm \left(1,1\right)$ have a neutral and a stable direction.
        \item for $B>A+2$, $\left(0,0\right)$ is an unstable node and  $\pm\left(1,1\right)$ are saddle points.
    \end{itemize}
\item Depending on the values of $A$ and $B$, there may be two additional equilibria. In particular, three situations may arise:
\begin{itemize}
    \item when $A-1<B<A+2$, there exists $\beta>0$ such that the points $\pm \left(x,\beta x\right)$ are fixed points for \eqref{micro:dyn:nonoise:AB}, with $0<x<1$ and $\beta <1$. That is, the equilibria are $\left(0,0\right)$, $\pm\left(1,1\right)$ and $\pm \left(x,\beta x\right)$, symmetrically located in the first and the third quadrants. The fixed points $\pm \left(x,\beta x\right)$ are saddle points.
    \item when $B=A+2$, no other fixed points are present apart from $\left(0,0\right)$ and $\pm \left(1,1\right)$.
    \item when $B>A+2$, there exists $\beta>0$ such that $\pm \left(x,\beta x\right)$ are  fixed points for \eqref{micro:dyn:nonoise:AB}, with $x>1$ and $\beta >1 $. That is, system \eqref{micro:dyn:nonoise:AB} has five equilibria: $\left(0,0\right)$, $\pm\left(1,1\right)$ and $\pm \left(x,\beta x\right)$, symmetrically located in the first and the third quadrants. The fixed points $\pm \left(x,\beta x\right)$ are stable nodes.
\end{itemize}
\end{itemize}

\begin{table}
    \centering
    \begin{tabular}{|c||c|c|c|}
	\hline
 	\multicolumn{4}{|c|}{\cellcolor{Gray}\footnotesize{Parameter regime $A>1$, $B>A-1$}} \\
 	\hline\hline
    \cellcolor{LightGray} \diagbox[width=7em]{ \footnotesize{Parameters}}{\footnotesize{Equilibria}} & \cellcolor{LightGray} \footnotesize{$\left(0,0\right)$} & \cellcolor{LightGray} \footnotesize{$\pm\left(1,1\right)$} & \cellcolor{LightGray} \footnotesize{$\pm\left(x, \beta x\right)$} \\
    \hline\hline
    \cellcolor{LightGray} \footnotesize{${A-1 < B < A+2}$} &  \footnotesize{unstable node} &  \footnotesize{stable nodes} &  \footnotesize{$0<x<1$, $0<\beta<1$,} \\ 
    \cellcolor{LightGray} &&&  \footnotesize{saddle points} \\[.3cm]
    \cellcolor{LightGray} \footnotesize{$B = A+2$} &  \footnotesize{unstable node} &  \footnotesize{one negative and one null eigenvalue} & $-$ \\[.3cm]
    \cellcolor{LightGray} \footnotesize{$B > A+2$} &  \footnotesize{unstable node} &  \footnotesize{saddle points} &  \footnotesize{$x>1$, $\beta >1$, }\\
     \cellcolor{LightGray} &&&  \footnotesize{stable nodes}\\
    \hline
    \end{tabular}
    \caption{\footnotesize{Regime $A>1$, $B> A-1$. Nature of the fixed points of system \eqref{micro:dyn:nonoise:AB} for different values of the parameters.}}
    \label{tab:eq_points}
\end{table}

The depicted scenarios are summarized in Table~\ref{tab:eq_points}. We refer the reader to Appendix~\ref{appendix:C} for a detailed proof. 
In Fig.~\ref{fig:PhDiaNoNoiseNoInternalInteractions}, we display numerically obtained phase portraits for specific values of the parameters in the three cases $A-1<B<A+2$, $B=A+2$ and $B>A+2$. In all these cases, numerical investigations strongly corroborate the absence of limit cycles for system \eqref{micro:dyn:nonoise:AB}. 

We remark that the main results of this paper, given in Sections~\ref{sec:prop:chaos} and \ref{sec:gauss:approx}, hold for all $A,\, B >0$, as one can see from the proofs in the Appendices. Furthermore, qualitatively analogous behaviors were numerically observed in the case $0<A\leq 1$, $B>0$, when extra fixed points for system \eqref{micro:dyn:nonoise:AB} may exist.

\subsection{The Fokker-Planck equation}\label{sec:FP}
\begin{figure}
    \centering
    \subfigure{\includegraphics[width=0.49\textwidth]{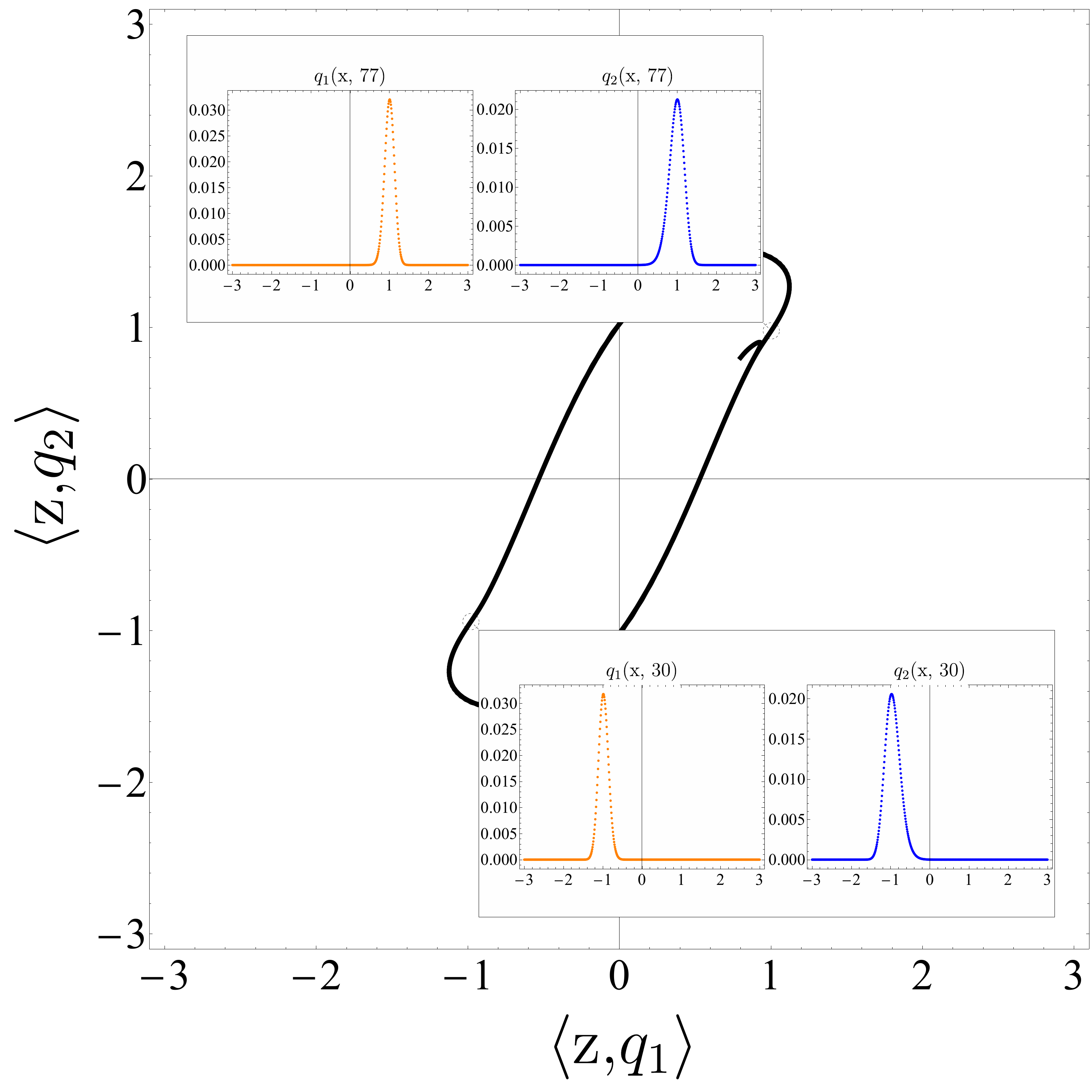}}
    \caption{\footnotesize{Temporal evolution of the average positions $\langle z, q_1\rangle$ and $\langle z, q_2\rangle$ of the two populations in the thermodynamic limit. Parameter values: $A=2$ and $B=2.5$; the other regimes are analogous.  The insets show the densities $q_1$ (orange) and $q_2$ (blue) at some times during the simulation.}}
    \label{fig:FP_simulations}
\end{figure}

The long-time behavior of the law of the solution to system \eqref{prop:chaos}  may be investigated by considering the corresponding Fokker-Planck equation, that reads as 
\begin{align}\label{fokker:planck}
\frac{\partial q_1}{\partial t} &=\frac{\sigma^2}{2} \, \frac{\partial^2 q_1}{\partial z^2}-\frac{\partial}{\partial z}\left\{\left[(1-\alpha\theta_{11}- \left(1-\alpha\right)\theta_{12})z-z^3\right]q_1\right\}\nonumber\\
&-\alpha \theta_{11} \langle z,q_1\rangle\frac{\partial q_1}{\partial z}-  \left(1-\alpha\right)\theta_{12}\langle z,q_2\rangle\frac{\partial q_1}{\partial z}\nonumber\\[.3cm]
\frac{\partial q_2}{\partial t} &=\frac{\sigma^2}{2} \, \frac{\partial^2 q_2}{\partial z^2}-\frac{\partial}{\partial z}\left\{\left[(1-\alpha \theta_{21}-(1-\alpha)\theta_{22})z-z^3\right]q_2\right\}\nonumber\\
&-\alpha \theta_{21}\langle z,q_1\rangle\frac{\partial q_2}{\partial z}- (1-\alpha)\theta_{22}\langle z,q_2\rangle\frac{\partial q_2}{\partial z},
\end{align}
where time and space dependencies have been left implicit for simplicity of notation. Here $\langle z,q_i \rangle := \int z q_i(z;t) dz$, with $i=1,2$. The regularizing effect of the second-order partial derivatives guarantees that, for $t \in [0,T]$, the laws of $x(t)$ and $y(t)$ have respective densities $q_1(\cdot;t)$ and $q_2(\cdot;t)$ solving \eqref{fokker:planck}.
By using the finite element method \cite{quarteroni}, we performed numerical simulations of system \eqref{fokker:planck} starting from the initial distributions $q_1(z;0)=q_2(z;0)=\delta_{0.8}(z)$. These initial conditions correspond to what we did in Section~\ref{sec:numerical:simulations}, where we initialized the particles of both groups at $z=0.8$ in the simulations of the microscopic system. We observed that $q_1$ and $q_2$ both assume a bell shape
during the simulation, while the average positions of the two populations, $\langle z, q_i\rangle$ ($i=1,2$), computed numerically, display an oscillatory behavior.  
We show the results of these simulations in Fig.~\ref{fig:FP_simulations}. 
The above considerations justify the idea of the Gaussian approximation for system \eqref{prop:chaos} that will be analyzed in the following section.

\subsection{Small-noise approximation}\label{sec:gauss:approx} 
In this section we derive a small-noise approximation of system \eqref{prop:chaos}. In particular, motivated by what we observed in Section~\ref{sec:FP}, we build a pair of independent Gaussian processes $\left(\left(\tilde{x}(t), \tilde{y}(t)\right), 0 \leq t \leq T \right)$ that closely follows $((x(t),y(t)), 0 \leq t \leq T)$, solution to \eqref{prop:chaos}, when the noise is small. Although such an approximation holds rigorously true in the limit of vanishing noise, numerical simulations suggest it remains valid also beyond the assumption  $\sigma \ll 1$ and that it explains the qualitative behavior of system \eqref{micro:dyn} shown in Section~\ref{sec:numerical:simulations}. We give the precise statement of our result below, whereas the proof is postponed to Appendix~\ref{appendix:D}. Here we remark that it is possible to take $(\tilde{x}(t), 0 \leq t \leq T)$ independent of $(\tilde{y}(t), 0 \leq t \leq T)$ because of the specific form of the equations in \eqref{prop:chaos}, that do not have mixed terms (i.e. of the type $x^n\,y^m$).

The first step towards the 
Gaussian approximation of \eqref{prop:chaos} is the derivation of the equations of the moments of $x(t)$ and $y(t)$ in system \eqref{prop:chaos}. Since the approximation will be given by a pair of independent processes, we can avoid computing mixed moments (see Appendix~\ref{appendix:D}).
By applying It{\^o}'s rule to system  \eqref{prop:chaos}, we can obtain the SDEs solved by $x^p(t)$ and $y^p(t)$ for any $p \geq 1$. This yields
\begin{align}\label{previous}
dx^{p} &= \sigma p x^{p-1} dw_1 + \left[ -p x^{p+2} + p x^p -\alpha\theta_{11} p \left(x-\mathbb{E}[x]\right)x^{p-1} \right.\nonumber\\
&\left.- \left(1-\alpha\right)\theta_{12} p \left(x-\mathbb{E}[y]\right)x^{p-1} + \tfrac{\sigma^2}{2}p(p-1)x^{p-2} \right]dt \nonumber \\[.2cm]
dy^{p} &= \sigma p y^{p-1} dw_2 + \left[ -p y^{p+2} + p y^p  -\alpha\theta_{21} p \left(y-\mathbb{E}[x]\right)y^{p-1}\right.\nonumber \\
&\left. -\left(1-\alpha\right)\theta_{22} p \left(y-\mathbb{E}[y]\right) y^{p-1}+ \tfrac{\sigma^2}{2}p (p-1) y^{p-2}\right]dt.
\end{align}
Let $m_{p}^{x}(t)=\mathbb{E}[x^p(t)]$ and $m_{p}^{y}(t)=\mathbb{E}[y^p(t)]$ be the $p$-th moments of the variables $x(t)$ and $y(t)$ solving system \eqref{prop:chaos}, respectively. Taking the expectation in \eqref{previous},
we obtain
\begin{align}\label{eq:moment}
\frac{dm_{p}^{x}}{dt} &= -p m_{p+2}^{x}+p m_{p}^{x}-\alpha \theta_{11}p\left(m_{p}^{x}-m_{1}^{x} \, m_{p-1}^{x}\right) \nonumber\\
&-\left(1-\alpha\right)\theta_{12}p\left(m_{p}^{x}-m_{1}^{y} \, m_{p-1}^{x}\right)+\tfrac{\sigma^2}{2}p(p-1)m_{p-2}^{x}\nonumber\\[.2cm]
\frac{dm_{p}^{y}}{dt} &= -p m_{p+2}^{y}+p m_{p}^{y}-\alpha\theta_{21} p\left(m_{p}^{y}-m_{1}^{x} \, m_{p-1}^{y}\right) \nonumber\\
&-(1-\alpha)\theta_{22}p\left(m_{p}^{y}-m_{1}^{y}\, m_{p-1}^{y}\right)+\tfrac{\sigma^2}{2}p(p-1)m_{p-2}^{y}.
\end{align}
Since the $p$-th moments in \eqref{eq:moment} depend on the $(p+2)$-th moments, the system is infinite dimensional - and hence hardly tractable - unless higher-order moments of $x(t)$ and $y(t)$ are functions of the first moments. The latter would be the case if $x(t)$ and $y(t)$ were Gaussian processes. In general, the processes $x(t)$ and $y(t)$ are neither Gaussian nor independent, however we prove in Appendix~\ref{appendix:D} that it is possible to build a Gauss-Markov process $\left(\left(\tilde{x}(t), \tilde{y}(t)\right), 0 \leq t \leq T \right)$, with independent components, which stays close to $\left(\left(x(t), y(t)\right), 0 \leq t \leq T \right)$ when the noise size is small. We have the following theorem.

\begin{thm}\label{thm:Gaussian_approx}
Fix $T>0$. Let $\left(\left(x(t), y(t)\right), 0 \leq t \leq T \right)$ solve Eq.~\eqref{prop:chaos} with deterministic initial conditions $x(0) = x_0$ and $y(0) = y_0$. There exists a Gaussian Markov process $\left(\left(\tilde{x}(t), \tilde{y}(t)\right), 0 \leq t \leq T \right)$ with $\tilde{x}(0)=x_0$ and $\tilde{y}(0)=y_0$ satisfying the properties:
\begin{enumerate}\label{enum:Gaussian_approx}
\item The first two moments of $\tilde x(t)$ and $\tilde y(t)$ satisfy the respective equations in \eqref{eq:moment} for $p=1,2$.
\item For all $T > 0$, there exists a constant $C_T>0$ such that, for every $\sigma>0$, it holds
  \begin{equation*}
        \mathbb{E}\left[\sup_{t\in[0,T]}\left\{ \left| x(t) -\tilde x(t)\right|+\left| y(t) -\tilde y(t)\right|\right\}\right] \leq C_T \sigma^2.
    \end{equation*}
This means that the processes $\left(\tilde{x}(t), 0 \leq t \leq T \right)$ and $\left(\tilde{y}(t), 0 \leq t \leq T\right)$ are simultaneously $\sigma$-closed to the solutions of \eqref{prop:chaos}.
\end{enumerate}
\end{thm}

Since $\tilde x(t)$ and $\tilde y(t)$ are Gaussian, their higher-order moments are polynomial functions of the first two moments. In particular, the laws of $\tilde x(t)$ and $\tilde y(t)$ are completely determined by the dynamics of the respective mean and variance.  Thus, rather than studying the infinite dimensional system \eqref{eq:moment}, it suffices to analyze the subsystem describing the time evolution of the mean and the variance of each approximating process. 
We will show in Appendix~\ref{appendix:D} that such a system is
    \begin{align}\label{eq:odes}
         \frac{dm_1}{dt} &= -m_1^3 + m_1 (1-3 v_1) -A(m_1-m_2)\nonumber\\
         \frac{dm_2}{dt} &=-m_2^3 + m_2  (1-3 v_2)+ B\left(m_2 - m_1\right)\nonumber\\
         \frac{dv_1}{dt} &= -6 v_1^2 -6 m_1^2 v_1 + 2 v_1 - 2\alpha \theta_{11} v_1 -2 A v_1+ \sigma^2\nonumber\\
         \frac{dv_2}{dt} &= -6 v_2^2 -6 m_2^2 v_2 + 2 v_2 + 2 B v_2 - 2\left(1-\alpha\right)\theta_{22} v_2 + \sigma^2, 
    \end{align}
where $m_1(t)$ (resp. $m_2(t)$) is the expectation of $\tilde{x}(t)$ (resp. $\tilde{y}(t)$) and $v_1(t)$ (resp. $v_2(t)$) is the variance of $\tilde{x}(t)$ (resp. $\tilde{y}(t)$). As before, we have set $A\coloneqq \left(1-\alpha\right)\theta_{12}$ and \mbox{$B\coloneqq -\alpha\theta_{21}$}. 

For the values of $\theta_{11}$, $\theta_{22}$, $A$ and $B$ considered in this paper (i.e., $\theta_{11}=\theta_{22}=8$ and $A$ and $B$ as reported in Table~\ref{tab:comparison_fft_poincare}), the dynamical system \eqref{eq:odes} features a \emph{subcritical Hopf bifurcation} \cite{Per2001} at the equilibrium $(m_1,m_2,v_1,v_2)=\left(0,0,\tilde v_1,\tilde v_2\right)$, for a critical value $\sigma_c=\sigma_c(\theta_{11}, \theta_{22}, A,B)$ of the noise size, as reported in Appendix~\ref{Hopf}. In other words, when the noise intensity \emph{decreases} to cross the threshold value $\sigma_c$, the fixed point $\left(0,0,\tilde v_1,\tilde v_2\right)$ changes its nature from stable to unstable and, at the same time, a stable limit cycle appears. Thus, in an intermediate range of noise size, system \eqref{eq:odes} displays stable rhythmic oscillations that disappear for $\sigma=0$. Indeed, when $\sigma=0$, $v_1=v_2=0$ is a fixed point of the subsystem formed by the third and fourth equations in \eqref{eq:odes}. As a consequence, the zero-noise limit of  the first two equations in \eqref{eq:odes} reduces to the noiseless version of system \eqref{prop:chaos}, which does not display any oscillatory behavior.

\begin{figure}[h!]
    \centering
    \begin{center}
    \subfigure{\includegraphics[width=0.45\textwidth]{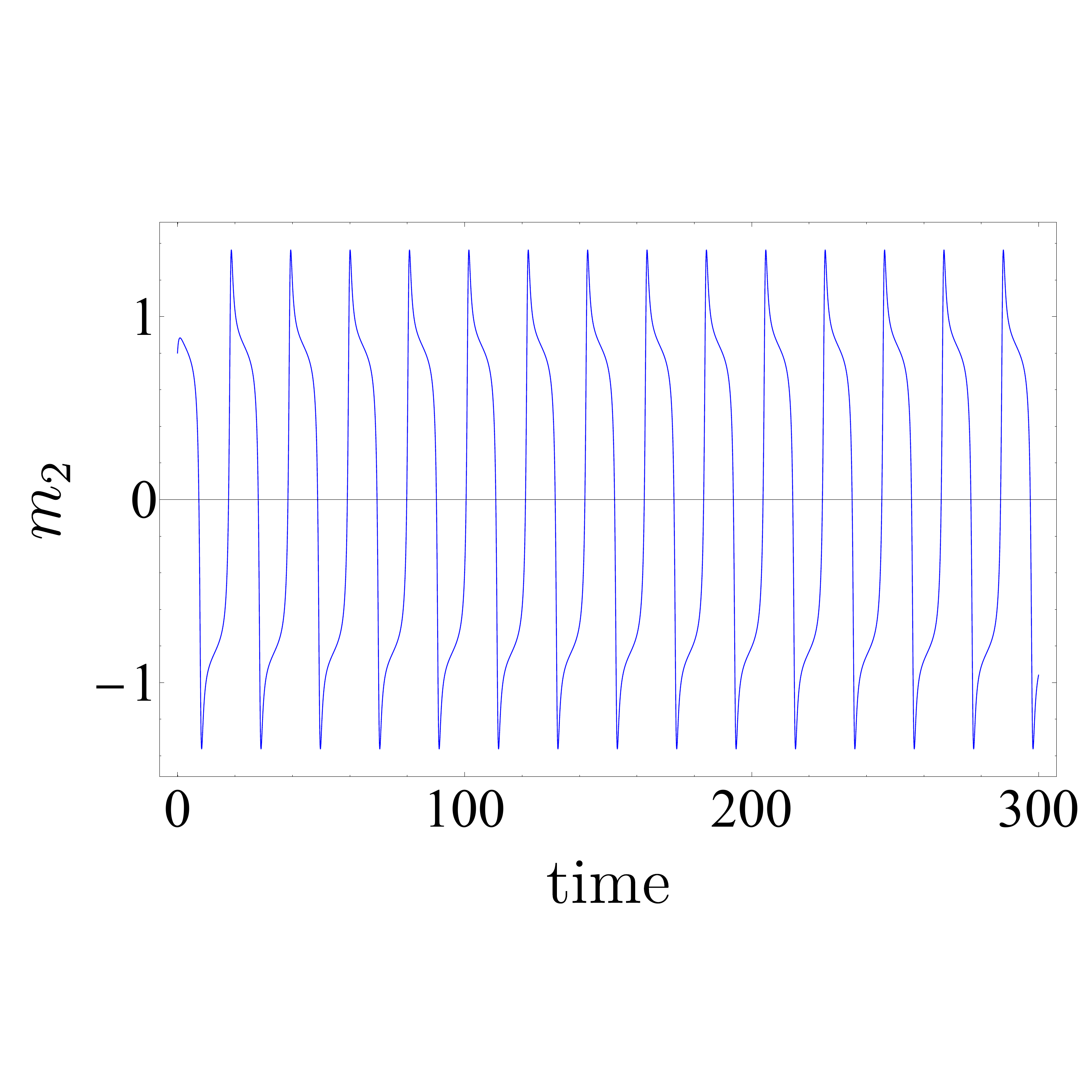}}
    \subfigure{\includegraphics[width=0.45\textwidth]{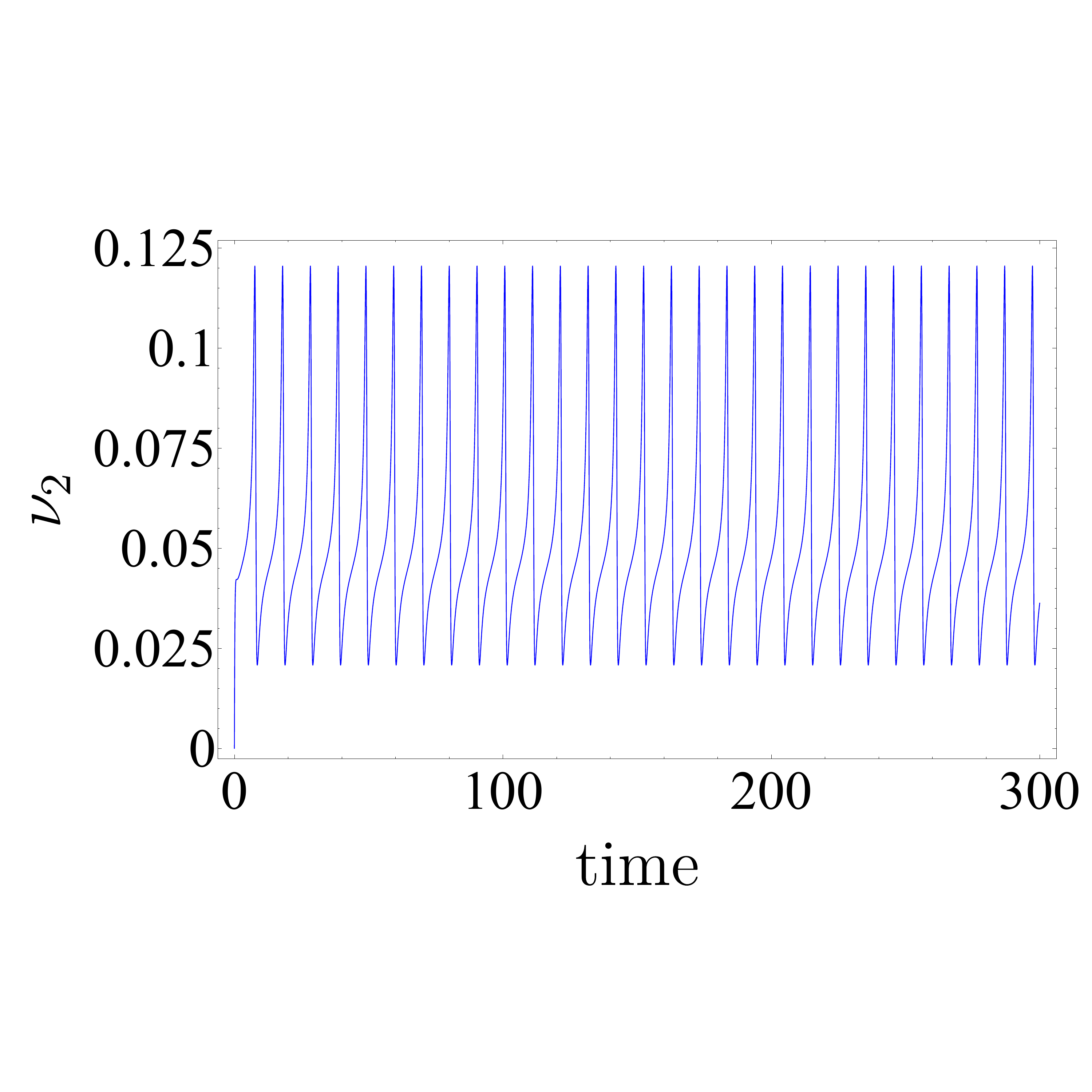}}
    \end{center}

    \vspace{-50pt}
    \begin{center}
    \subfigure{\includegraphics[width=0.45\textwidth]{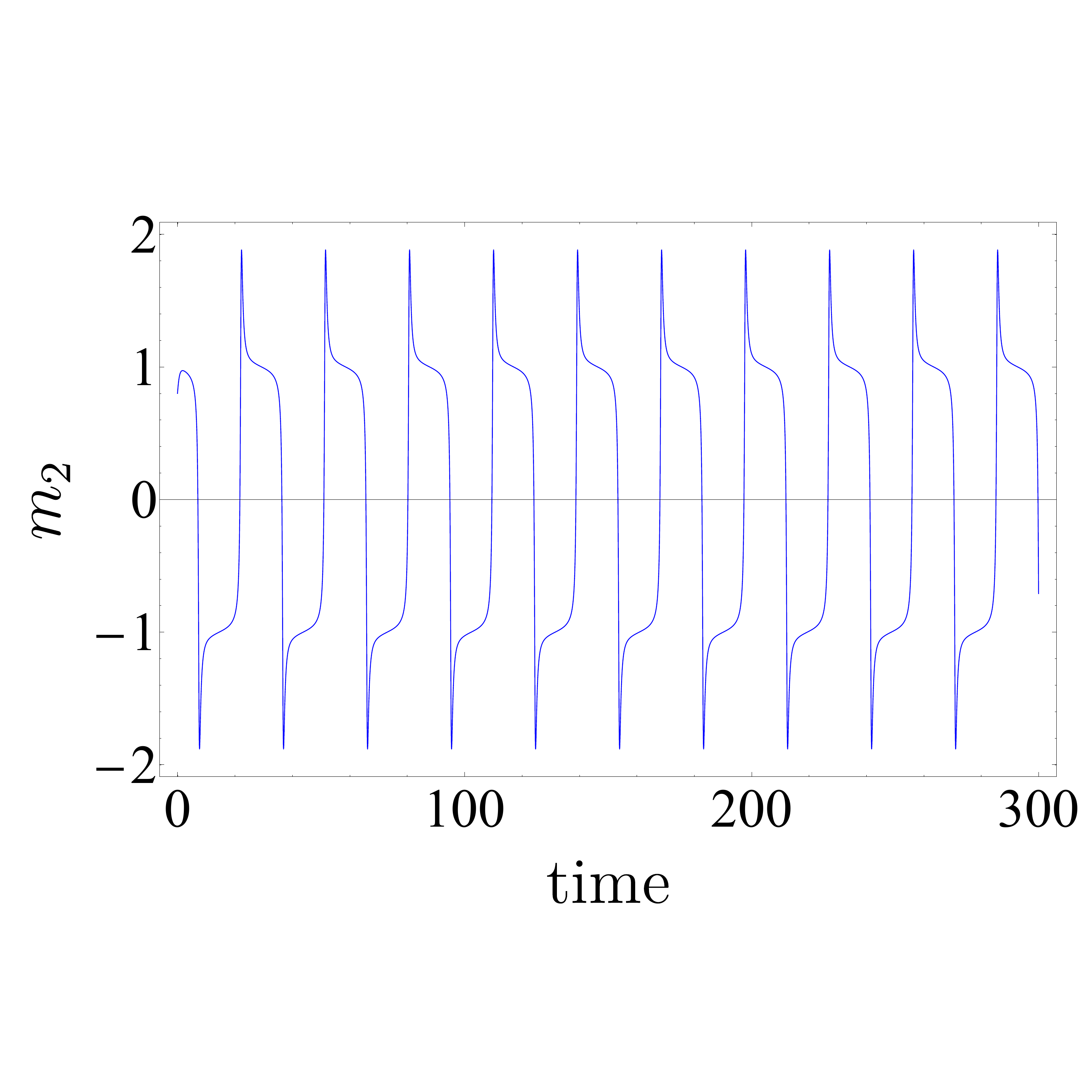}}
    \subfigure{\includegraphics[width=0.45\textwidth]{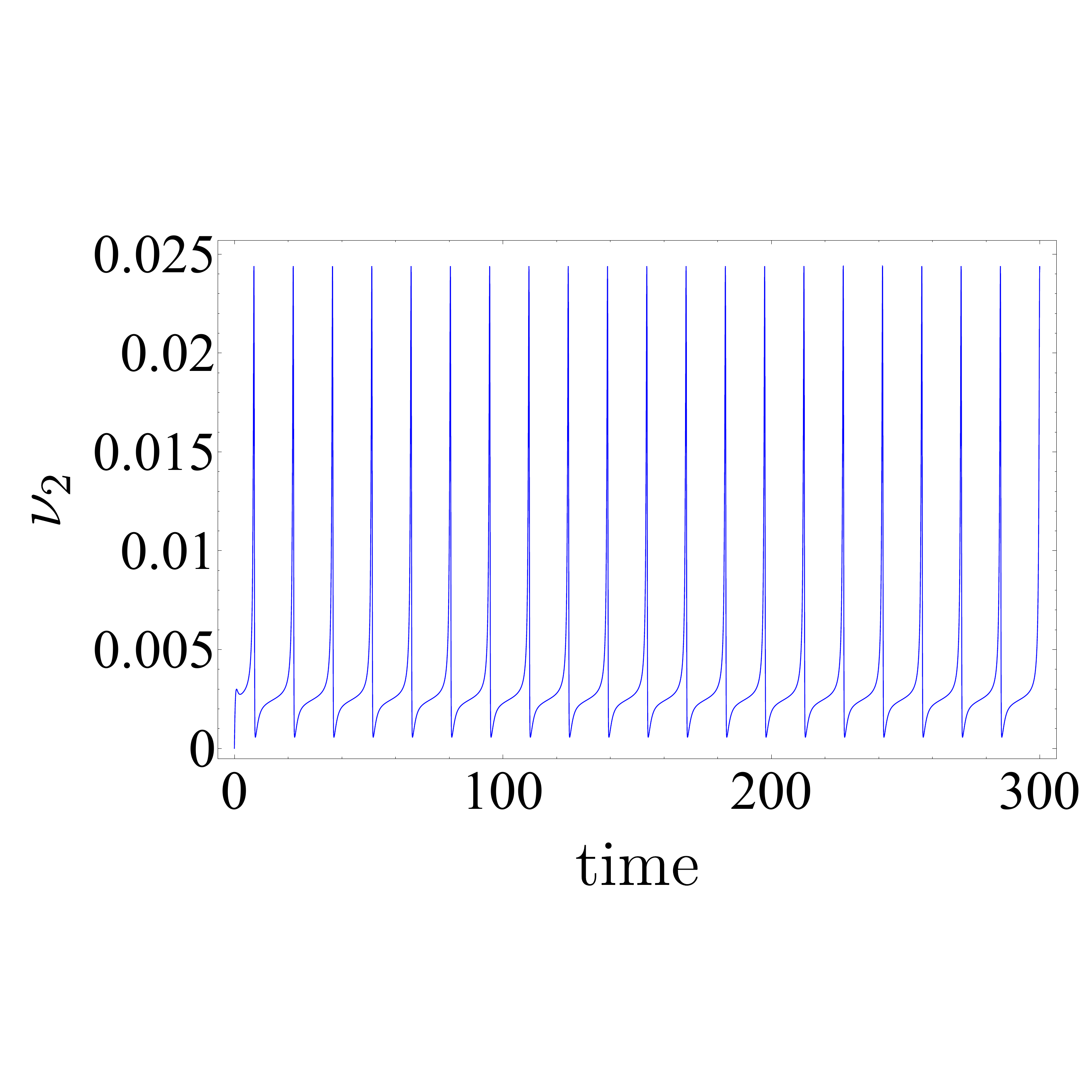}}
    \end{center}

    \vspace{-50pt}
    \begin{center}
    \subfigure{\includegraphics[width=0.45\textwidth]{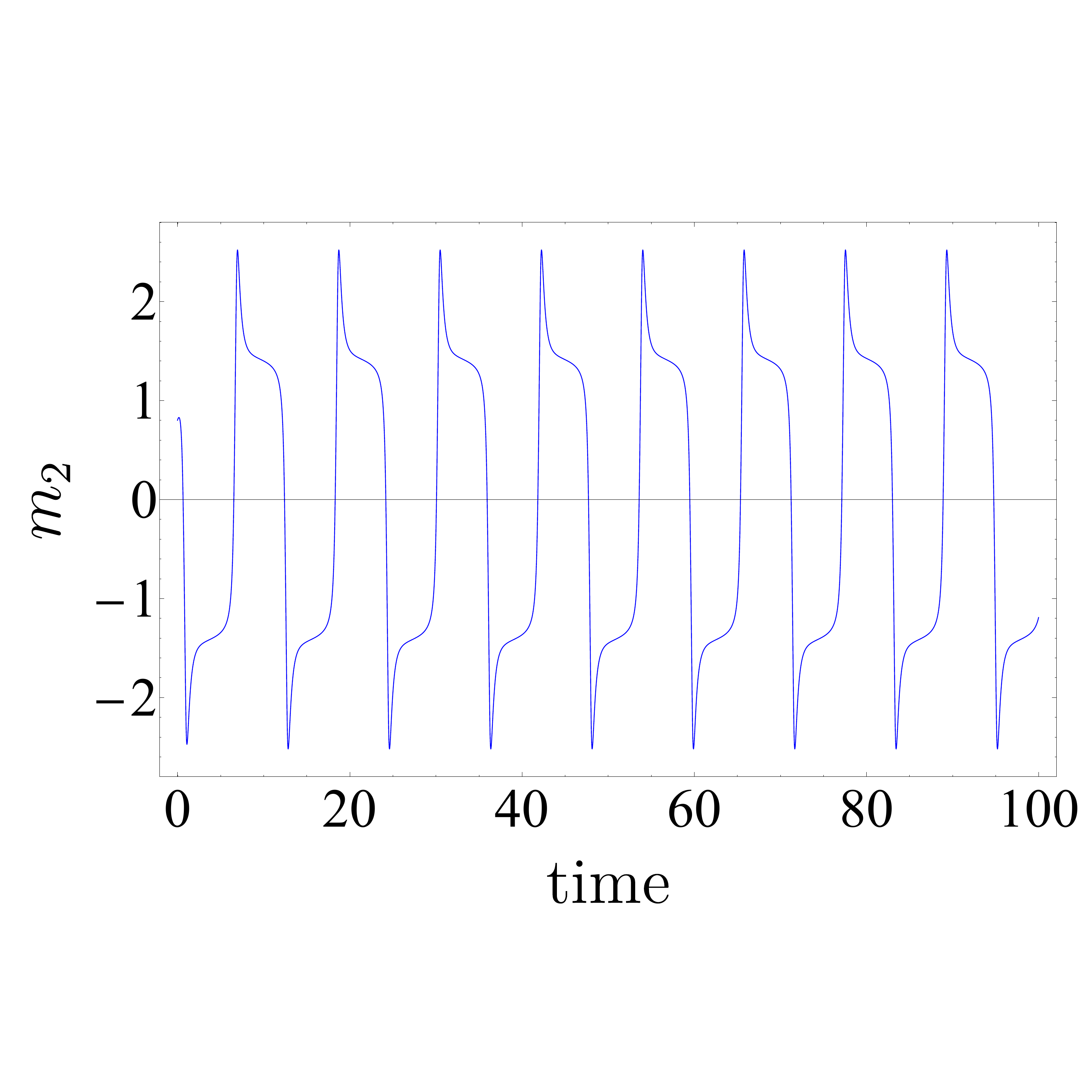}}
    \subfigure{\includegraphics[width=0.45\textwidth]{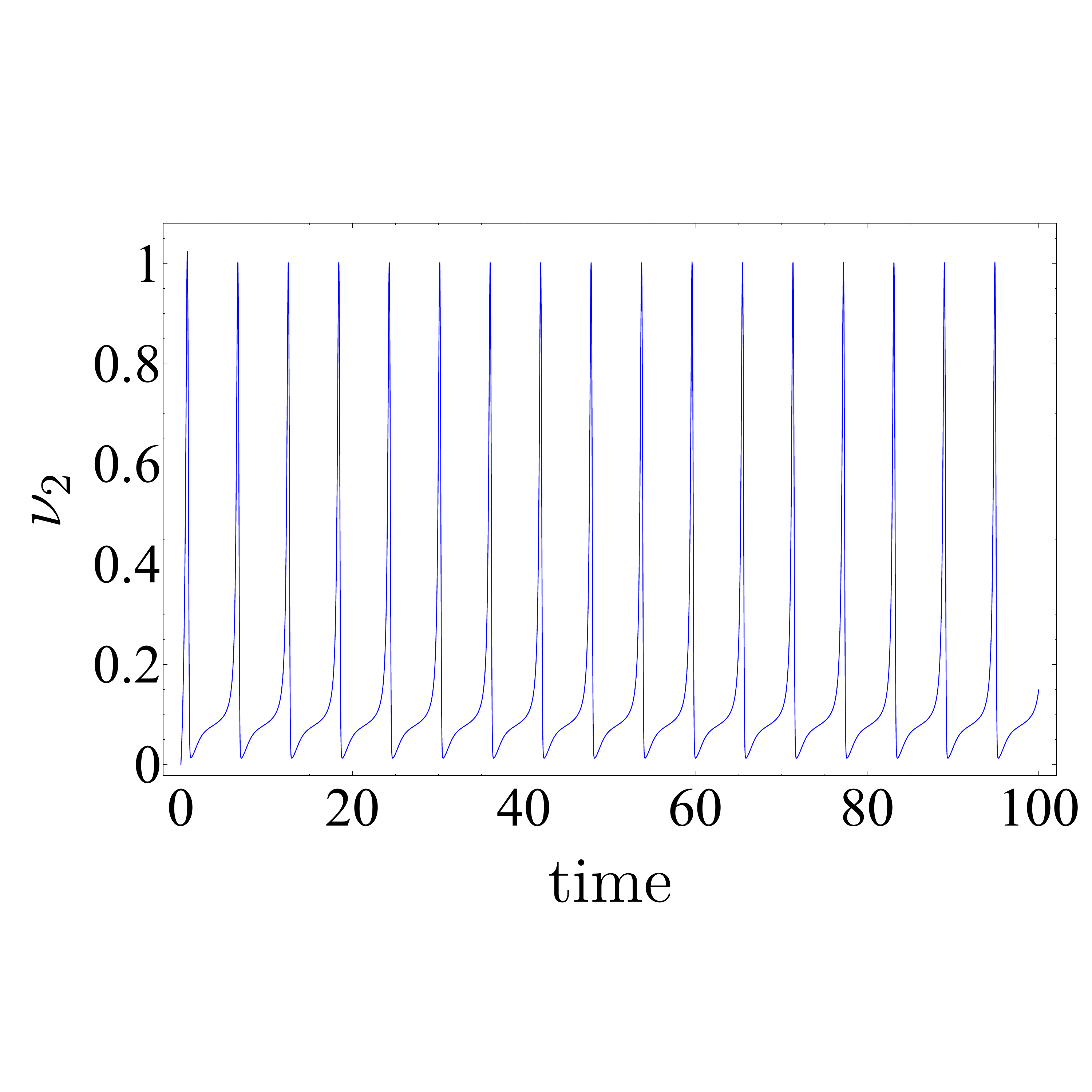}}
    \end{center}
    
    \vspace{-25pt}
    \caption{\footnotesize{Time evolution of the mean $m_2$ and the variance $v_2$ according to the dynamical system \eqref{eq:odes}. In all cases, we considered $10^6$ iterations with a time-step $dt=0.005$, $\alpha = 0.5$, $\theta_{11}=\theta_{22}=8$. From top to bottom: $A-1<B<A+2$, in particular, $A=2$ and $B=2.5$, with $\sigma=0.5$; $B=A+2$, in particular, $A=2$ and $B=4$, with $\sigma=0.1$; $B>A+2$, in particular, $A=2$ and $B=7$, with $\sigma=0.6$.}}
\label{fig:Gauss:m2:v2}
\end{figure}

\begin{figure}[h!]
    \centering
    \begin{center}
    \subfigure
{\includegraphics[width=0.325\textwidth]{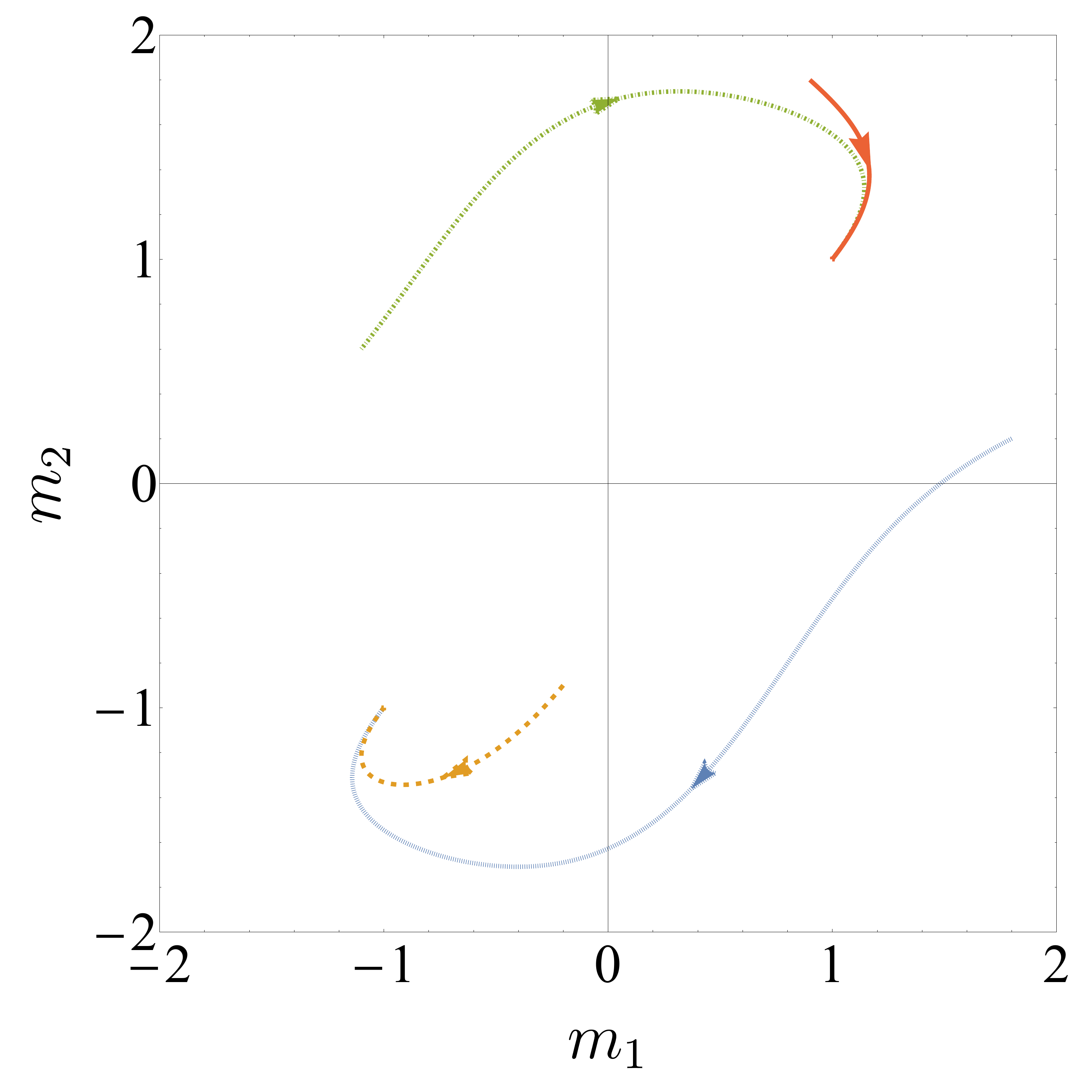}}
    \subfigure
{\includegraphics[width=0.325\textwidth]{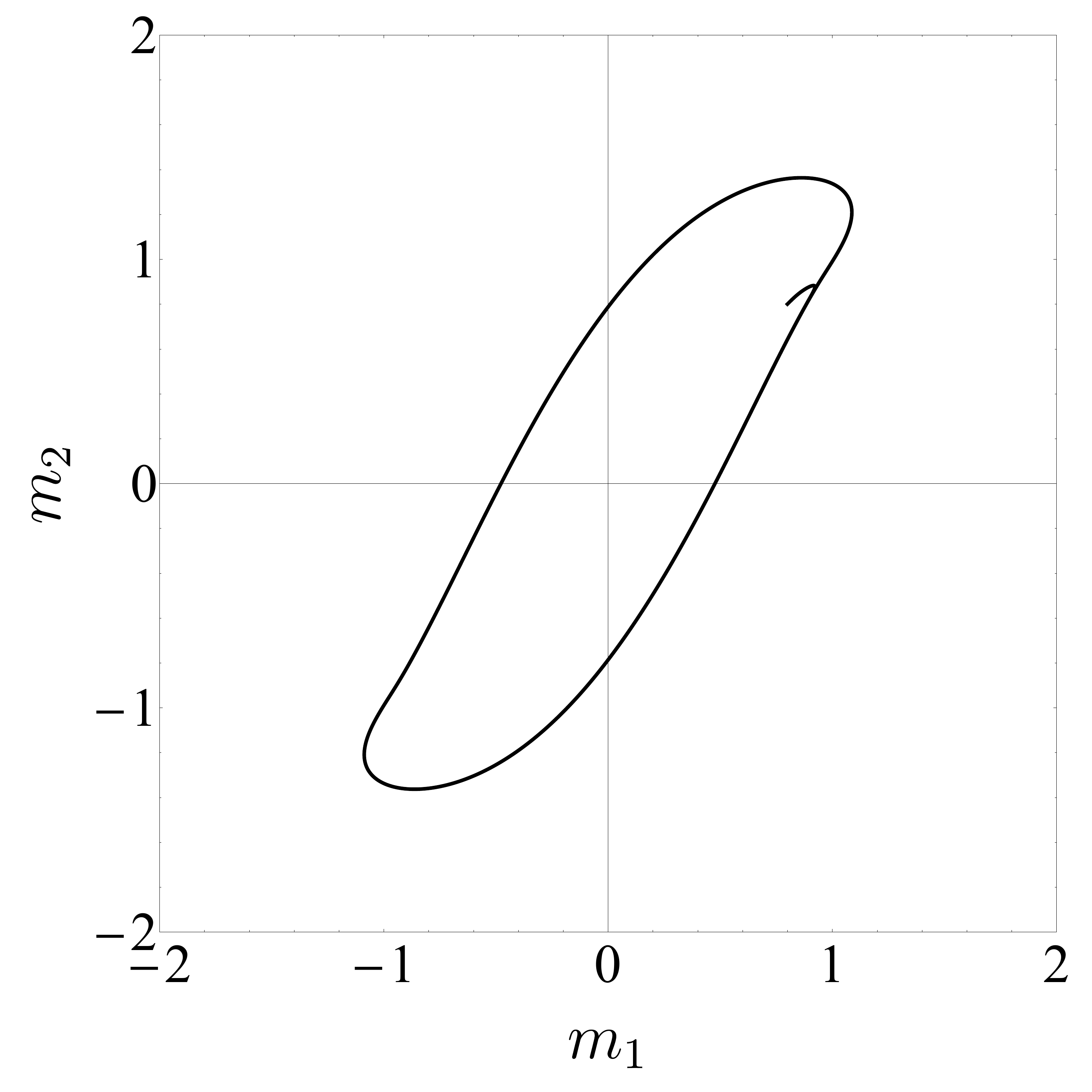}}
    \subfigure
{\includegraphics[width=0.325\textwidth]{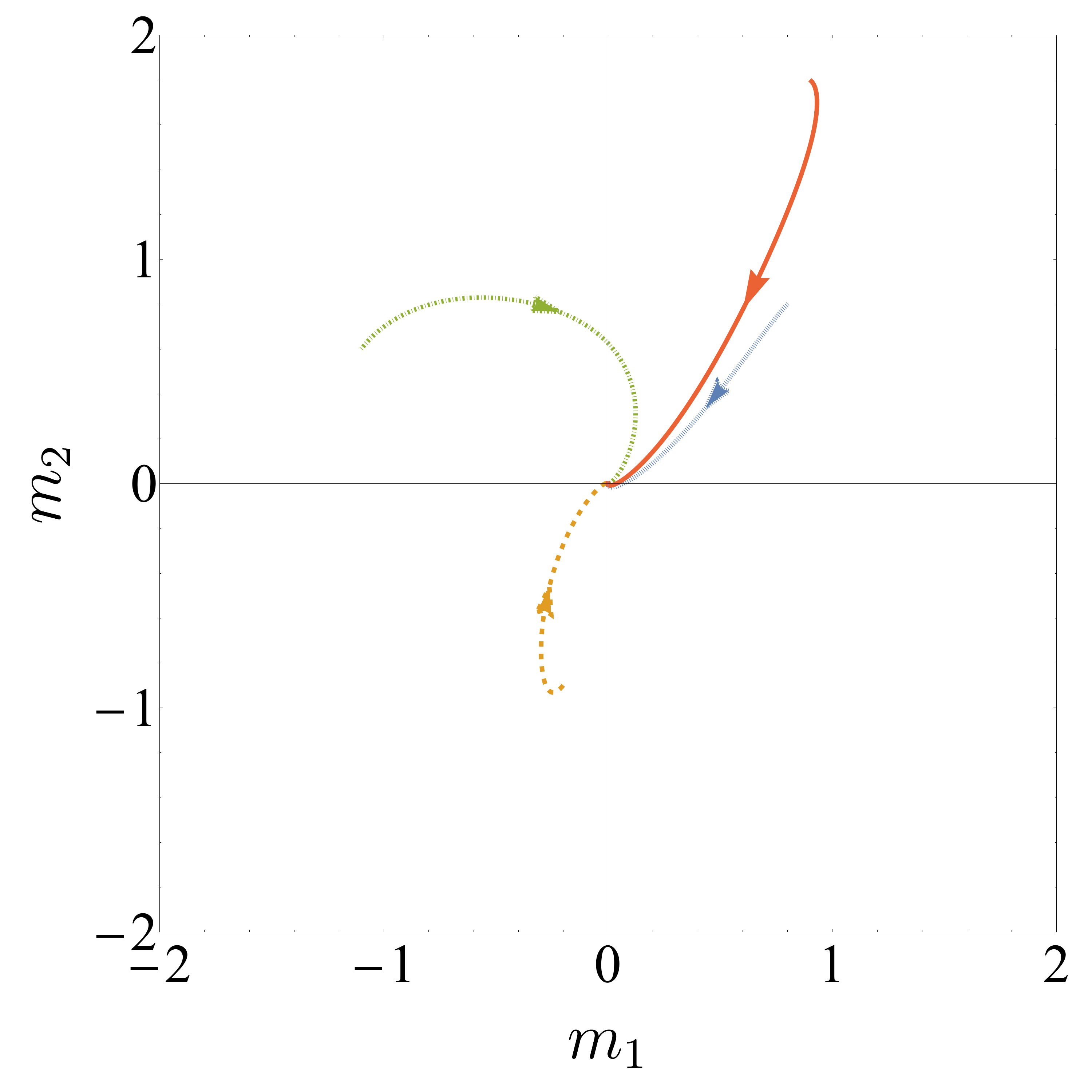}}
    \end{center}

    \vspace{5pt}
    \begin{center}
    \subfigure
{\includegraphics[width=0.325\textwidth]{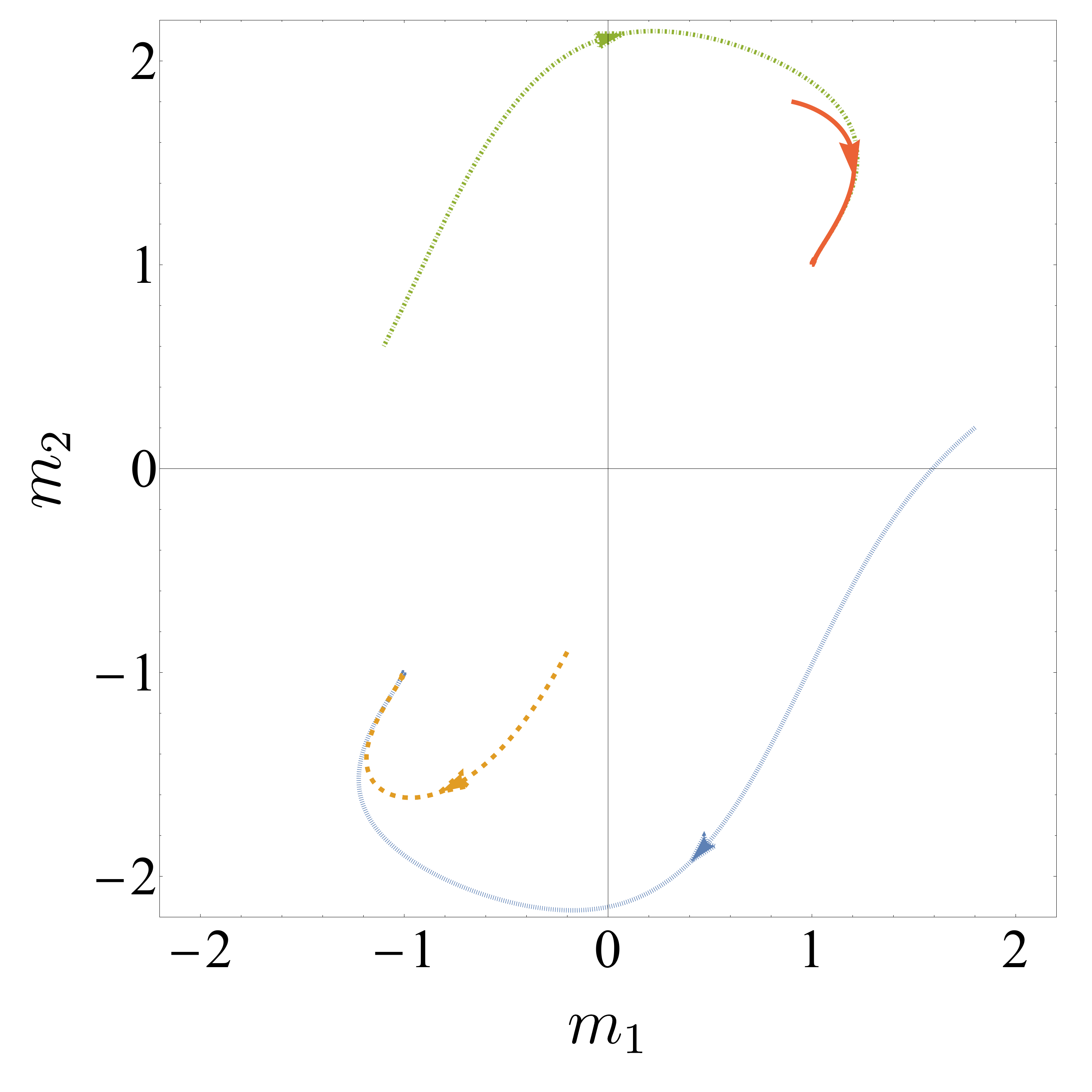}}
    \subfigure
{\includegraphics[width=0.325\textwidth]{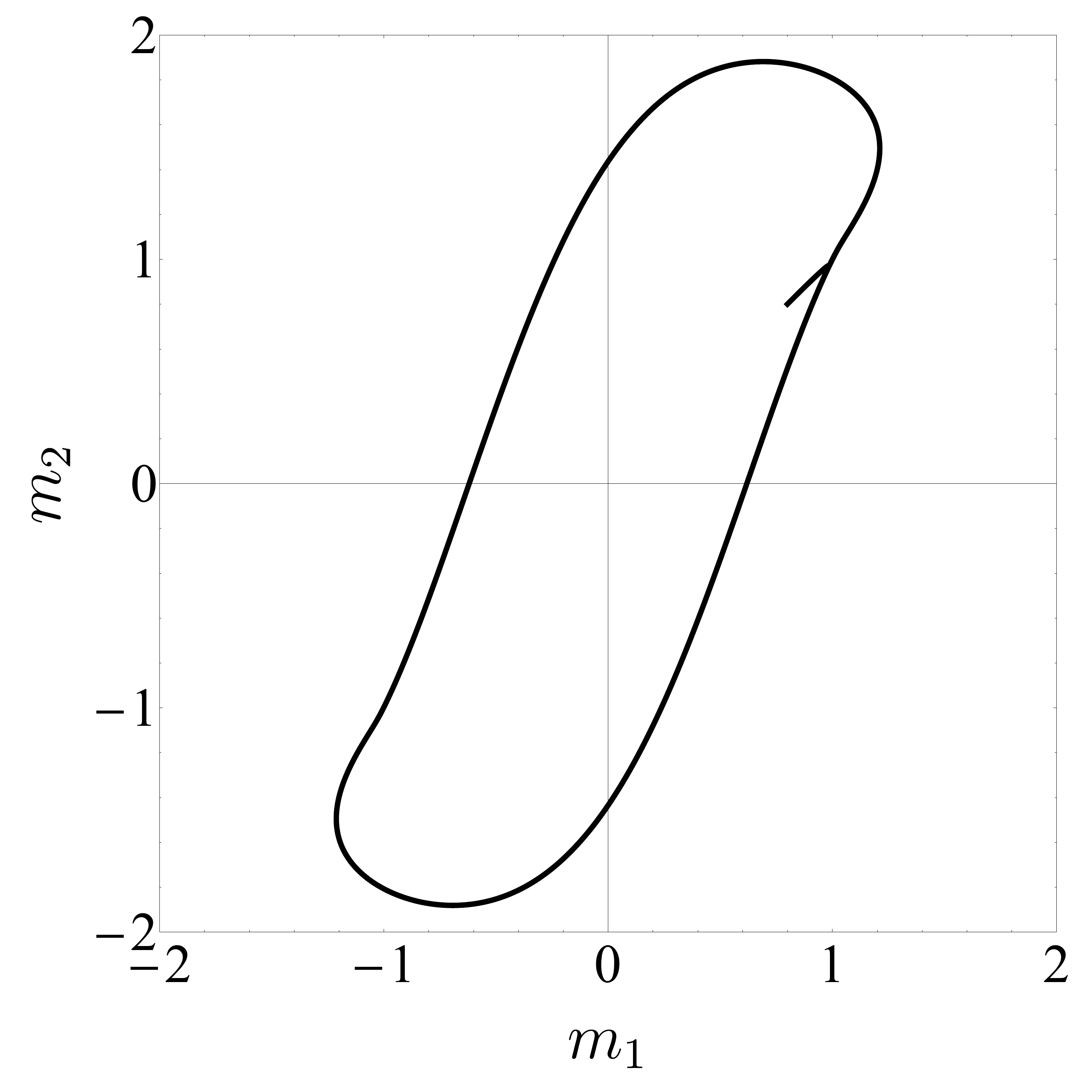}}
    \subfigure
{\includegraphics[width=0.325\textwidth]{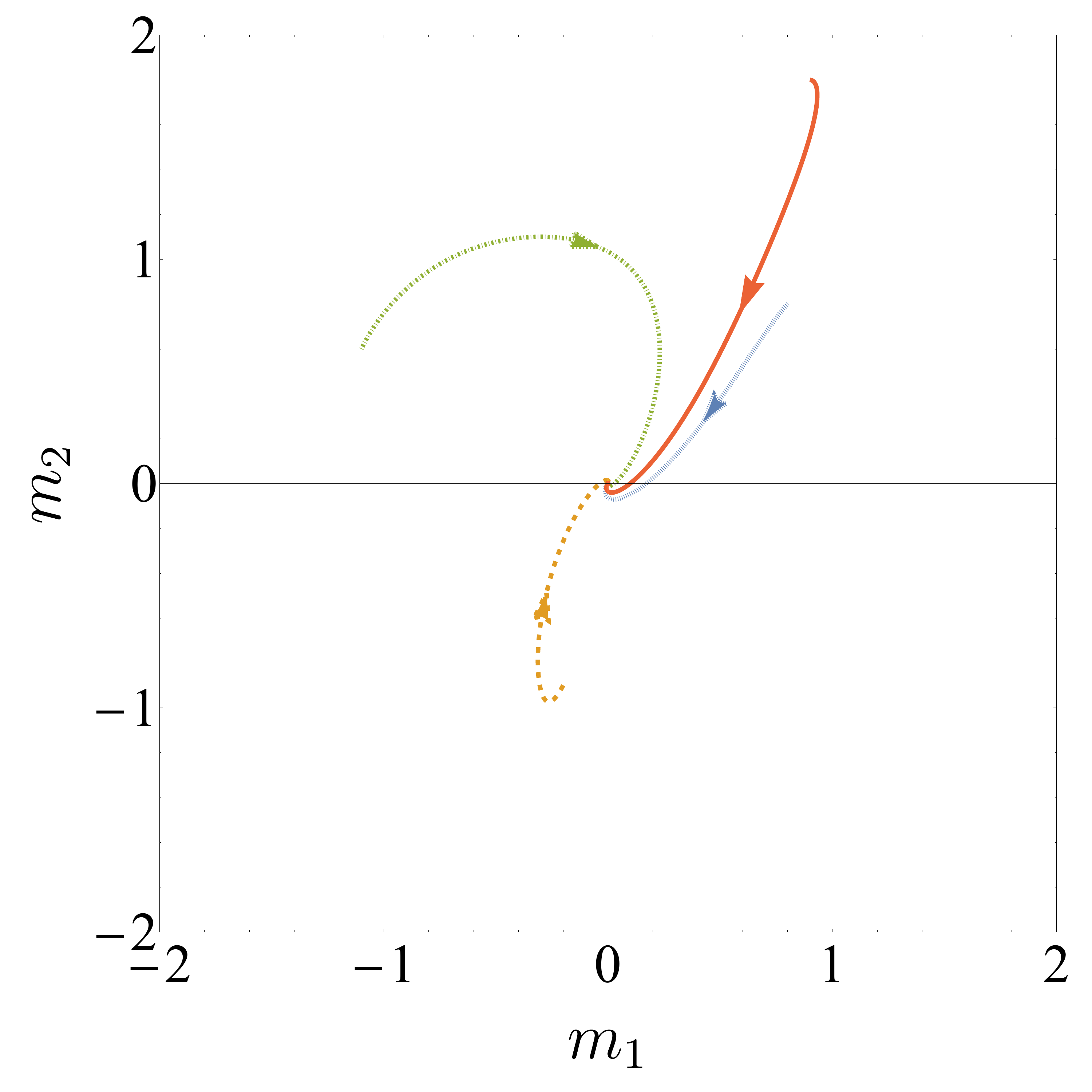}}
    \end{center}
    
    \vspace{5pt}
    \begin{center}
    \subfigure
{\includegraphics[width=0.325\textwidth]{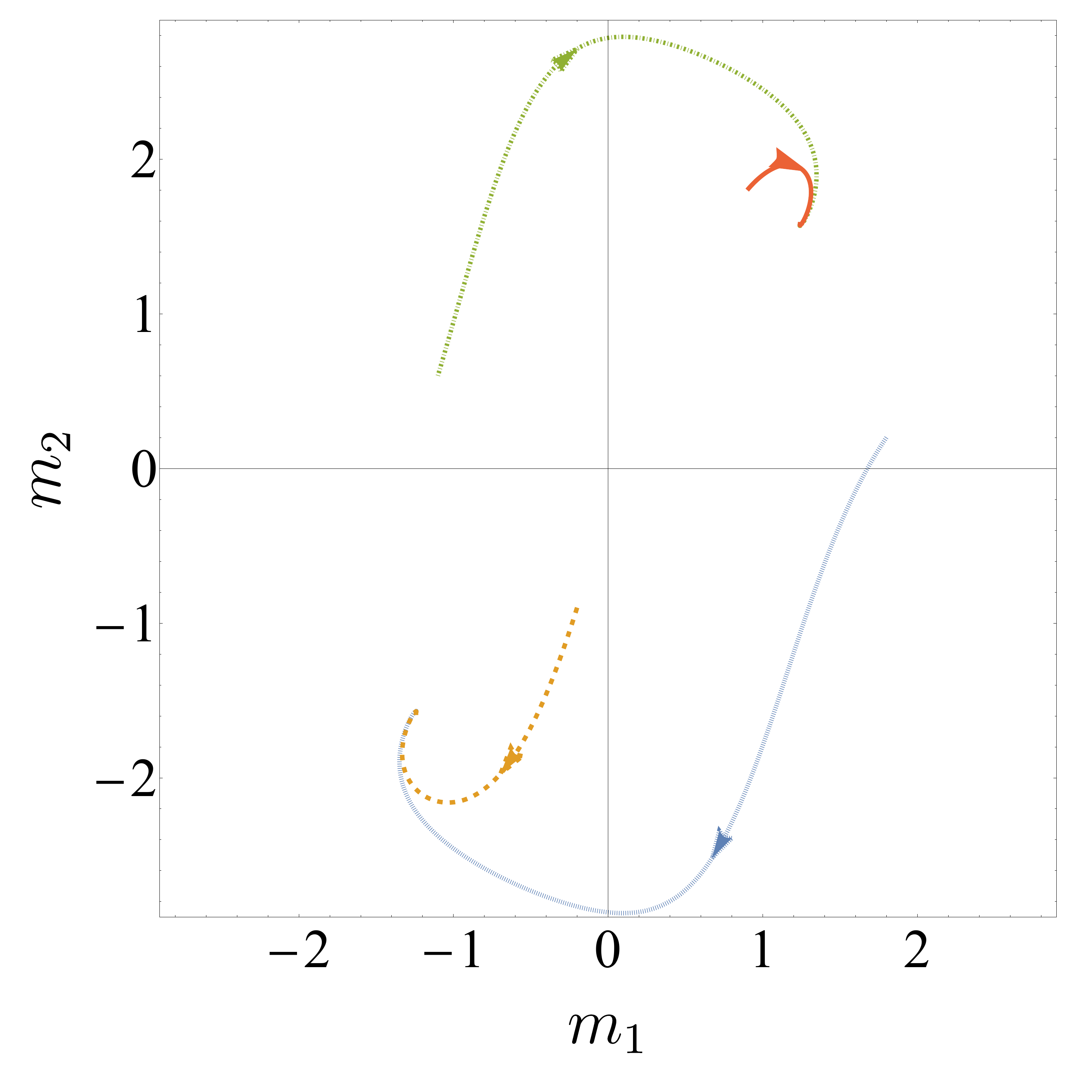}}
    \subfigure
{\includegraphics[width=0.325\textwidth]{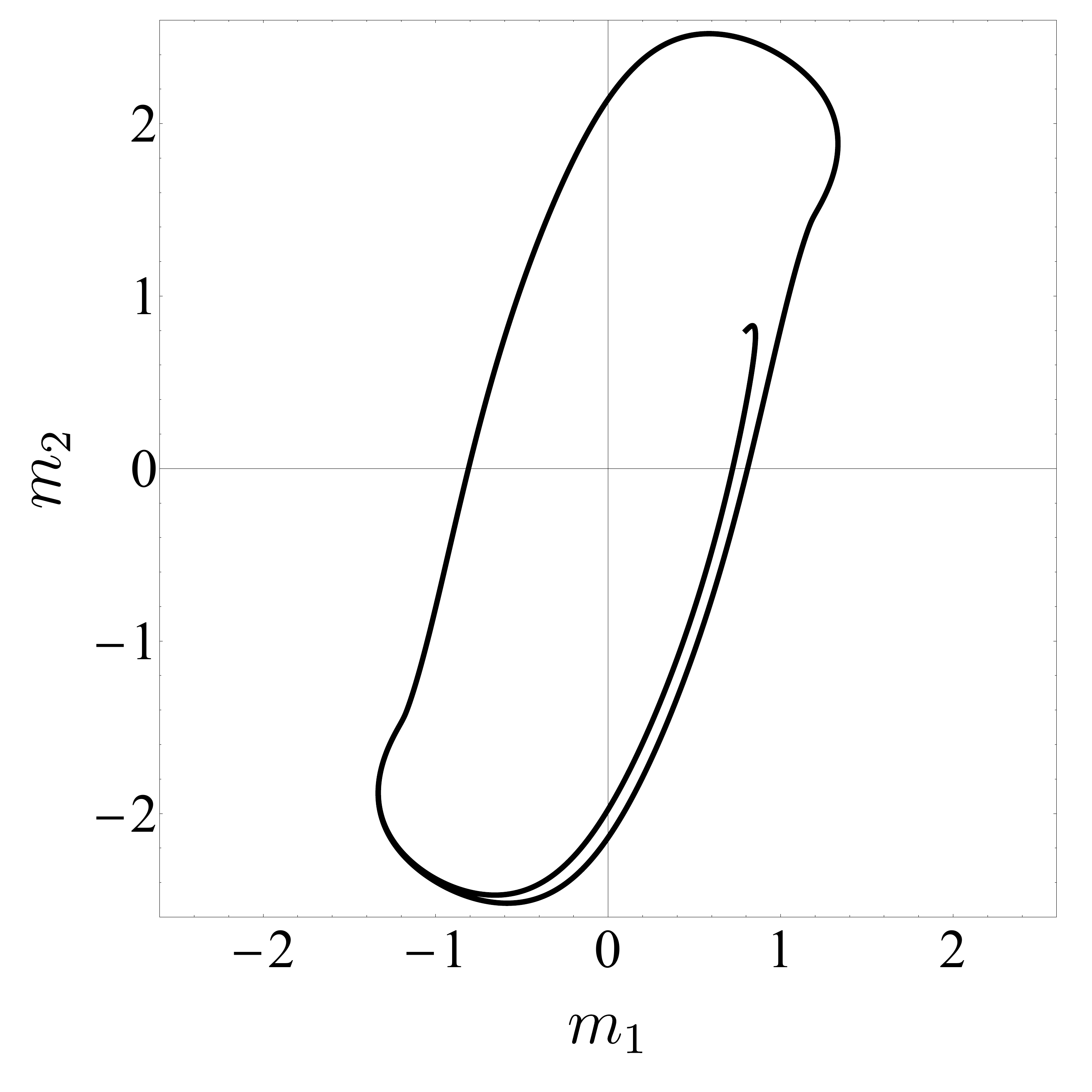}}
    \subfigure
{\includegraphics[width=0.325\textwidth]{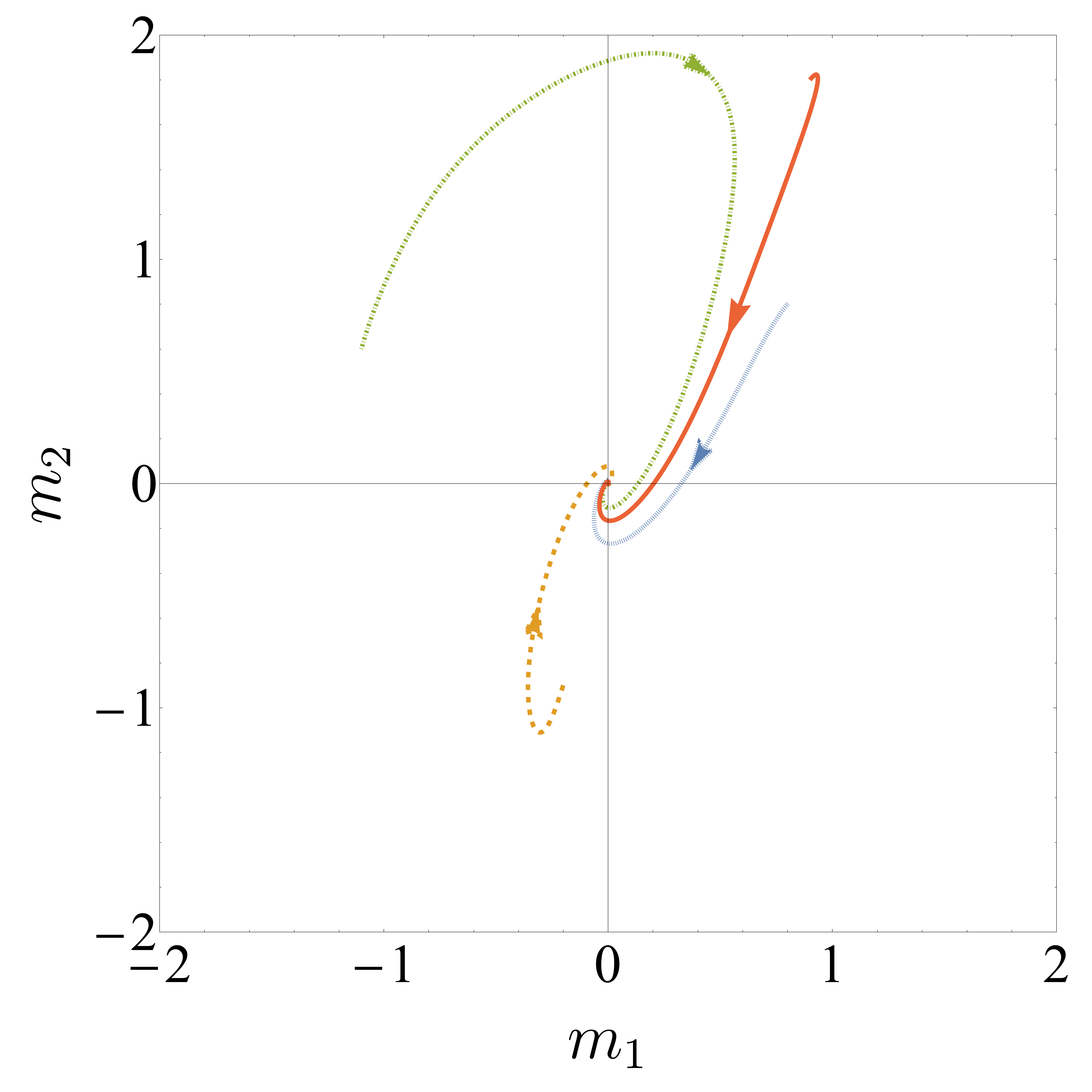}}
    \end{center}
    \caption{\footnotesize{Projected dynamics of system \eqref{eq:odes} in the $\left(m_1,\, m_2\right)$ plane. In all cases, we considered $10^6$ iterations with a time-step $dt=0.005$, $\alpha = 0.5$, $\theta_{11}=\theta_{22}=8$. From top to bottom: $A-1<B<A+2$, in particular, $A=2$ and $B=2.5$; $B=A+2$, in particular, $A=2$ and $B=4$; $B>A+2$, in particular, $A=2$ and $B=7$. In the first column we plot the trajectories of the system in the zero-noise case ($\sigma=0$), in the second column we consider an intermediate intensity for the noise ($\sigma=0.5$, $0.1$ and $0.6$ respectively) and in the third one we set $\sigma=5$.}}
    \label{fig:sigmas_Gauss}
\end{figure}

Simulations of system \eqref{eq:odes}, with values of $A$ and $B$ as in Table \ref{tab:comparison_fft_poincare}  and Fig.~\ref{fig:small_vs_big_noise}, gave the results shown in Fig.~\ref{fig:Gauss:m2:v2} and Fig.~\ref{fig:sigmas_Gauss}, where rhythmic oscillations for intermediate values of noise were detected. 

Our analysis shows that the behavior of system \eqref{micro:dyn} for different noise sizes is well described, at least qualitatively, by the Gaussian approximation  \eqref{eq:odes}. 


\subsection*{Acknowledgment}
We thank Paolo~Dai~Pra and Daniele~Tovazzi for useful discussions, comments and suggestions. MF thanks Serena Bonvicini that partially studied the model numerically in her master thesis.  EM acknowledges financial support from Progetto Dottorati - Fondazione Cassa di Risparmio di Padova e Rovigo. LA, FC and MF did part of this work during a stay at the Institut Henri Poincaré - Centre Emile Borel in occasion of the trimester ``Stochastic Dynamics Out  of Equilibrium''. They thank this institution for hospitality and support.

\begin{appendices}
Here we detail the proofs of the statements in the main text.  

\section{Well-posedness of system \eqref{prop:chaos}}\label{appendixA}

\begin{thm}\label{thm:prop:chaos:existence}
Fix $T>0$. For any initial condition $\left(x(0), y(0)\right) = (\xi_x, \xi_y)$, with $\xi_x$, $\xi_y$ real random variables having finite first moment and  being independent of the Brownian motions $(w_i(t), 0 \leq t \leq T)_{i=1,2}$, system \eqref{prop:chaos} has a unique strong solution.
\end{thm}

\begin{proof}
We follow the argument in \cite{Daw83}, based on a Picard iteration. We define recursively two sequences of stochastic processes $\left(x_n(t), 0 \leq t \leq T\right)$ and $\left(y_n(t), 0 \leq t \leq T\right)$, indexed by $n\geq 1$, via their It{\^o}'s differentials
    \begin{align*}
        dx_n(t) &= \left\{-x^3_n(t) +x_n(t) - \alpha \theta_{11}\Big(x_n(t)-\mathbb{E}\left[x_{n-1}(t)\right]\Big) \right. \\
        &\left.-\left(1-\alpha\right)\theta_{12}\Big(x_n(t)-\mathbb{E}\left[y_{n-1}(t)\right]\Big)\right\}dt +\sigma dw_1(t) \\[.3cm]
        dy_n(t) &= \left\{-y^3_n(t)+y_n(t) -\alpha\theta_{21}\Big(y_n(t)-\mathbb{E}\left[x_{n-1}(t)\right]\Big) \right.\\
        &\left.-\left(1-\alpha\right)\theta_{22}\Big(y_n(t)-\mathbb{E}\left[y_{n-1}(t)\right]\Big)\right\}dt + \sigma dw_2(t),
    \end{align*}
all with the same initial condition $(x_n(0),y_n(0))=(\xi_x,\xi_y)$. By subtracting two subsequent sets of equations, written in integral form, for every $t \in [0,T]$, we obtain
\begin{align}
x_{n+1}(t)-x_n(t) &= \int_0^t \left(x_{n+1}(s)-x_n(s)\right) \left[1-f_n(s)-\alpha\theta_{11}-\left(1-\alpha\right)\theta_{12}\right]ds\nonumber\\
&+\int_{0}^t\left( \alpha \theta_{11}\mathbb{E}\left[x_n(s)-x_{n-1}(s)\right]+ \left(1-\alpha\right)\theta_{12} \mathbb{E}\left[y_n(s)-y_{n-1}(s)\right]\right)ds \label{eq:x1}  
\end{align}
\begin{align}
y_{n+1}(t)-y_n(t) &= \int_0^t \left(y_{n+1}(s)-y_n(s)\right) \left[1-\alpha\theta_{21}-g_n(s)-\left(1-\alpha\right)\theta_{22}\right]ds\nonumber\\
&+\int_{0}^t\left(\alpha\theta_{21} \mathbb{E}\left[x_n(s)-x_{n-1}(s)\right]+\left(1-\alpha\right)\theta_{22} \mathbb{E}\left[y_n(s)-y_{n-1}(s)\right]\right)ds,  \label{eq:y1}
\end{align}
where we have employed the identity $a^3 -b^3 = \left(a-b\right)(a^2 + b^2 + ab)$ and we have set $f_n(s)\coloneqq x^2_{n+1}(s)+x^2_n(s)+x_{n+1}(s) x_n(s)$ and $g_n(s)\coloneqq y^2_{n+1}(s)+y^2_n(s)+y_{n+1}(s) y_n(s)$. Observe that $f_n(t), g_n(t) \geq 0$ for all $t \in [0,T]$.

Eq.~\eqref{eq:x1} and Eq.~\eqref{eq:y1} are of the form $\varphi(t) = \int_0^t \varphi(s)H(s)ds + \int_0^t Q(s)ds$, where $\varphi(t)$ is given by $x_{n+1}(t)-x_n(t) $ and $y_{n+1}(t)-y_n(t)$, respectively. The solution to an equation of this form can be explicitly written as $\varphi(t) = \varphi(0) +\int_{0}^t Q(s) e^{\int_{s}^{t} H(r) dr} ds$. In our case, $\varphi(0) = 0$ since $x_n(0)=\xi_x$ and $y_n(0)=\xi_y$ for all $n\geq 1$ by assumption. 
Therefore, for $t \in[0,T]$, the solutions of Eq.~\eqref{eq:x1} and Eq.~\eqref{eq:y1} are
\begin{align}
x_{n+1}(t)-x_n(t) &= \int_0^t \left\{\alpha\theta_{11}\mathbb{E}\left[x_n(s)-x_{n-1}(s)\right]+ \left(1-\alpha\right)\theta_{12} \mathbb{E}\left[y_n(s)-y_{n-1}(s)\right]\right\} \nonumber\\ 
    &\qquad\quad \times\: e^{\int_s^t \left(1-f_n(r)-\alpha\theta_{11}-\left(1-\alpha\right)\theta_{12}\right)dr } ds \label{eq:existence:x:1} \\[.3cm]
    y_{n+1}(t)-y_n(t) &= \int_0^t \left\{\alpha\theta_{21} \mathbb{E}\left[x_n(s)-x_{n-1}(s)\right]+\left(1-\alpha\right)\theta_{22} \mathbb{E}\left[y_n(s)-y_{n-1}(s)\right]\right\}\nonumber\\ 
    &\qquad\quad \times\:e^{\int_s^t \left(1-\alpha\theta_{21}-g_n(r)-\left(1-\alpha\right)\theta_{22}\right)dr} ds. \label{eq:existence:y:1}
\end{align}

We now get into the core of the proof. \smallskip\\

\emph{Step~$1$: auxiliary property.} We will show that, for every $T>0$, the sequences $\left(\mathbb{E}\left[x_n(t)\right], 0 \leq t \leq T\right)_{n \geq 1}$ and $\left(\mathbb{E}\left[y_n(t)\right], 0 \leq t \leq T\right)_{n\geq 1}$ are Cauchy sequences in the space $\mathcal{C}\left([0,T]\right)$, equipped with the supremum norm
$$d(f,g)\coloneqq \sup_{t\in [0,T]} \left| f(t)-g(t)\right| \quad \forall f,\,g \in \mathcal{C}\left([0,T]\right). $$
As a consequence, since $(\mathcal{C}([0,T]), d)$ is a complete metric space, we will obtain convergence to elements $(m_x(t), 0 \leq t \leq T), (m_y(t), 0 \leq t \leq T) \in \mathcal{C}\left([0,T]\right)$. 

\sloppy We take the absolute value and the expectation in both Eq.~\eqref{eq:existence:x:1} and Eq.~\eqref{eq:existence:y:1}. If we denote by $\phi_n(t) \coloneqq \sup_{s \in [0,t]} \mathbb{E}\left[\left|x_{n+1}(s)-x_n(s) \right|\right]$ and $\psi_n(t) \coloneqq \sup_{s\in [0,t]}\mathbb{E}\left[\left|y_{n+1}(s)-y_n(s) \right| \right]$, from Eq.~\eqref{eq:existence:x:1} we obtain
\begin{equation}\label{inequality:phi_n}
        \phi_n(t) \leq \tilde C_t\int_0^t \phi_{n-1}(r)dr + \tilde D_t\int_0^t \psi_{n-1}(r)dr,
\end{equation}
for some positive constants $\tilde C_t$ and $\tilde D_t$. From Eq.~\eqref{eq:existence:y:1} we get an analogous inequality for $\psi_{n}(t)$. The inequalities are valid for all $t \in [0,T]$. 
Iteratively employing inequality \eqref{inequality:phi_n} gives
\begin{equation*}
\phi_n(T) \leq C_T \phi_1(T) \frac{T^{n-1}}{(n-1)!} + D_T \psi_1(T) \frac{T^{n-1}}{(n-1)!}, 
\end{equation*}
for suitable positive constants $C_T$ and $D_T$. Similarly we bound $\psi_n(T)$. Hence, $\phi_n(T)$ and $\psi_n(T)$ go to zero as $n\to +\infty$. Thus, it follows that $\left(\mathbb{E}\left[x_n(t)\right], 0 \leq t \leq T\right)_{n\geq 1}$ and $\left(\mathbb{E}\left[y_n(t)\right], 0 \leq t \leq T\right)_{n\geq 1}$ are Cauchy sequences in $\mathcal{C}\left([0,T]\right)$ and converge to the continuous limits $(m_x(t), 0 \leq t \leq T)$ and $(m_y(t), 0 \leq t \leq T)$, respectively. \smallskip\\

\emph{Step~$2$: existence of the solution to \eqref{prop:chaos}.} Consider the following system of stochastic differential equations
    \begin{align}\label{eq:intermediate:C}
        dx(t)  &= \left[-x^3(t) +x(t) -\alpha\theta_{11} \Big(x(t)-m_x(t)\Big)- \left(1-\alpha\right)\theta_{12}\Big(x(t)-m_y(t)\Big) \right]dt + \sigma dw_1(t) \nonumber\\[.2cm]
        dy(t)  &= \left[-y^3(t) +y(t) -\alpha\theta_{21} \Big(y(t)-m_x(t)\Big) - \left(1-\alpha\right) \theta_{22}\Big(y(t)-m_y(t)\Big) \right]dt + \sigma dw_2(t),
    \end{align}
with initial condition $(x(0),y(0))=(\xi_x,\xi_y)$. Since the functions $m_x$ and $m_y$ are bounded for every $t \in [0,T]$, existence and uniqueness of a strong solution for \eqref{eq:intermediate:C} follows from a Khasminskii's test with norm-like function $V(x,y)=\frac{x^4}{4}+\frac{x^2}{2}+\frac{y^4}{4}+\frac{y^2}{2}$. See \cite{Kha12,MeTw93}. 

Let $((x(t),y(t)), 0 \leq t \leq T)$ be the unique strong solution for \eqref{eq:intermediate:C}. We construct the differences
\begin{align*}
x_{n+1}(t) - x(t)  &= \int_0^t
     \left(x_{n+1}(s)-x(s) \right) \left[ 1-f_n(s)-\alpha\theta_{11}-\left(1-\alpha\right)\theta_{12}\right] ds \nonumber\\
        & + \int_{0}^t \left(\alpha\theta_{11}\mathbb{E}\left[x_n(s)-m_x(s)\right] + \left(1-\alpha\right)\theta_{12}\mathbb{E}[y_n(s)-m_y(s) ] \right)ds \\[.4cm]
        y_{n+1}(t) - y(t) &= \int_0^t \left(y_{n+1}(s)-y(s) \right) \left[ 1 -\alpha\theta_{21} -g_n(s) -\left(1-\alpha\right)\theta_{22}\right] ds\nonumber\\
        & + \int_{0}^t \left(-\alpha\theta_{21} \mathbb{E}\left[x_n(s)-m_x(s)\right]
        -\left(1-\alpha\right)\theta_{22} \mathbb{E}[y_n(s)-m_y(s) ] \right)ds,
\end{align*}
with $f_n(s)\coloneqq x^2_{n+1}(s)+x^2(s)+x_{n+1}(s) x(s)$ and $g_n(s)\coloneqq y^2_{n+1}(s)+y^2(s)+y_{n+1}(s) y(s)$, and we repeat the same argument as in Step~$1$. As a consequence, we get that 
\[
(\mathbb{E}\left[x_n(t) \right], 0 \leq t \leq T) \xrightarrow{\:\, n \to \infty \:\,} (\mathbb{E}\left[x(t)\right], 0 \leq t \leq T)
\]
and 
\[
(\mathbb{E}\left[y_n(t)\right], 0 \leq t \leq T) \xrightarrow{\:\, n \to \infty \:\,} (\mathbb{E}\left[y(t)\right], 0 \leq t \leq T)
\]
in $\mathcal{C}\left([0,T]\right)$. Therefore, since we have already showed that 
\[
(\mathbb{E}\left[x_n(t) \right], 0 \leq t \leq T) \xrightarrow{\:\, n \to \infty \:\,} (m_x(t), 0 \leq t \leq T)
\]
and 
\[
(\mathbb{E}\left[y_n(t)\right], 0 \leq t \leq T) \xrightarrow{\:\, n \to \infty \:\,} (m_y(t), 0 \leq t \leq T),
\]
we have that $\mathbb{E}\left[x(t)\right] = m_x(t)$ and $\mathbb{E}\left[y(t)\right] = m_y(t)$ for all $t \in [0,T]$. Hence, system \eqref{eq:intermediate:C} coincides with system \eqref{prop:chaos} and its solution $\left((x(t), y(t)), 0 \leq t \leq T \right)$ provides a solution for \eqref{prop:chaos}. \smallskip\\

\emph{Step~$3$: uniqueness of the solution to \eqref{prop:chaos}.} Let $((u(t),v(t)), 0 \leq t \leq T)$ be another solution to \eqref{prop:chaos}. We write the integral equations for $x(t)-u(t)$ and $y(t)-v(t)$ and we use them to get estimates for $\phi(t)=\vert \mathbb{E}[x(t)-u(t)]\vert$ and $\psi(t)=\vert \mathbb{E}[y(t)-v(t)]\vert$. By mimicking the computations above, we obtain
\[
\phi(t) \leq \hat{C}_{T} \int_0^t (\phi(s)+\psi(s)) \, ds 
\quad \text{ and } \quad
\psi(t) \leq \hat{D}_{T} \int_0^t (\phi(s)+\psi(s)) \, ds,
\]
for suitable positive constants $\hat{C}_{T}$ and $\hat{D}_{T}$. Summing up the two previous inequalities and using Gronwall's lemma yields $\phi(t)+\psi(t) \leq 0$ for all $t \in [0,T]$. This shows that $\mathbb{E}[x(t)]=\mathbb{E}[u(t)]$ and $\mathbb{E}[y(t)]=\mathbb{E}[v(t)]$ for all $t \in [0,T]$. Thus, $((x(t),y(t)), 0 \leq t \leq T)$ and $((u(t),v(t)), 0 \leq t \leq T)$ are both solutions to \eqref{eq:intermediate:C} with the same pair $(m_x,m_y)$ and the same initial condition. It follows that $(x(t),y(t))=(u(t),v(t))$ for all $t \in [0,T]$.
\end{proof}

\section{Proof of Theorem \ref{thm:prop:caos}}\label{appendix:B}

The proof of Theorem~\ref{thm:prop:caos} is standard and it relies on a coupling method \cite{CoDaPFo15, Luisa2018}. We want to prove \eqref{prop_chaos_on_average}. Without loss of generality, we take $\mathcal{I}=\{1,\ldots, k_1\}$ and $\mathcal{J}=\{1,\ldots, k_2\}$. Since
\begin{equation*}
\begin{split}
\mathbb{E}\left[\sup_{t\in[0,T]}\left|\mathbf{z}_{k_1,k_2}^{(N)}(t)-\mathbf{z}_{k_1,k_2}(t)\right|\right] 
&\leq 
\sum_{j=1}^{k_1} \mathbb{E}\left[\sup_{t\in [0,T]} \left| x^{(N)}_j(t)-x_j(t)\right| \right] \\
&+ \sum_{j=1}^{k_2} \mathbb{E}\left[\sup_{t\in [0,T]} \left| y^{(N)}_j(t)-y_j(t)\right| \right],
\end{split}
\end{equation*}
to conclude it suffices to show that each of the $k_1+k_2$ terms goes to zero in the limit $N\rightarrow\infty$. 

In the sequel we will consider only the term  $\mathbb{E}\left[\sup_{t\in [0,T]} \left| x^{(N)}_1(t)-x_1(t)\right| \right]$; the other terms can be dealt with similarly. We only sketch our computations as they use the very same tricks as in the proof of Theorem~\ref{thm:prop:chaos:existence} in Appendix~\ref{appendixA}. Since the processes $\left(x_1^{(N)}(t), 0 \leq t \leq T\right)$ and $(x_1(t), 0 \leq t \leq T)$ are initiated at the same position, we obtain
\begin{equation}\label{eq:1x}
x^{(N)}_1(t)-x_1(t) = \int_{0}^t \left[ \left(x^{(N)}_1(s)-x_1(s) \right)\left(1-f(s)-\alpha\theta_{11}-(1-\alpha)\theta_{12}\right) + \mu(s) \right] ds,
\end{equation}    
where we have employed the identity $a^3-b^3 = (a-b)(a^2 +b^2+ ab)$ and we have set $f(s) \coloneqq \left(x^{(N)}_1(s)\right)^2+x_1^2(s) +x^{(N)}_1(s) x_1(s)$, and where
$$
\mu(s) \coloneqq \alpha \theta_{11} \left( m^{(N)}_1(s) - \mathbb{E}[x_1(s)]\right) + (1-\alpha)\theta_{12}\left( m^{(N)}_2(s) 
- \mathbb{E}[y_1(s)]\right).
$$
Observe that Eq.~\eqref{eq:1x} is of the form $\varphi(t) = \int_{0}^t \varphi(s)H(s) ds + \int_{0}^t Q(s) ds$, with $\varphi(t)=x^{(N)}_1(t)-x_1(t)$.
Therefore, as the solution of the latter equation is $\varphi(t) = \varphi(0) +\int_{0}^t Q(s) e^{\int_{s}^{t} H(r) dr} ds$ and, in our case, $x^{(N)}_1(0)=x_1(0)$ by assumption, for every $t \in [0,T]$, we can estimate
\begin{equation*}
\left| x^{(N)}_1(t) - x_1(t) \right| \leq \int_0^t \left|\mu(s) \right| e^{\int_{s}^{t}\left(1-f(r)-\alpha\theta_{11}-\left(1-\alpha\right)\theta_{12}\right) dr} ds \leq C_{T}\int_{0}^{t}\left| \mu(s)\right| ds,
\end{equation*}
for some positive constant $ C_{T}$. At this point, taking the supremum and the expectation of both sides of the previous inequality, for $\tilde{t} \in [0,T]$, we obtain 
\begin{equation}\label{eq:2x}
\mathbb{E}\left[\sup_{t \in [0,\tilde t]}\left| x^{(N)}_1(t) - x_1(t)\right| \right] \leq C_{T}\int_{0}^{\tilde t} \mathbb{E}\left[\left| \mu(s)\right| \right] ds.
\end{equation}
We need an upper bound for $\mathbb{E}\left[\left| \mu(s)\right| \right]$. We have
\begin{align*}
\mathbb{E}\left[\left| \mu(s)\right|\right] 
& \leq \mathbb{E}\left[ \alpha\theta_{11}\left| m^{(N)}_1(s)- \mathbb{E}[x_1(s)] \right| + (1-\alpha)\theta_{12}\left| m^{(N)}_2(s)- \mathbb{E}[y_1(s)] \right|\right] 
\end{align*}
and, by adding and subtracting $\frac{1}{N_1}\sum_{i=1}^{N_1}x_i(s)$ (resp. $\frac{1}{N_2}\sum_{i=1}^{N_2}y_i(s)$) inside the first (resp. second) absolute value, we get
\begin{align*}
\mathbb{E}\left[\left| \mu(s)\right|\right] 
&\leq \alpha \theta_{11} \left\{ \frac{1}{N_1}\sum_{i=1}^{N_1}\mathbb{E}\left[\left| x^{(N)}_i(s)-x_i(s)\right|\right] + \mathbb{E}\left[\left|\frac{1}{N_1}\sum_{i=1}^{N_1}x_i(s) - \mathbb{E}[x_1(s)] \right|\right] \right\} \\
&+ (1-\alpha)\theta_{12} \left\{ \frac{1}{N_2}\sum_{i=1}^{N_2}\mathbb{E}\left[\left| y^{(N)}_i(s)-y_i(s)\right|\right] + \mathbb{E}\left[\left| \frac{1}{N_2} \sum_{i=1}^{N_2}y_i(s) - \mathbb{E}[y_1(s)]\right|\right]\right\}.
\end{align*}
Since the limiting variables $(x_i(t))_{i = 1,\dots, N_1}$ and $(y_i(t))_{i = 1,\dots, N_2}$ are i.i.d. families and have uniformly bounded second moments for all $t \in [0,T]$ (due to the well-posedness of system \eqref{prop:chaos}), the standard CLT assures that there exists a positive constant $K_T$ such that, uniformly for all $s \in [0,T]$, it holds  
$$
\mathbb{E}\left[\left| \frac{1}{N_1} \sum_{i=1}^{N_1}x_i(s) - \mathbb{E}[x_1(s)]\right|\right]\leq \frac{K_T}{\sqrt{N_1}} $$
and
$$\mathbb{E}\left[\left| \frac{1}{N_2} \sum_{i=1}^{N_2}y_i(s) - \mathbb{E}[y_1(s)]\right|\right] \leq \frac{K_T}{\sqrt{N_2}}.
$$
Moreover, we have
$$\mathbb{E}\left[\left| x^{(N)}_i(s)-x_i(s)\right|\right] \leq \mathbb{E}\left[ \sup_{r \in [0,s]}\left| x^{(N)}_i(s)-x_i(s)\right|\right]$$
and
$$\mathbb{E}\left[\left| y^{(N)}_i(s)-y_i(s)\right|\right] \leq \mathbb{E}\left[ \sup_{r \in [0,s]}\left| y^{(N)}_i(s)-y_i(s)\right|\right].$$
These last terms are in fact independent of the index $i$ due to the symmetry of the system which, in turn, is due to the choice of the initial conditions and the mean field assumption. Thus, recalling that $\alpha=\frac{N_1}{N}$, we obtain
\begin{align}\label{eq:3x}
\mathbb{E}\left[\left| \mu(s)\right|\right] 
&\leq \alpha\theta_{11} \mathbb{E}\left[\sup_{r \in [0,s]}\left| x^{(N)}_1(r)-x_1(r)\right|\right] + \frac{\sqrt{\alpha} \, \theta_{11}K_T}{\sqrt{N}} \nonumber \\
&+ (1-\alpha)\theta_{12}\mathbb{E}\left[\sup_{r \in [0,s]}\left| y^{(N)}_1(r)-y_1(r)\right|\right] + \frac{\sqrt{1-\alpha} \, \theta_{12} K_T}{\sqrt{N}}.
\end{align}
By employing Eq.~\eqref{eq:3x} in Eq.~\eqref{eq:2x}, it is easily seen that there exists a constant $D$, depending on $T$ and on the parameters $\alpha, \theta_{11},  \theta_{22},  \theta_{12}$, and $\theta_{21}$, such that 
\begin{align}\label{eq:4x}
\mathbb{E}\left[\sup_{t\in [0,\tilde t]}\left|x^{(N)}_1(t)-x_1(t) \right|\right]
&\leq D \int_{0}^{\tilde t}\mathbb{E}\left[\sup_{r\in [0,s]}\left| x^{(N)}_1(r) -x_1(r)\right| \right] ds \nonumber\\
&+ D\int_{0}^{\tilde t} \mathbb{E}\left[\sup_{r\in [0,s]}\left| y^{(N)}_1(r) -y_1(r)\right| \right] ds  \nonumber\\
&+ \frac{D}{\sqrt{N}},
\end{align}
Similarly, we obtain also  
\begin{align}\label{eq:4y}
\mathbb{E}\left[\sup_{t\in [0,\tilde t]}\left|y^{(N)}_1(t)-y_1(t) \right|\right] &\leq D \int_{0}^{\tilde t}\mathbb{E}\left[\sup_{r\in [0,s]}\left| x^{(N)}_1(r) -x_1(r)\right| \right] ds\nonumber\\
&+ D\int_{0}^{\tilde t} \mathbb{E}\left[\sup_{r\in [0,s]}\left| y^{(N)}_1(r) -y_1(r)\right| \right] ds \nonumber\\
&+ \frac{D}{\sqrt{N}}.
\end{align}
If we define
\begin{equation*}
g\left(\tilde{t} \right) \coloneqq  \mathbb{E}\left[\sup_{t \in[0,\tilde{t}]}\left| x^{(N)}_1(t)-x_1(t)\right|\right] + \mathbb{E}\left[\sup_{t \in[0,\tilde{t}]}\left| y^{(N)}_1(t)-y_1(t)\right|\right],
\end{equation*}
by summing up the inequalities \eqref{eq:4x} and \eqref{eq:4y}, we obtain
\begin{equation*}
 g\left(\tilde t\right) \leq2\: D \int_{0}^{\tilde t} g(s) ds + 2\:\frac{ D}{\sqrt{N}}.
\end{equation*}
An application of Gronwall's lemma leads to the conclusion. Indeed, we get the inequality $g(T) \leq \frac{2De^{2DT}}{\sqrt{N}}$, whose right-hand side goes to zero as $N \to +\infty$.

\section{Equilibria of the noiseless dynamics}\label{appendix:C}
We consider system \eqref{micro:dyn:nonoise:AB} and study the nature of its fixed points, depending on the values of the parameters $A$ and $B$. Throughout this analysis, we use the basic theory of dynamical systems as it can be found for instance in \cite{Per2001}. Moreover, unless otherwise specified, we assume for the moment $A>1$ and $B>0$.

\begin{enumerate}
\item The fixed points $\left(0,0\right)$ and $\pm\left(1,1\right)$ are present for any values of $A$ and $B$.  
\begin{enumerate}
\item The linearized system around the origin has the eigenvalues
$$\lambda_1=1\mbox{ and }\lambda_2=1-A+B=1-\gamma,$$
where $\gamma \coloneqq A-B$.
Thus $\left(0,0\right)$ is a saddle when $\gamma> 1$, it has an unstable and a neutral direction for $\gamma= 1$ and it is an unstable node otherwise.
\item Eigenvalues of the linearized system around $\pm\left(1,1\right)$ are
$$\lambda_1=-2 \mbox{ and } \lambda_2=-2-A+B=-2-\gamma.$$
As a consequence, $\pm\left(1,1\right)$ are stable nodes for $\gamma> -2$. They have a neutral and a stable direction when $\gamma=-2$ and they are saddle points otherwise.
\end{enumerate}
In conclusion:  i) for $\gamma<-2$, $\left(0,0\right)$ is unstable and  $\pm\left(1,1\right)$ are saddle points; ii) when  $-2<\gamma<1$, $\left(0,0\right)$ is unstable and  $\pm\left(1,1\right)$ are stable nodes; iii) for $\gamma > 1$,  $\left(0,0\right)$ is a saddle point and $\pm\left(1,1\right)$ are stable nodes.
\item Depending on the values of the parameters $A$ and $B$, there might be two additional equilibria. We search for equilibria of the form $(x,\beta x)$ with $\beta\neq 0$, $x\neq 0$. Notice that all the possible equilibria except for $(0,0)$ have such a form.\\
Substituting $y=\beta x$ in the first equation of \eqref{micro:dyn:nonoise:AB}, we get 
\begin{equation}\label{eq:x(beta)}
\bar x_{\beta}=\pm\left(\sqrt{1-A\left(1-\beta\right)},\,\beta\sqrt{1-A\left(1-\beta\right)}\right),
\end{equation}
subject to the condition 
\begin{equation}\label{eq:conditionA>1}
\beta >\frac{A-1}{A}.   
\end{equation}
Notice that no $\beta<0$ fulfills \eqref{eq:conditionA>1}, since we have assumed $A>1$. 
Therefore, system \eqref{micro:dyn:nonoise:AB} can possibly have extra fixed points of the form $\left(x, \beta x\right)$ only if they lie in the first or the third quadrant.\\
The second equation in \eqref{micro:dyn:nonoise:AB} leads to the fixed point equation
\begin{equation}\label{eq:f(beta)=beta}
\beta = f\left(\beta\right) \quad \text{ with } \quad f\left(\beta\right) := \sqrt{\frac{1-B\frac{1-\beta}{\beta}}{1-A\left(1-\beta\right)}}.
\end{equation}
Observe that, for $\beta = 1$, we recover the equilibria $\pm\left(1,1\right)$. Therefore, fixed points of the type $\bar{x}_{\beta}$ may exist if condition \eqref{eq:conditionA>1} is satisfied and Eq.~\eqref{eq:f(beta)=beta} has real solutions, that is if
\begin{equation}\label{eq:two_conditions}
\beta>\max\left\{\frac{A-1}{A}, \frac{B}{1+B}\right\},
\end{equation}
which is equivalent to
\begin{equation*}
    \begin{cases}
    \beta>\frac{A-1}{A} = \frac{B}{1+B}, & \text{if } B=A-1\\
    \beta>\frac{A-1}{A}, & \text{if } B<A-1\\
    \beta>\frac{B}{1+B}, & \text{if } B>A-1.
    \end{cases}
\end{equation*}
Therefore, we have the following.
\begin{itemize}
    \item If $B=A-1$, Eq.~\eqref{eq:f(beta)=beta} becomes
    $\beta = \frac{1}{\sqrt{\beta}}$, whose unique solution is $\beta =1$.
    In this case, $\gamma = 1$, so $\pm\left(1,1\right)$ are stable nodes and $\left(0,0\right)$ has a zero eigenvalue, thus it is not a hyperbolic fixed point and the linearization cannot give information about the phase portrait close to it.
    The dynamical system \eqref{micro:dyn:nonoise:AB} can be rewritten as 
    \begin{align}\label{eq:origin:gamma1}
    \dot{x} &= -x^{3} - x\left(A-1\right) +A y\nonumber\\
    \dot{y} &= -y^{3} + x -A\left(x-y\right).
    \end{align}
Observe that the linear terms in both the components of the vector field in Eq.~\eqref{eq:origin:gamma1} are positive above the line $y=\frac{A-1}{A}x$ and negative below it. Thus, the third-order terms can be neglected close to the equilibrium and the linearization gives an accurate sketch of the phase portrait of the system locally.
    Along the line $y=\frac{A-1}{A}x$, that is the eigendirection of the zero eigenvalue, only the third-order terms count and we get $\dot{x}<0$, $\dot{y}<0$ in the first quadrant and $\dot{x}>0$, $\dot{y}>0$ in the third one.
    Overall, the eigendirection of the zero eigenvalue is locally stable.

    \item If $B<A-1$, due to condition \eqref{eq:two_conditions}, we expect to find solutions to Eq.~\eqref{eq:f(beta)=beta} only for $\beta>\frac{A-1}{A}$. 
    Observe that $f\left(\beta\right)$ has a vertical asymptote to positive infinity as $\beta$ approaches $\frac{A-1}{A}$ and it has a horizontal asymptote to zero as $\beta$ grows to infinity. 
Moreover, $\frac{\partial f}{\partial \beta}\left(\beta\right) = 0$ when 
\begin{equation}\label{beta:pm}
\beta_\pm= \frac{A B \pm\sqrt{A B\left(B-(A-1)\right)}}{A\left(1+B\right)},
\end{equation}
which are  complex for $B<A-1$. Therefore,  $f\left(\beta\right)$ is strictly decreasing. 
Thus, its graph cannot have more than one intersection with the line $y=\beta$ and this unique intersection must be at $\beta =1$.
Overall, when $B<A-1$, no fixed points of the type $\left(x, \beta x\right)$ are present, except for $\pm \left(1,1\right)$.
In this setting, we already established that $\left(0,0\right)$ is a saddle point and $\pm\left(1,1\right)$ are stable nodes.

\item If $B>A-1$, we have already pointed out that $\left(0,0\right)$ is unstable, while the nature of $\pm\left(1,1\right)$ can change according to $\gamma$ being less than or greater than $-2$. From condition \eqref{eq:two_conditions}, it follows that we have to look for solutions to Eq.~\eqref{eq:f(beta)=beta} for $\beta > \frac{B}{B+1}$.
\sloppy Observe that $f\left(\frac{B}{1+B}\right) = 0$ and $\lim_{\beta\to +\infty}f\left(\beta\right) = 0$.
The points $\beta_\pm$, given in \eqref{beta:pm}, where $\frac{\partial f}{\partial\beta}\left(\beta\right) = 0$, are real and distinct in this case. Moreover, $\beta_{-}<\frac{B}{1+B}$. Hence, the function $f\left(\beta\right)$ has only a critical point (maximum) at $\beta=\beta_{+} > \frac{B}{1+B}$, so it may cross the line $y=\beta$ once, twice or never. 
    In particular, since we know the solution $\beta=1$ to be always present, the intersections might coincide ($\beta=1$ itself) or be distinct ($\beta=1$ and a second intersection, for $\beta$ greater or less than $1$).
    
    We distinguish three subcases.
    \begin{itemize}
        \item The intersection at $\beta=1$ is the unique solution to Eq.~\eqref{eq:f(beta)=beta} if the graphs of $y=\beta$ and $y=f\left(\beta\right)$ are tangent at that point, i.e., if $\frac{\partial f}{\partial \beta}\left(1\right)=1$.
        This holds only if $B=A+2$.  
        In this case, the analysis of the linearized system tells us that $\pm\left(1,1\right)$ have a negative and a zero eigenvalue and to check stability one has to take into account higher-order terms.

We study only the point $\left(1,1\right)$, the analysis of $-(1,1)$ being similar. To make our computations easier, we translate the vector field so that the fixed point $\left(1,1\right)$ is shifted to $\left(0,0\right)$. We make the change of variables \mbox{$\hat{x}=x-1$} and $\hat{y}=y-1$. In the new coordinates $(\hat{x},\hat{y})$, system \eqref{micro:dyn:nonoise:AB} becomes
\begin{align}\label{micro:dyn:nonoise:eta=2:trasl}
\dot{\hat{x}} &= - \left(A+2\right)\hat{x}+A \hat{y} -3\hat{x}^{2}-\hat{x}^{3}\nonumber\\
\dot{\hat{y}} &= - \left(A+2\right)\hat{x}+A \hat{y} -3\hat{y}^{2}-\hat{y}^{3}.
\end{align}
The line $\hat{y}=\frac{A+2}{A}\hat{x}$, along which the first-order terms in \eqref{micro:dyn:nonoise:eta=2:trasl} vanish, that is the eigendirection of the zero eigenvalue, always lies above the line $\hat{y}=\hat{x}$, that is the eigendirection of the non-zero eigenvalue. 
The first-order terms in Eq.~\eqref{micro:dyn:nonoise:eta=2:trasl} are positive above the line $\hat{y}=\frac{A+2}{A}\hat{x}$ and negative below it. 
So, out of this line, higher-order terms can be neglected close to the origin, whereas the second-order terms become non-negligible as soon as we are along that line, where it is immediate to see that the vector field points downward-left.
        
\item If $0<\frac{\partial f}{\partial \beta}\left(1\right) <1$, i.e., if $B<A+2$, two intersections are present, one at $\beta=1$ and one at $\beta=\beta_{\times}(A,B) < 1$.
The smallest solution to Eq.~\eqref{eq:f(beta)=beta} gives rise to two extra fixed points of the type \eqref{eq:x(beta)}, with both coordinates smaller than $1$ in absolute value.
        
\item If $\frac{\partial f}{\partial \beta}\left(1\right) >1$, i.e., $B>A+2$, in addition to the intersection at $\beta=1$, we have a second intersection at $\beta=\beta_{\times}(A,B)>1$ and we get two fixed points of the type \eqref{eq:x(beta)}, with both coordinates greater than $1$ in absolute value. 
\end{itemize}
\end{itemize}
\end{enumerate}

In what follows we restrict to the case $B>A-1$ and we examine in more detail what happens in the three cases that we considered in the main text  and that are shown in Fig.~\ref{fig:PhDiaNoNoiseNoInternalInteractions}. The point $(0,0)$ is an unstable fixed point in all the scenarios. We will give information on the other equilibria.

\noindent \textbf{Case 1.}  If $A=2$, $B=2.5$, Eq.~\eqref{eq:f(beta)=beta} has the solutions $\beta = 1$ and $\beta = \beta_{\times}(A,B) < 1$, numerically obtained. From Eq.~\eqref{eq:x(beta)}, we obtain respectively the fixed points $\pm\bar{x}_{1} = \pm\left(1,1\right)$ and $\pm\bar{x}_{\beta_{\times}}=\pm\left(0.78,0.63\right)$.
The eigenvalues of the linearized system around $\pm \bar{x}_{1}$ are both real and negative, implying that the points are stable nodes.
The fixed points $\pm\bar{x}_{\beta_{\times}}$ turn out to be saddle points. The phase portrait numerically obtained for this first case is shown in Fig.~\ref{fig:PhDiaNoNoiseNoInternalInteractions:first}.

\noindent \textbf{Case 2.}  If $A=2$, $B=4$, Eq.~\eqref{eq:f(beta)=beta} has the unique solution $\beta = 1$, so the only fixed points, apart from $\left(0,0\right)$, are $\pm \bar{x}_1$.
We refer the reader to the analysis above, which holds for any $A>1$, $B=A+2$, and to Fig.~\ref{fig:PhDiaNoNoiseNoInternalInteractions:second}.

\noindent \textbf{Case 3.}  If $A=2$, $B=7$, Eq.~\eqref{eq:f(beta)=beta} has two solutions: $\beta=1$ and $\beta=\beta_{\times}(A,B) >1$. 
The fixed points $\pm \bar x_{1}$ can be easily seen to be saddle points, while the fixed points $\pm \bar x_{\beta_{\times}}=\pm\left(1.24, 1.58\right)$ have complex conjugate eigenvalues with negative real part, thus they are stable spirals. 
The phase portrait in this case is shown in Fig.~\ref{fig:PhDiaNoNoiseNoInternalInteractions:third}.

\section{Proof of Theorem~\ref{thm:Gaussian_approx}}\label{appendix:D}

The proof of  Theorem~\ref{thm:Gaussian_approx} follows the strategy used in \cite{CoDaPFo15}, where an analogous result for a one population system of mean field interacting particles with dissipation is given.

Recall that $((x(t),y(t)), 0 \leq t \leq T)$ is the unique solution to system \eqref{prop:chaos}.
If $Z\sim N(\mu,v)$ is a Gaussian random variable with mean $\mu$ and variance $v$, we have that $\mathbb{E}(Z^3)= \mu^3 + 3 \mu v$ and $\mathbb{E}(Z^4)= \mu^4 + 6 \mu^2 v+3 v^2$. Therefore, plugging these identities into Eq.~\eqref{eq:moment}, we get differential equations for the mean and the variance of two processes having the same first and second moments as $x(t)$ and $y(t)$, but with Gaussian-like higher-order moments. This is how we obtained system \eqref{eq:odes}. Therefore, part~1 of Theorem~\ref{thm:Gaussian_approx} is true by construction. 

We turn to the second part of the statement. Notice that system \eqref{eq:odes} has a unique global solution, since the four-dimensional vector field is continuous in each variable and has continuous partial derivatives at each point. 
\sloppy Let \mbox{$((m_1(t),m_2(t),v_1(t),v_2(t)), t \geq 0)$} be the unique solution to \eqref{eq:odes}, with initial conditions $m_1(0)=x(0)$, $m_2(0)=y(0)$, $v_1(0)=v_2(0)=0$, and set $V_i(t) \coloneqq \sigma^{-2} \, v_i(t)$ ($i=1,2$). 
The first step of the proof is to define two centered Gaussian processes, $(\xi_1(t),0 \leq t \leq T)$ and \mbox{$(\xi_2(t), 0 \leq t \leq T)$}, so that \mbox{$\mathbb{E}[\xi^2_i(t)]=V_i(t)$} for all $t \in [0,T]$ ($i=1,2$).
We consider the process $(\xi_1(t), 0 \leq t \leq T)$ first. If we write its differential as that of a generic It{\^o}'s process, i.e. $d\xi_1(t) = \psi(t)dt + \phi(t) \, dw_1(t)$, with $\phi$, $\psi$ suitable functions and $(w_1(t), 0 \leq t \leq T)$ a standard Brownian motion, by It{\^o}'s formula we get
$$ 
d\xi^2_1(t) = \left( 2\xi_1(t)\psi(t) +\phi^2(t)\right)dt + 2\xi_1(t)\phi(t) \, dw_1(t).
$$
In turn, we obtain
$$\frac{d\mathbb{E}[\xi^2_1(t)]}{dt} = 2\mathbb{E}[\xi_1(t)\psi(t)]+\mathbb{E}[\phi^2(t)]$$
and we can impose $\phi(t)=1$ and $\xi_1(t)\psi(t) =\xi^2_1(t) \tilde \psi(t)$, where $\tilde\psi(t)$ is a deterministic factor such that $\mathbb{E}[\xi^2_1(t)]$ satisfies the equation for $V_1(t)$, obtained from Eq.~\eqref{eq:odes}. 
Namely, we must require that 
$\tilde\psi(t) = -3 \sigma^2 V_1(t) -3 \left(m_1(t)\right)^2+1-\alpha\theta_{11}-A$. 
With straightforward modifications, we also obtain a differential characterization for the process $(\xi_2(t), 0 \leq t \leq T)$. Putting everything together, we get the following system of SDEs
\begin{align}\label{eq:zs}
     d\xi_1(t) &= \left(-3 \sigma^2 V_1(t) -3 m_1^2(t) +1-\alpha\theta_{11}-A\right)\xi_1(t) \, dt + dw_1(t)\nonumber  \\
     d\xi_2(t) &= \left(-3 \sigma^2 V_2(t)  -3 m_2^2(t) +1+B-\left(1-\alpha\right)\theta_{22}\right)\xi_2(t) \, dt + dw_2(t)\nonumber \\
     \xi_1(0) &= \xi_2(0) = 0.
\end{align}
The processes $(\xi_i(t), 0 \leq t \leq T)_{i=1,2}$ are both Gaussian Markov processes,
with zero mean and such that 
$\mbox{Var}\left[\xi_i(t)\right] = V_i(t)$ for all $t \in [0,T]$ ($i=1,2$). Moreover, they are well-defined, since we have uniqueness of the solution for system \eqref{eq:odes}.\\
Now we define two new processes: 
\begin{equation}\label{eq:approximating_Gaussian_process}
\tilde{x}(t) \coloneqq m_1(t) + \sigma \, \xi_1(t)
\quad \text{ and } \quad
\tilde{y}(t) \coloneqq m_2(t) +\sigma \, \xi_2(t),
\end{equation}
which can be easily seen to be Markovian and Gaussian. Moreover, their respective means $m_1(t)$, $m_2(t)$ and their respective variances $v_1(t)=\sigma^2 \, V_1(t)$, $v_2(t)=\sigma^2 \, V_2(t)$ satisfy Eq.~\eqref{eq:odes} by construction. As a consequence, the processes $(\tilde{x}(t), 0 \leq t \leq T)$ and $(\tilde{y}(t), 0 \leq t \leq T)$ have first and second moments satisfying Eq.~\eqref{eq:moment}.\\
To conclude the proof we need to upper bound the right-hand side of the following inequality
\begin{multline}\label{eq:Gauss_approx_3_rewritten}
        \mathbb{E}\left[\sup_{t\in[0,T]}\left\{\left| x(t) -\tilde x(t)\right|+\left|  y(t) - \tilde y(t)\right|\right\}\right] \\
        \leq \mathbb{E}\left[\sup_{t\in[0,T]}\left|  x(t) -\tilde x(t)\right|\right] + \mathbb{E}\left[\sup_{t\in[0,T]}\left| y(t) -\tilde y(t)\right|\right].
\end{multline}
By using together Eq.~\eqref{eq:zs}, Eq.~\eqref{eq:approximating_Gaussian_process} and Eq.~\eqref{eq:odes}, we obtain
\begin{align*}
d\tilde x(t) &= d m_1(t) +\sigma \, d \xi_1(t)\\[.2cm]
&= \left[-m_1^3(t) - 3 \sigma^2 m_1(t) V_1(t) + m_1(t) - A \left(m_1(t) - m_2(t)\right) \right]dt\\
&+ \sigma \xi_1(t) \left[-3 \sigma^2 V_1(t)  - 3 m_1^2(t) + 1 -\alpha \theta_{11} -A\right]dt + \sigma \, dw_1(t) \\[.2cm]
&= \left[-\tilde{x}^3(t) +\sigma^3 \xi^3_1(t) +3 \sigma^2 m_1(t)  \xi^2_1(t) + \left(1-3 \sigma^2 V_1(t) -A\right)\tilde{x}(t)\right] dt\\
&+ \left[A m_2(t) -\alpha\theta_{11}\sigma \xi_1(t)\right]dt + \sigma \,  dw_1(t)
\end{align*}
and, analogously,
\begin{align*}
d\tilde y(t) &= \left[ -\tilde{y}^3(t) +\sigma^3 \xi^3_2(t) + 3 \sigma^2 m_2(t)  \xi^2_2(t) + \left(1+B -3 \sigma^2 V_2(t) \right) \tilde{y}(t)\right]dt\\
&+ \left[-B m_1(t) -\left(1-\alpha\right)\theta_{22}\sigma \xi_2(t) \right] dt + \sigma \, dw_2(t).
\end{align*}
At this point we follow the very same steps we used before in Appendices~\ref{appendixA} and \ref{appendix:B}. We have 
\begin{align}\label{eq:Gauss_deltax}
x(t) -\tilde x(t) &= \int_{0}^{t}\left(x(s)-\tilde x(s)\right)\left[1-f_1(s)-\alpha \theta_{11}-A\right]ds \nonumber\\
&- \sigma^2\int_0^t \left(\sigma \xi^3_1(s) +3 m_1(s) \xi^2_1(s) -3 m_1(s) V_1(s) -3\sigma V_1(s) \xi_1(s) \right)ds,
\end{align}
with $f_1(s) = x^2(s)+\tilde{x}^2(s)+x(s)\tilde x(s)$. Eq.~\eqref{eq:Gauss_deltax} is of the form $\varphi(t) = \int_0^t\varphi(s) H(s) ds +\int_0^t Q(s) ds$, with $\varphi(t)=x(t)-\tilde{x}(t)$. As $\varphi(0)=0$, the solution to Eq.~\eqref{eq:Gauss_deltax} is given by $\varphi(t) = \int_0^t Q(s) e^{\int_s^t H(r)dr} ds$, where
$$H(s) = 1-f_1(s)-\alpha\theta_{11}-A$$
and 
$$Q(s) = -\sigma^{2}\left(\sigma \xi^3_1(s) +3 m_1(s) \xi^2_1(s) -3 m_1(s) V_1(s) -3 \sigma V_1(s) \xi_1(s) \right).$$
Hence, we have
\begin{equation*}
\left|x(t)-\tilde x(t) \right| \leq \int_0^t \left|Q(s)\right| e^{\int_s^t \left(1-f_1(r)-\alpha\theta_{11}-A\right) dr} ds
\end{equation*}
and, therefore,
\begin{align}\label{eq:Gauss_x}
\mathbb{E}\left[\sup_{t\in[0,T]}\left|x(t)-\tilde x(t) \right|\right]  
&\leq  \mathbb{E}\left[\sup_{t\in[0,T]}\int_0^t \left|Q(s)\right| e^{\int_s^t \left(1-f_1(r)-\alpha\theta_{11}-A\right) dr} ds\right] \nonumber\\[.2cm]
&\leq  \mathbb{E}\left[\int_0^T|Q(s)| \sup_{t\in[0,T]}e^{\int_{s}^t\left(1-f_1(r)-\alpha\theta_{11}-A\right)dr} ds \right]\nonumber\\[.2cm]
&\leq  C_T\int_0^T\mathbb{E}\left[\left|Q(s)\right|\right]ds\nonumber \\
&\leq  \tilde C_T \sigma^2,
\end{align}
for some $\tilde{C}_T>0$. The last inequality follows from the fact that we have introduced $\tilde{Q}(s)=\sigma^{-2}Q(s)$, that, being 
a polynomial function of a Gauss-Markov process, has a time-locally bounded $L^1$-norm.\\
An analogous estimate holds for the second term in the right-hand side of \eqref{eq:Gauss_approx_3_rewritten}, so that
\begin{equation}\label{eq:Gauss_y}
\mathbb{E}\left[\sup_{t \in [0,T]} \left| y(t)- \tilde y(t)\right| \right]  \leq \tilde D_T \sigma^2,
\end{equation}
for a suitable positive constant $\tilde{D}_T$. Putting together Eq.~\eqref{eq:Gauss_approx_3_rewritten}, Eq.~\eqref{eq:Gauss_x} and Eq.~\eqref{eq:Gauss_y} yields the conclusion of the proof of part~2 in Theorem~\ref{thm:Gaussian_approx}.
%

\section{Subcritical Hopf bifurcation}\label{Hopf}

We consider the dynamical system \eqref{eq:odes} with $\theta_{11}=\theta_{22}=8$ and $A$ and  $B$ chosen in the regime given by assumption \textbf{(H)}. We study the nature of the equilibrium point $(m_1,m_2,v_1,v_2)=\left(0,0,\tilde v_1,\tilde v_2\right)$, where
\begin{equation*}
        \begin{split}
            & \tilde v_1 = \frac{1}{6} \left(-3-A +\sqrt{\left(-3-A\right)^2 + 6 \sigma^2} \right)\\
            & \tilde v_2 = \frac{1}{6}\left(-3+B + \sqrt{\left(-3+B\right)^2 + 6\sigma^2}\right),
        \end{split}
\end{equation*}
as the noise intensity $\sigma>0$ varies. The following analysis reveals the presence of a subcritical Hopf bifurcation at the point  $(m_1,m_2,v_1,v_2)=\left(0,0,\tilde v_1,\tilde v_2\right)$, in all the three parameter regimes examined in the main text (see Table~\ref{tab:comparison_fft_poincare} and Fig.~\ref{fig:PhDiaNoNoiseNoInternalInteractions}). Recall that a Hopf bifurcation occurs when a stable periodic orbit arises from an equilibrium point as, at some critical value of the parameter, it loses stability. Subcritical means that - as in the present case - such a transition happens when moving the parameter from larger to smaller values. A Hopf bifurcation can be detected by checking whether a pair of complex eigenvalues of the linearized system around the equilibrium crosses the imaginary axis as the parameter changes. See Theorem~2, Chapter~4.4 in \cite{Per2001}. We briefly analyze our case.

The Jacobian matrix relative to system \eqref{eq:odes} at $\left(0,0,\tilde v_1, \tilde v_2\right)$ reads 
\begin{equation*}\label{eq:J}
    J (\sigma) = 
        \begin{bmatrix}
            1-A-3\tilde v_1 & A & 0 & 0 \\
            -B & 1+B-3\tilde v_2 & 0 & 0 \\
            0 & 0 & -6 -2 A-12\tilde v_1 & 0  \\
            0 & 0 & 0 & -6 + 2 B -12 \tilde v_2
        \end{bmatrix}
    \end{equation*}
and its eigenvalues are
\begin{allowdisplaybreaks}
\begin{align*}
            &\lambda_1 = \frac{1}{4}\left(10-A+B-\sqrt{\left(A+3\right)^2 + 6\sigma^2}-\sqrt{\left(B-3\right)^2+ 6\sigma^2} \right.\\
            &\quad \left. - \sqrt{2}\left[A^2+A\left(\sqrt{\left(A+3\right)^2 + 6\sigma^2}-\sqrt{\left(B-3\right)^2+6\sigma^2}-7 B + 3\right)\right.\right.\\
            & \left.\left.\quad-\left(\sqrt{\left(B-3\right)^2+6\sigma^2}-B\right)\left(\sqrt{\left(A+3\right)^2 + 6\sigma^2}+B\right) - 3 B + 6\sigma^2 +9\right]^{\frac{1}{2}}
            \right)\\[.2cm]
            & \lambda_2 = \frac{1}{4}\left(10-A+B-\sqrt{\left(A+3\right)^2 + 6\sigma^2}-\sqrt{\left(B-3\right)^2+ 6\sigma^2} \right.\\
            &\quad \left. + \sqrt{2}\left[A^2+A\left(\sqrt{\left(A+3\right)^2 + 6\sigma^2}-\sqrt{\left(B-3\right)^2+6\sigma^2}-7 B + 3\right)\right.\right.\\
            &\left.\left. \quad-\left(\sqrt{\left(B-3\right)^2+6\sigma^2}-B\right)\left(\sqrt{\left(A+3\right)^2 + 6\sigma^2}+B\right) - 3 B + 6\sigma^2 +9\right]^{\frac{1}{2}}
            \right)\\[.2cm]
            & \lambda_3 = -6 -2 A-12\tilde v_1\\[.2cm]
            & \lambda_4 =-6 + 2 B -12 \tilde v_2.\\
\end{align*}
\end{allowdisplaybreaks}

\begin{figure}[h!]
        \centering
        \subfigure[]{\includegraphics[scale=0.18]{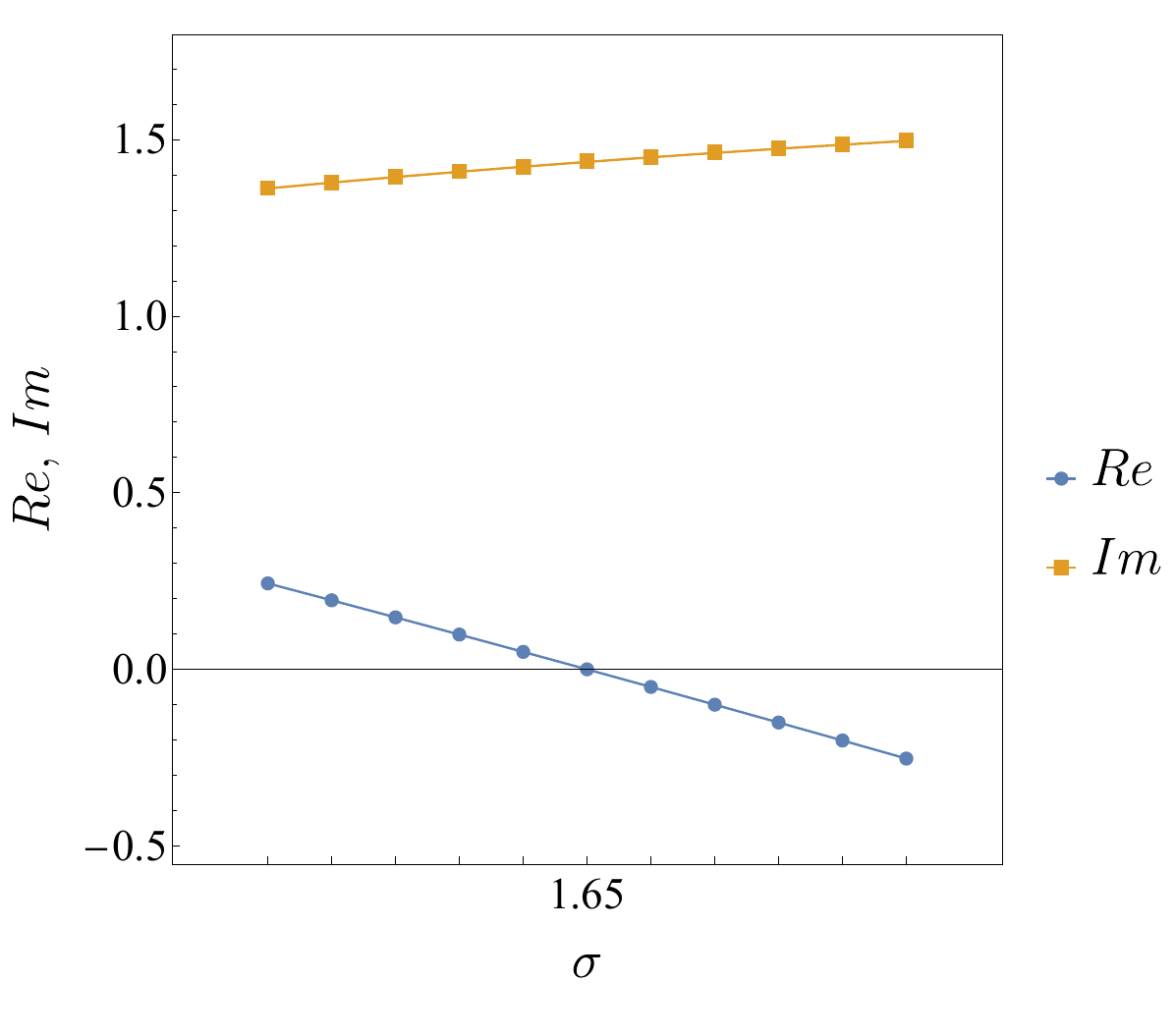}}
        \subfigure[]{\includegraphics[scale=0.18]{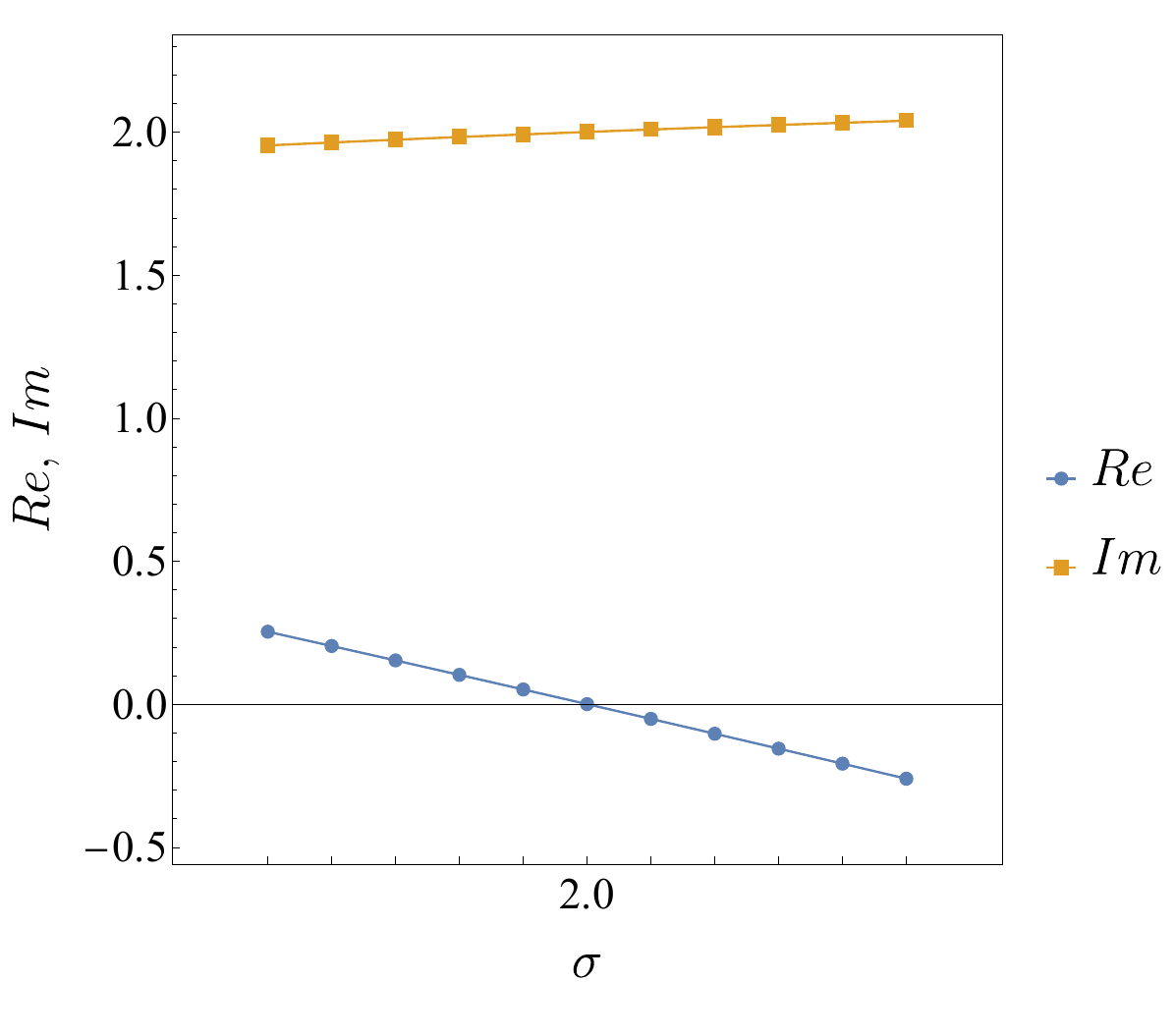}}
        \subfigure[]{\includegraphics[scale=0.18]{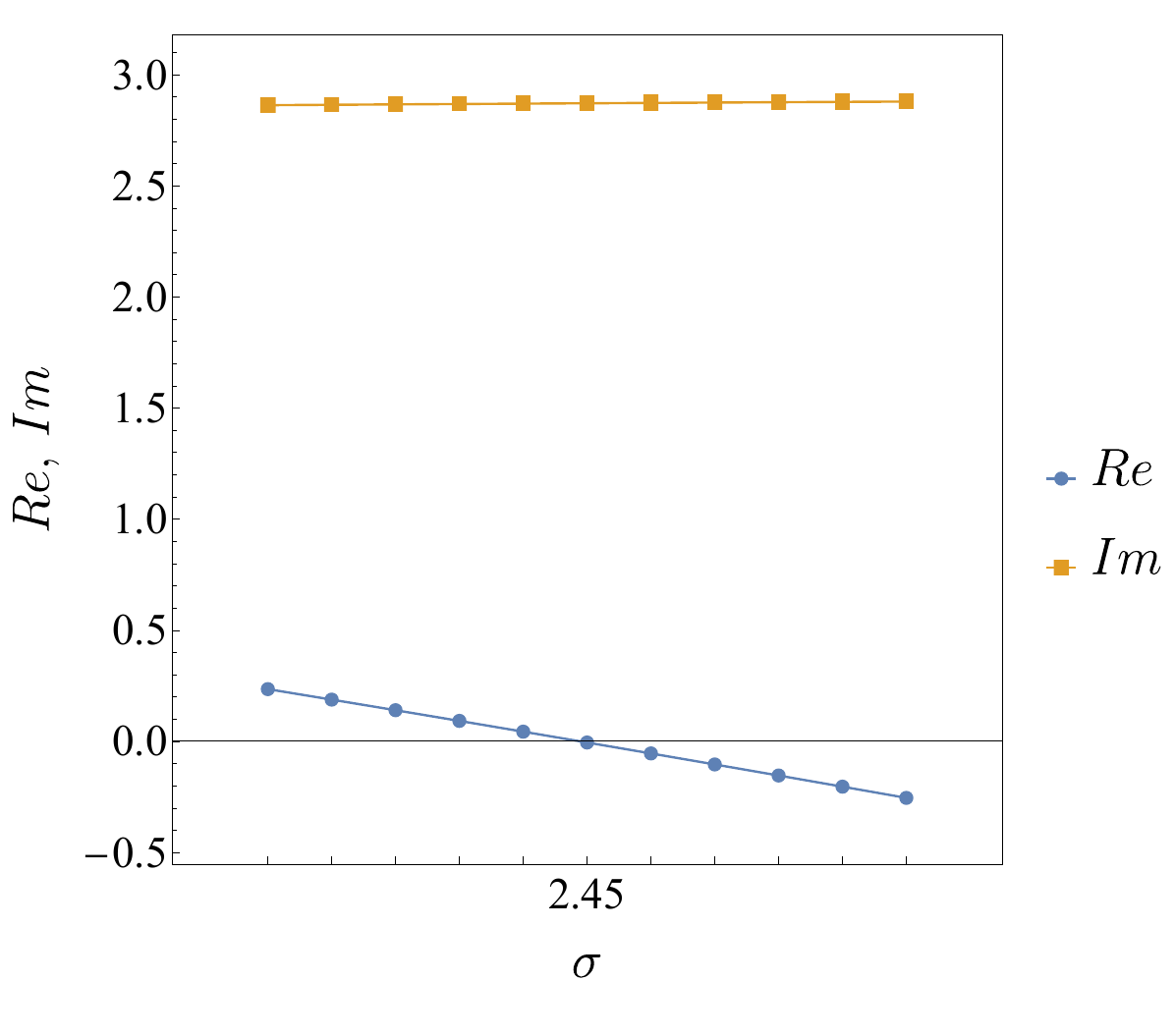}}
        \caption{\footnotesize{Real parts (blue) and absolute values of the  imaginary parts (orange) of the eigenvalues $\lambda_1$, $\lambda_2$ of the Jacobian matrix $J(\sigma)$ as functions of $\sigma$ with: (a) $A=2$, $B=2.5$, $\sigma$ ranging from 1.4 to 1.9 with step 0.05; (b) $A=2$, $B=4$, $\sigma$ ranging from 1.75 to 2.25 with step 0.05; (c) $A=2$, $B=7$, $\sigma$ ranging from 2.2 to 2.7 with step 0.05.}}
        \label{fig:eigs}
    \end{figure}

The eigenvalues $\lambda_3$ and $\lambda_4$ are negative for all $\sigma>0$. The eigenvalues $\lambda_1$ and $\lambda_2$ are complex conjugate when (a) $A=2$ and $B=2.5$, (b) $A=2$ and $B=4$, (c) $A=2$ and $B=7$, and  we checked numerically they cross the imaginary axis with negative derivative at the respective critical values (a) $\sigma_c\simeq 1.65$, (b) $\sigma_c\simeq 2$, (c) $\sigma_c\simeq 2.45$ (see Fig.~\ref{fig:eigs}). These evidences confirm the presence of a subcritical Hopf bifurcation for the regimes examined in the main text.   
 
\end{appendices}






\bibliographystyle{abbrv}

\vspace{25pt}
\textit{E-mail addresses: } elisa.marini@math.unipd.it, luisa.andreis@polimi.it, francesca.collet@univr.it, marco.formentin@unipd.it

\end{document}